\documentclass[10pt]{article}

\usepackage[dvips]{graphicx}
\usepackage{makeidx}
\usepackage{amsmath}
\usepackage[russian,english]{babel}

\newcommand {\etavec}{\boldsymbol{\eta}}
\newcommand {\omegavec}{{\boldsymbol{\omega}}}

\newcommand {\alphavec}{\boldsymbol{\alpha}}
\newcommand {\xivec}{{\boldsymbol{\xi}}}
\newcommand {\varthetavec}{{\boldsymbol{\vartheta}}}
\newcommand {\nuvec}{\boldsymbol{\nu}}
\newcommand {\psivec}{\mbox{\boldmath $\psi$}}

\newcommand {\thetavec}{{\boldsymbol{\theta}}}

\newfont{\pseudocode}{cmtt10}

\newtheorem{thm}{Theorem}
\newtheorem{cor}{Corollary}
\newtheorem{lem}{Lemma}
\newtheorem{prop}{Proposition}

\newtheorem{assume}{Assumption}

\newtheorem{rem}{Remark}
\newcommand{\bfx}{\mathbf{x}}
\newcommand{\bff}{\mathbf{f}}
\newcommand{\bfz}{\mathbf{z}}
\newcommand{\bfy}{\mathbf{y}}
\newcommand{\bfh}{\mathbf{h}}

\newcommand{\bfg}{\mathbf{g}}

\newcommand{\bfq}{\mathbf{q}}

\newcommand{\dpsi}{{\dot\psi}}
\newcommand{\pd}{{\partial}}

\makeindex
\begin{document}
% for Improved Performance!

%\twocolumn

\baselineskip 0.7cm

\title{\bf Adaptive Algorithms in Finite Form}
%\title{\bf Adaptive Algorithms in a Finite Form for Nonlinear Dynamic Plants}
\author{Ivan Tyukin\thanks{Laboratory for Perceptual Dynamics, RIKEN (Institute for Physical and Chemical Research)
                           Brain Science Institute, 2-1, Hirosawa, Wako-shi, Saitama, 351-0198, Japan, e-mail: \{tyukinivan,ceesvl\}@brain.riken.jp}
%\thanks{The authors are grateful to Prof. V.A.Terekhov for his useful comments.},
Danil Prokhorov\thanks{Ford Research Laboratory, Dearborn, MI,
48121, USA, e-mail: dprokhor@ford.com} and Cees van
        Leeuwen$^\ast$}
\date{\today}
\date{}
\maketitle{}
%\thanks{Laboratory for Perceptual Dynamics, RIKEN (Institute for Physical and Chemical Research)
%                           Brain Science Institute, 2-1, Hirosawa, Wako-shi, Saitama, 351-0198, Japan, e-mail: \{tyukinivan,ceesvl\}@brain.riken.jp}
%\thanks{Ford Research Laboratory, Dearborn, MI, 48121, USA, e-mail:
%dprokhor@ford.com}

\begin{abstract}

We propose a new method for the design of adaptation algorithms
that guarantees a certain prescribed level of performance and
applicable to systems with nonconvex parameterization. The main
idea behind the method is two-fold. First, we augment the tuning
error function and design the adaptation scheme in the form of
ordinary differential equations. The resulting augmentation is
allowed to depend on state derivatives.  Second, we find a
suitable realization of the designed adaptation scheme in an
algebraic-integral form.  Due to their explicit dependence on the state of the original system,
such adaptation schemes are referred to as adaptive algorithms in {\it
finite form}, in contrast to (conventional) algorithms in differential form.
Sufficient conditions for the existence of finite
form realizations are proposed. It is shown that our method to
design algorithms in finite form is applicable to a broad class of
nonlinear systems including systems with nonconvex
parameterization  and low-triangular systems.
%These conditions necessitate
%finding a solution of a system of partial differential equations,
%which in general is not an easy task. Several full-state feedback
%finite form realizations are presented and illustrated with
%examples.
\end{abstract}

%\vspace{10mm}


% \vspace{10mm}
%{ \small
%{\bf Corresponding author:}{  \\Ivan Tyukin \\
%                              Laboratory for Perceptual Dynamics,\\
%                              RIKEN Brain Science Institute,\\
%                              2-1, Hirosawa, Wako-shi, Saitama,\\
%                              351-0198, Japan\\
%                              phone: +81-48-462-1111 extension 7436\\
%                              fax:   +81-48-467-7236\\
%                              e-mail: tyukinivan@brain.riken.jp}
%}
\section{Introduction}

Significant progress in adaptive control theory has been made in
the areas of linear and nonlinear systems
\cite{Jac,Kokotovich95,Narendra89,Sastry89}, plants with relative
degree greater than one \cite{Kreisselmeier77, Morse92,
MarinoTomei93}, and systems with nonconvex parameterization
\cite{Annaswamy99,Lin}. However, there is still room for further
developments, as there are important unresolved problems regarding
the issue of performance, especially in the presence of nonconvex
parameterization.

As expressed by asymptotic stability of
adaptive systems \cite{Sastry89}, robustness, and good transient
behavior \cite{Narendra_1988_exc},
suitable performance can be proven under the
requirement of persistent excitation.
%, such guarantees cannot be given for the more general
%case.
As it is generally observed in practice, insufficient excitation
results in absence of asymptotic stability and, as its consequence,
in poor parameter convergence, sensitivity to small disturbances
%\cite{Tsakalisa_1996}
and poor transient performance. A suitable performance criterion
is needed  to assure the efficiency and quality of the system.
%These restrictions share, in our understanding, the same
%underlying problem: the lack of adequate performance criteria.
As substitutes for performance criteria, most of the available
results in direct adaptive control
%under normal conditions
without restrictive persistent excitation requirements limit
themselves to $L_2$ and $L_\infty$\footnote{Function $\nu:
R_+\rightarrow R$ is said to belong to $L_2$ iff
$L_2(\nu)=\int_0^{\infty}\nu^2(\tau)d\tau<\infty$. The value
$\sqrt{L_2(\nu)}$ stands for the $L_2$ norm of $\nu(t)$. Function
$\nu: R_{+}\rightarrow R$ belongs to and $L_\infty$ iff
$L_{\infty}(\nu)=\sup_{t\geq0}\|\nu(t)\|<\infty$, where
$\|\cdot\|$ is the Euclidean norm. The value of $L_{\infty}(\nu)$
stands for the $L_{\infty}$ norm of $\nu(t)$.} norm bounds of the
tracking errors. For more sophisticated performance measures like
the LQ criterion, some results are available
\cite{French2000,French2002}. These results, however, deal either
with too narrow a class of uncertain systems \cite{French2000} or
present only a comparison
%suggest only a qualitative comparison
between adaptive and robust backstepping \cite{French2002} without
suggesting new adaptation schemes.
% Few attempts have been made to
%show by  rigorous analysis how to enable an adaptation scheme to
%meet a certain pre-specified performance criterion.
On the other
hand, when improvement {\it heuristics} are suggested like in
\cite{Narendra94}, no exact performance criterion is provided that
can explicitly be computed {\it a-priori}, except probably the
bounds on $L_2$ and $L_{\infty}$ norms for the tracking errors.

Another unresolved issue in conventional adaptive control theory
is nonconvex parameterization of the plant model. The available
approaches encourage
%the designer
to compensate for the
nonlinearity
%(at least, in part)
by using an additional damping
term
%, or high-gain feedback
\cite{Annaswamy99,Lin,Annaswamy2000}.
These techniques guarantee existence of a solution to the control
problem, but are limited in terms of their practical value because
of their high-gain nature.
%Very recent results on nonparametric
%adaptation
% \cite{Ortega02} can also be applied to nonlinearly parameterized
%systems.
%Nevertheless,
Moreover, none of these approaches provides a performance measure
beyond $L_\infty$ or $L_2$ norm bounds with respect to the state
or tracking error.

% merely can provide
%integrability of the squared error.

An important impediment to solving these two major problems at
once, we believe, is the lack of sufficient information in the
conventional adaptive schemes.
% Our aim, therefore, is to clean the
%way for a strategy to provide the information needed.
% to
%improve the performance and deal with nonconvex parameterization.
One way to provide the algorithms with extra information is to
augment the tuning errors.
% Many adaptive control schemes use error
%augmentation to make the estimation error dependent on controller
%parameters or computable derivatives of the known signals.
This idea is inherent to both Morse's adaptive controllers
\cite{Morse92} and those based on Kreisselmeier's observers
\cite{Kreisselmeier77} when dealing with  plants with relative
degree greater than one. These augmented errors then are used in
conventional gradient schemes
\begin{equation}\label{alg0}
\dot{\hat{\thetavec}}=-\Gamma
{\psi}(\bfx,t)\mathcal{A}(\bfx,\hat{\thetavec},t),
\end{equation}
where $\bfx\in R^n$ is a state (or output) vector,
$\hat{\thetavec}\in R^d$ is a vector of controller parameters,
$\mathcal{A}(\bfx,\hat{\thetavec})$ is an operator that depends on
the particular problem at hand, ${\psi}(\bfx,t)$ is the error
function and gain $\Gamma>0$. Improved performance of these
augmented controllers for Morse's high-order tuners is reported in
\cite{Datta_1995,Ortega_1993}. In these papers the plants are
assumed to be linear, but no additional criteria are
provided except for $L_2$ and $L_\infty$ norm bounds with respect to the tracking errors.
%It is shown in \cite{Datta_1995}, however,  that these
%bounds can be attenuated by the controller parameters and
%adaptation gain if we neglect the influence of non-zero initial
%conditions. Similar results apply to adaptive backstepping
%controllers \cite{Krstic_1993}.  The price for such an improvement
%is that the adaptation gains should be sufficiently large.

The above-mentioned limitations of performance and applicability
motivate us to search for a new augmentation that uses additional
information about the system dynamics, such as, e.g., state
derivatives.  As a result, new properties in the
system can be created.
% (in addition to readily achievable finiteness of $L_2$ and
%$L_\infty$ norm bounds for the tracking errors).
On the other hand, we wish to find physically realizable
algorithms %,
%with such
%an augmentation,
%i.e. computational procedures
that do not require measurements of
any unknown signals, derivatives, or parameters. In order to meet
these seemingly contradictory requirements we propose to
%drop the usual restrictions of realizability of
%adaptation algorithms in differential form, ----------------------
%when
%would like to search for such algorithms which if being written in
%differential form (\ref{alg0}) may result in the desired
%augmentation.
%
%with controller parameters $\hat{\thetavec}$ that can be computed
%at any time instant without measuring the unknown signals or
%parameters.
%
%The question of existence of the adaptive control algorithms
%guaranteeing improved performance and applicability to nonconvexly
%parameterized models is not trivial in general.
%We propose a structure  first and then search for the algorithms
%with improved properties within an already proposed class of
%algorithms.
%Instead of searching for the desired adaptive
%algorithm in the conventional form defined by (\ref{alg0}),
%We suggest  to extend algorithms (\ref{alg0}) as follows:
%when
extend conventional algorithms (\ref{alg0}) as follows:
\begin{eqnarray}\label{fin_form_ours}
%\hat{\thetavec}(\bfx(t),t)=-\Gamma\left(\mathcal{A}_1(\bfx,t)+\int_{0}^{t}\mathcal{A}_2(\bfx(\tau),\hat{\thetavec}(\tau)
\hat{\thetavec}(\bfx(t),t)&=&\hat{\thetavec}_P(\bfx,t)+\hat{\thetavec}_I(t);
\ \dot{\hat{\thetavec}}_I=\mathcal{A}_2(\bfx,\hat{\thetavec},t), \
\ \hat{\thetavec}_P(\bfx,t)=\mathcal{A}_1(\bfx,t).
\end{eqnarray}
%It is obvious that algorithms (\ref{fin_form_ours}) contain
%algorithms (\ref{alg0}) as  a special case.
Notice that if functions $\hat{\thetavec}(\bfx,t)$ are written in
differential form (\ref{alg0}), they may depend on unknown
parameters and unmeasured signals, e.g, state derivatives. Thus,
the equivalent description of adaptive algorithm
(\ref{fin_form_ours}) in differential form
%(\ref{alg0})
may
produce an augmentation that is in fact derivative-dependent,
thereby providing the
%adaptation
algorithm with more information about the plant uncertainties.
These observations lead to quite unexpected consequences. Instead
of restricting the design procedure to the algorithms given in
differential form (\ref{alg0}),
%to those algorithms that can
%be realized in the form of equation (\ref{alg0}),
it becomes possible to design adaptation algorithms in two steps.
First, search for the desired augmentation  to obtain the required
 adaptive control properties. At this stage it does not matter whether the augmentation is
uncertainty-dependent or not.
%\footnote{It was
%shown, for instance, in \cite{tpt2002_at, tpt2003_tac} that the
%derivative-dependent algorithms written in differential form
%(\ref{alg0}) are able to deal with a class nonconvexly
%parameterized plant.}.
Once a suitable tuning error is chosen, second: find a realization
of the algorithm in the form of integral-algebraic equation
(\ref{fin_form_ours}). Such a realization will be termed {\it
algorithm in finite form}.
%\footnote{To the best of our
%knowledge, such terminology in adaptive control literature has
%been proposed first by A. Fradkov in \cite{Fradkov90} in order to
%denote the speed-gradient based adaptive control schemes given by
%integral-algebraic, not differential equations.}.

According to our knowledge, algorithms (\ref{fin_form_ours}) have
been introduced for direct adaptive control of nonlinear systems
in 1986 in \cite{Fradkov86} and then were reintroduced later in
\cite{Ortega02,Ortega2003}.
%, where they were
%termed algorithms in finite form.
Their distinctive performance properties  and extended
applicability, however, were not appreciated at the time.
% reported
%in \cite{Fradkov86,Fradkov90}.
The efficiency of the proposed algorithms in
\cite{Fradkov86,Fradkov90} was limited by restrictive
pseudo-gradient assumptions on $\hat{\thetavec}_P(\bfx,t)$ (for
the details see \cite{Fradkov86,Fradkov90}).
%, although the control goal reaching  was shown
%to be guaranteed for those plants that satisfy the convexity
%restriction \cite{Jac}.
Nevertheless, it has been reported recently that algorithms
(\ref{fin_form_ours}) may be able to deal with nonconvex
parameterization (see for example \cite{tpt2003_tac}, Lemma 1, p.
558; \cite{Ortega02,t_fin_formsA&T2003}) and guarantee improved
transient performance \cite{Stotsky93,t_fin_formsA&T2003}.
% -- (it
%has been shown that for the class of algorithms
%(\ref{fin_form_ours}) $\dpsi\in L_2$ is also guaranteed).
Some preliminary analysis of the distinctive properties of
algorithms (\ref{fin_form_ours}) is available, for example, in
\cite{t_fin_formsA&T2003,ECC_2003}. Additional support and
motivation for algorithms (\ref{fin_form_ours}) can be found in
\cite{Ortega2003}, where the authors introduced their adaptation
schemes from immersion and invariance principles. However, the
main problem,  with the current study of algorithms
(\ref{fin_form_ours}) is that there is no systematic method that
allows us to design these algorithms with guaranteed improvements
in performance and, at the same time, achieve applicability to
systems with nonlinear parameterization for a sufficiently broad
class of nonlinear dynamical systems.

In our present work we suggest a new method to design adaptive
algorithms in finite form (\ref{fin_form_ours}) that guarantee
improved performance and in addition are applicable to a class of
nonlinearly parameterized plants. The method is systematic and is
based on two fundamental ideas in adaptive control theory: {\it
augmentation} of the error and {\it embedding}
%\footnote{During
%revision and preparation of the manuscript the authors became
%aware of publication \cite{Ortega2003} by A.Astolfi and R.Ortega,
%where algorithms (\ref{fin_form_ours}) were introduced on the
%basis of immersion principle. Though the authors pointed out
%advantages of the algorithms in finite forms (\ref{fin_form_ours})
%over those given by (\ref{alg0}), no formal proofs of the
%transient behavior improvements were given. This fact additionally
%support relevance of our current study. In addition, our method is
%motivated primarily by performance and applicability limitations
%of the existing adaptation schemes, whereas method presented in
%\cite{Ortega2003} aims general stabilization problem.}
the original system dynamics into one of a higher order. These
ideas are embodied in two independent stages of the design. The
first stage is {\it augmentation} of the tuning error for
algorithms in the conventional differential form (\ref{alg0}) in
order to ensure improved performance and extended applicability of
these algorithms. The resulting augmentation may not
necessarily be independent on the uncertainties or
time-derivatives of the state vector. The resulting augmentation,
however, should guarantee certain desired properties of the
%adaptive
system.

Based on the augmentation (possibly, derivative-dependent)
obtained in the first stage,
the second stage of the design method should be to find functions $\mathcal{A}_1(\bfx,t)$,
$\mathcal{A}_2(\bfx,\hat{\thetavec},t)$ which guarantee that
algorithms (\ref{fin_form_ours}) realize the desired adaptation
scheme. We show that this problem may require finding a
solution of a system of partial differential equations.
% The
%results are presented in Part I of the paper.
It is well-known that such a solution may not exist in general. To
avoid this problem, we consider several special cases of plant
models, with their structures satisfying sufficient
conditions for the existence of a solution.  As soon as these
basic structures are found, we {\it embed} the original system
into a system of higher order for which the solution is known to
exist.
%These results are given in Part II of the paper.
The embedding is to be made in such a way that the extended system
belongs to one of already established basic classes that guarantee
existence of the solution to the realization problem. With
embedding we shall be able to obtain adaptation schemes that
guarantee not only square integrability of the error but also
integrability of its first derivatives as  well as square
integrability of control efforts injected into the system due to
the parametric uncertainties. In addition, we provide the
conditions for which the decrease of the parametric uncertainties
and exponential convergence into a neighborhood of the target
manifold are guaranteed without the restrictive assumption of
persistent excitation. Last but not the least, our new adaptive
schemes can also be applied to systems with nonlinear
parameterization.

The layout of the paper is as follows. In Section 2 we specify the
class of nonlinear dynamical systems under consideration and
select the desired augmentation.
%Some results of this section are
%extensions of known ones and are presented here for consistency.
Section 3 contains the main results of the paper.
% We show that,
%for a given augmentation, the realization problem is solvable for
%a broad class of nonlinear systems. We provide sufficient
%conditions which guarantee the existence of solutions.
In Section 4 we present an example of the design and results of
computer simulations.
%Section 5 contains an example of the design and
%results of computer simulations.
Section 5 concludes the paper.

\section{Problem Formulation and Preliminary Results}

Let the plant mathematical model be given as follows:
%\begin{eqnarray}\label{plant_ins2}
%\dot{\bfx}&=&\bff(\bfx) + \varthetavec(\bfx,\thetavec) +
%\bfg(\bfx)u,
%\end{eqnarray}
%or
\begin{eqnarray}\label{plant_ins}
\dot{x}_i&=&f_i(\bfx)+g_i(\bfx)u, \ \ \ \ i=1,\dots,m\nonumber \\
\dot{x}_j&=&f_j(\bfx)+\nu_{j-m}(\bfx,\thetavec)+g_j(\bfx)u, \ \ \
\ j=m+1,\dots,n,
\end{eqnarray}
where $\bfx\in R^n$ is a state vector, $f_i, \ g_i: R^n\rightarrow
R$, $f_i, g_i \in C^1$, $\thetavec\in \Omega_{\theta}\subset
R^{d}$ is a vector of unknown parameters, $\nu_i: R^{n}\times
R^{d}\rightarrow R$, $\nu_i\in C^1$, $u$ is a control input. %Thus,
Let us define functions $\bff(\cdot)$, $\bfg(\cdot)$,
$\nuvec(\cdot,\cdot)$:
\[
\bff(\bfx)=(f_1(\bfx),\dots,f_n(\bfx))^{T}, \ \
\bfg(\bfx)=(g_1(\bfx), \dots,g_n(\bfx))^{T}, \ \
%\varthetavec(\bfx,\thetavec)=
\nuvec(\bfx,\thetavec)=(\nu_1(\bfx,\thetavec),\dots,\nu_{n-m}(\bfx,\thetavec))^{T}
%(\vartheta_1(\bfx,\thetavec),\dots\vartheta_m(\bfx,\thetavec),\vartheta_{m+1}(\bfx,\thetavec),\dots\vartheta_{n}(\bfx,\thetavec))^{T}=
%(0,\dots,0,\nu_{1}(\bfx,\thetavec),\dots,\nu_{n-m}(\bfx,\thetavec))^{T}
\]
%For the sake of compactness when dealing with the partial
%derivatives of a function  we will use the following notation:
%\[
%L_{\bff}\psi(\bfx)=\frac{\pd \psi}{\pd \bfx}\bff(\bfx).
%\]
%to denote the Lie derivative of function $\psi(\bfx)$ along the
%vector field $\bff(\bfx)$. If function
%$\psivec(\bfx)=(\psi_1(\bfx),\dots,\psi_n(\bfx))^{T}$ is a
%vector-function then symbol $L_{\bff}\psi(\bfx)$ denotes the
%following vector:
%\[
%L_{\bff}\psivec(\bfx)=\left(
%L_\bff\psi_1(\bfx),\dots,L_\bff\psi_n(\bfx)
%                         \right)^{T}
%\]

It will be useful sometimes to think of state vector
$\bfx\in\mathcal{L}\subseteq R^n$ as
%\begin{eqnarray}%\label{state_partition}
$\bfx=\bfx_1\oplus\bfx_2$, $\bfx_1=\left(x_1,x_2,\dots,x_m
\right)^T$, $\bfx_2=\left(x_{m+1},x_{m+2},\dots, x_n\right)^T$,
%\end{eqnarray}
where symbol $\oplus$ denotes  concatenation of two vectors
$\bfx_1\in \mathcal{L}_1\subseteq R^m$, $\bfx_2\in
\mathcal{L}_2\subseteq R^{n-m}$ and $\mathcal{L}$,
$\mathcal{L}_1$, $\mathcal{L}_2$ are linear spaces. The
time-derivative of $\bfx_1$ is independent on $\thetavec$, whereas
the time-derivative of vector $\bfx_2$ depends on unknown
parameters $\thetavec$ explicitly. Therefore we refer to the
spaces $\mathcal{L}_1$ and $\mathcal{L}_2$ as {\it
uncertainty-independent} and  {\it uncertainty-dependent
partitions} of system (\ref{plant_ins}), respectively. To denote
the right-hand sides of the partitioned system, we use the
following notations: $\bff_1=\left(f_1,\dots, f_m\right)^T$,
$\bff_2=\left(f_{m+1},\dots,f_n\right)^{T}$,
$\bfg_1=\left(g_1,\dots,g_m\right)^{T}$,
$\bfg_2=\left(g_{m+1},\cdots,g_n\right)^{T}$. Hence, the
partitioned system can be written as follows:
\begin{eqnarray}\label{partitioned_plant}
\dot{\bfx}_1&=& \bff_1(\bfx)+\bfg_1(\bfx)u, \ \  \dot{\bfx}_2=
\bff_2(\bfx)+\nuvec(\bfx,\thetavec)+\bfg_2(\bfx)u.
\end{eqnarray}
In analogy with the definition of independence of a function with
respect to the components $x_i$ of its argument $\bfx$, we would
like to define a notion of independence of the function with
respect to the partition. Let
$\bfx\in\mathcal{L}=\mathcal{L}_1\oplus\mathcal{L}_2$, and let function
$\omegavec(\bfx): R^{n}\rightarrow R^{n}$ be differentiable for
any $\bfx \in R^n$. Function $\omegavec(\bfx)$ is said to be {\it
independent on partition $\mathcal{L}_2$} iff ${\pd
\omegavec(\bfx_1\oplus\bfx_2)}/{\pd \bfx_2}=0$. We would also like
to extend the standard definition of the Lie derivatives to the
partitioned system. Given the following partition
$\bfx=\bfx_1\oplus\bfx_2\oplus\cdots\oplus\bfx_r$, we denote by
symbol $L_{\bff_i}\psi(\bfx)$ the following derivatives ${(\pd
\psi(\bfx)}/{\pd \bfx_i}) \ \bff_i(\bfx)=L_{\bff_i}\psi(\bfx)$,
where $\bff_i(\bfx)$ stands for the corresponding vector-function
in
$\bff(\bfx)=\bff_1(\bfx)\oplus\bff_2(\bfx)\oplus\cdots\oplus\bff_r(\bfx)$.
If function $\psivec(\bfx)=(\psi_1(\bfx),\dots,\psi_n(\bfx))^{T}$
is a vector-function then symbol $L_{\bff_i}\psivec(\bfx)$ denotes
the following vector: $L_{\bff_i}\psivec(\bfx)=\left(
L_{\bff_i}\psi_1(\bfx),\dots,L_{\bff_i}\psi_n(\bfx)
                         \right)^{T}$.

As in \cite{tpt2003_tac,tpt2002_at}, we define the control goal as
reaching asymptotically a target manifold. We assume that the
target manifold can be given by the following equality
$\psi(\bfx,t)=0$, where $\psi: \ R^{n}\times R \rightarrow R$,
$\psi(\bfx,t)\in C^1 $. Additional restrictions on the function
$\psi(\bfx,t)$ are formulated in Assumptions
\ref{boundedness_psi}, \ref{Regularity}.
\begin{assume}[Boundedness of the
Solutions]\label{boundedness_psi} Function $\psi(\bfx,t)$ is such
that for any $\delta>0$ there exists a function $\varepsilon: \
R_+\rightarrow R_+$ such that $|\psi(\bfx,t)|\leq \delta
\Rightarrow \|\bfx\|\leq \varepsilon(\delta)$ along system
(\ref{plant_ins}) solutions.
\end{assume}
Assumption \ref{boundedness_psi} simply states that any trajectory
of system (\ref{plant_ins}) belonging to a neighborhood of the
target manifold $\psi(\bfx,t)$ is bounded. Clearly, most of the
common goal criteria used in adaptive control satisfy this
property, for example, positive-definite functions $\psi(\bfx)$
for nonlinear systems and quadratic forms for linear ones.
%This property of the goal function is used in
%classical schemes to show that boundedness of a state follows from
%the boundedness of $\psi(\bfx,t)$. Once the boundedness of $\bfx$
%is established it is possible to show that $\psi(\bfx,t)\in L_2
%\Rightarrow \psi(\bfx,t)\rightarrow 0$ as $t\rightarrow \infty$
%under the assumptions that the system's right-hand side is locally
%bounded\footnote{Function $\bff(\bfx): R^{n}\rightarrow R^n$ is
%said to be locally bounded if for any $\|\bfx\|<\delta$ there
%exists such constant $D(\delta)>0$ that the following holds:
%$\|\bff(\bfx)\|\leq D(\delta)$.} and $\pd\psi(\bfx,t)/\pd t$ is
%uniformly bounded in $t$.
In general, however, in order  to show the boundedness of $\bfx$
it is not necessary for the function $\psi(\bfx)$ to be positive
definite. As an illustration, consider the following system:
\begin{eqnarray}\label{plnat_ins_ex}
\dot{x}_i&=&f_i(x_1,\dots,x_{i})+x_{i+1}, \ i=\{1,\dots,n-1\}\nonumber\\
\dot{x}_n&=&f_n(x_1,\dots,x_{n})+\nu(\bfx,\thetavec)+u.
\end{eqnarray}
Let
$\psi(\bfx)=x_n-p(x_1,\dots,x_{n-1})+f_{n-1}(x_1,\dots,x_{n-1})$,
 $f_i(\cdot), p(\cdot)\in C^{0}$ and furthermore, let the system
\begin{eqnarray}\label{plnat_ins_ex_psi=0}
\dot{x}_i&=&f_i(x_1,\dots,x_{i})+x_{i+1}, \ i=\{1,\dots,n-2\}\nonumber\\
\dot{x}_{n-1}&=&p(x_1,\dots,x_{n-1})+\upsilon
\end{eqnarray}
state be bounded for any $\upsilon\in L_\infty$ (i.e. system
(\ref{plnat_ins_ex_psi=0}) has  the {\it bounded input - bounded
state} property).
% it is known that {\it input-to-state stable}
% systems \cite{Sontag95} automatically
% satisfy this property
%, though the converse is not always true).
Then
%boundedness of the
%state follows immediately if $\upsilon(t)=\psi(\bfx(t))$ is
%bounded.
% (notice that neither
%stability of an equilibrium nor attractivity of a closed orbit is
%explicitly required in this case)
%Hence,
for the system of equations (\ref{plnat_ins_ex}), it is
sufficient that system (\ref{plnat_ins_ex_psi=0}) is
input-to-state stable with respect to input $\upsilon$ to satisfy
Assumption \ref{boundedness_psi}\footnote{For plants with
uncertainties in functions $f_i(\cdot)$, the time-dependent target
manifold $\psi(\bfx,t)=0$ may be required for control design. This
will be clarified later in the proof of Theorem
\ref{theorem_full_cascade} and in the examples section.} with
$\psi(\bfx)=x_n-p(x_1,\dots,x_{n-1})+f_{n-1}(x_1,\dots,x_{n-1})$.

\begin{assume}[Regularity]\label{Regularity} For any  $\bfx\in R^n$ and $t>0$ functions $\psi(\bfx,t)$ and $\bfg(x)$
satisfy the following inequality
$\left|L_{\bfg}\psi(\bfx,t))\right|>\delta_1>0$.
\end{assume}

Assumption \ref{Regularity} ensures the existence of feedback that
transforms the original system into that of the error model with
respect to the variable $\psi(\bfx,t)$. Let Assumption
\ref{Regularity} hold; consider
\begin{eqnarray}\label{plant_psi_ins}
\dpsi=L_{\bff}\psi(\bfx,t)+L_{\nuvec(\bfx,\thetavec)}\psi(\bfx,t)+\left(L_{\bfg}\psi(\bfx,t)\right)
u + {\pd \psi(\bfx,t)}/{\pd t}.
\end{eqnarray}
Because of Assumption \ref{Regularity} there exists the control
input
\begin{eqnarray}\label{control}
u(\bfx,\hat{\thetavec},t)=(L_{\bfg(\bfx)}\psi(\bfx,t))^{-1}\left(-\varphi(\psi)-L_{\bff}\psi(\bfx,t)-L_{\nuvec(\bfx,\hat{\thetavec})}\psi(\bfx,t)-{\pd\psi(\bfx,t)}/{\pd
t}\right),
\end{eqnarray}
where $\hat{\thetavec}\in\Omega_{\hat\theta}\subset R^d$ -- is a
vector of controller parameters that transforms
(\ref{plant_psi_ins}) into
\begin{eqnarray}\label{plant2_psi_ins}
\dpsi=-\varphi(\psi) + z(\bfx,\thetavec,t) -
z(\bfx,\hat{\thetavec},t)
\end{eqnarray}
where
$z(\bfx,\thetavec,t)=L_{\nuvec(\bfx,\thetavec)}\psi(\bfx,t)$. Let
the closed loop system satisfy some additional requirements:
\begin{assume}[Certainty Equivalence]\label{certainty_equivalence_part1}
For any $\thetavec\in \Omega_{\theta}$ there exists
$\hat{\thetavec}^{\ast}\in \Omega_{\hat{\thetavec}}\subset R^d$,
such that for all $\bfx\in R^n$, $t\in R_+$  the following
equivalence holds
\begin{eqnarray}\label{certainty_equivalence_1}
%\dpsi+\varphi(\psi) + z(\bfx,\thetavec,t) -
%z(\bfx,\hat{\thetavec}^{\ast},t)=\dpsi+\varphi(\psi)=0
 \frac{\pd \psi(\bfx,t)}{\pd
\bfx}[\bff(\bfx)+\nuvec(\bfx,\thetavec)+\bfg(\bfx)u(\bfx,\hat{\thetavec}^{\ast},t))]+\frac{\pd
\psi(\bfx,t)}{\pd t}+\varphi(\psi)
 = \dpsi(\bfx,\thetavec,\hat{\thetavec}^{\ast},t)+\varphi(\psi)=0.
\end{eqnarray}
\end{assume}
It is clear  that if  Assumption \ref{Regularity} holds then
Assumption \ref{certainty_equivalence_part1} is automatically
satisfied. According to Assumptions
\ref{Regularity} and \ref{certainty_equivalence_part1}, it follows
that $z(\bfx,\hat{\thetavec}^{\ast},t)=z(\bfx,\thetavec,t)$ for
any $\bfx\in R^n$ and time $t>0$.
\begin{assume}[Stability of the Target Dynamics]\label{certainty_equivalence_part2}
Function $\varphi(\psi)$ in (\ref{certainty_equivalence_1})
satisfies
\begin{equation}\label{certainty_equavalence_2}
\varphi(\psi)\in C^0, \ \varphi(\psi)\psi>0 \ \forall \psi\neq0, \
\lim_{\psi\rightarrow\infty}\int_0^{\psi}\varphi(\varsigma)d\varsigma=\infty.
\end{equation}
\end{assume}
\begin{assume}[Monotonicity and Linear Growth Rate in Parameters]\label{alpha}
There exists function $\alphavec(\bfx,t): R^{n}\times R\rightarrow
R^d$ such that
$(z(\bfx,\hat{\thetavec},t)-z(\bfx,\hat{\thetavec}^{\ast},t))(\alphavec(\bfx,t)^{T}(\hat{\thetavec}-\hat{\thetavec}^{\ast}))>0
\  \forall \ z(\bfx,\hat\thetavec^{\ast},t)\neq
z(\bfx,\hat{\thetavec},t)$. Furthermore, the following holds:
$|z(\bfx,\hat{\thetavec},t)-z(\bfx,\hat{\thetavec}^{\ast},t)|\leq
D |\alphavec(\bfx,t)^{T}(\hat{\thetavec}-\hat{\thetavec}^{\ast})|,
\ D>0$.
\end{assume}
%{\bf Either $D\in
%R_+$ or $D>0$!  There is really no need for both of them!  I'd use $D>0$.  Please fix this type of redundancy
%in other places throughout the text.}

Assumption \ref{certainty_equivalence_part1} (certainty
equivalence or matching condition) simply states that for every
unknown $\thetavec^{\ast}\in\Omega_{\theta}$ there exists a vector
of controller parameters
$\hat{\thetavec}^{\ast}(\thetavec^{\ast})\in\Omega_{\hat{\thetavec}}$
such that the system dynamics with this control function satisfies
the following equation $\dpsi=-\varphi(\psi)$. Assumption
\ref{certainty_equivalence_part2} specifies the properties of
function $\varphi(\psi)$, thus stipulating asymptotic stability of
manifold $\psi(\bfx,t)=0$ for
$\hat{\thetavec}=\hat{\thetavec}^{\ast}$ and ensuring unbounded
growth of integral $\int_0^{\psi}\varphi(\varsigma)d\varsigma$ as
$\psi\rightarrow \infty$.
% Notice that, due to Assumption
%\ref{certainty_equivalence_part2}, function
%$\int_0^{|\psi|}\varphi(\varsigma)d\varsigma$ is continuous and
%strictly monotonic in $|\psi|$ and therefore has the monotonic
%inverse $\Lambda: R_+\rightarrow R_+$
%\begin{eqnarray}\label{inverse_int_varphi}
%\Lambda\left(\int_0^{|\psi|}\varphi(\varsigma)d\varsigma\right)=|\psi|.
%\end{eqnarray}
%This fact will be used later to estimate the performance of the
%adaptive system.
Assumption \ref{alpha} is given to specify the admissible
nonlinear parameterization of the controller. For
linearly parameterized plants this assumption is automatically
satisfied. Sometimes we will further restrict the nonlinear in
parameter functions by:
\begin{assume}\label{alpha_upper_bound}
There exists a positive constant $D_1>0$ such that for any $\bfx$,
$\hat\thetavec$, $\hat{\thetavec}^{\ast}$, $t>0$ the following
inequality holds:
$|z(\bfx,\hat{\thetavec},t)-z(\bfx,\hat{\thetavec}^{\ast},t)|\geq
D_1
|\alphavec(\bfx,t)^{T}(\hat{\thetavec}-\hat{\thetavec}^{\ast})|$.
\end{assume}
Throughout the paper we will assume that functions
$\alphavec(\bfx,t)$ and $u(\bfx,\hat{\thetavec},t)$ are both
bounded in $t$. For the sake of convenience and if not stated
overwise we will also assume that functions
$\alphavec(\bfx,t),\psi(\bfx,t)$ are differentiable as many times
as necessary if differentiation is required to design the
algorithm. In addition, we will use the term ''smooth functions''
to denote those functions that belong to $C^{\infty}$. Though it
is not necessary at all for us to require existence of infinitely
many derivatives of the functions that we refer to as smooth in
the paper, this notational agreement will free the presentation of
numerous insignificant details in the formulations. Furthermore,
along with already defined $L_2$ and $L_\infty$ norms we will use
the following notation:
$\|\hat{\thetavec}^{\ast}-\hat\thetavec(t)\|^{2}_{\Gamma^{-1}}=$
$(\hat{\thetavec}^{\ast}-\hat\thetavec(t))^{T}\Gamma^{-1}(\hat{\thetavec}^{\ast}-\hat\thetavec(t))$.

As mentioned in the introduction, we propose to design the
adaptive algorithms in two steps: 1) search for the suitable
augmentation ensuring the desired properties of the control
system, and 2) find the appropriate realization of this algorithm
in finite form. Therefore we start with the choice of tuning
errors $\tilde\psi(\bfx,t)$ and operators $\mathcal{A}(\bfx,t)$
for the class of algorithms given by formula (\ref{alg0}): $
\dot{\hat{\thetavec}}=\Gamma
\tilde{\psi}(\bfx,t)\mathcal{A}(\bfx,\hat{\thetavec})$. As a
candidate for the augmented error $\tilde\psi(\bfx,t)$ we select
the following $\tilde\psi(\bfx,t)=\dpsi+\psi(\bfx,t)$. It has been
shown in \cite{tpt2002_at,tpt2003_tac} that the algorithm
\begin{equation}\label{alg1_a}
\dot{\hat{\thetavec}}=\Gamma(\varphi(\psi)+\dpsi)\alphavec(\bfx,t),
\ \Gamma>0
\end{equation}
with control (\ref{control}) guarantee that
$\psi(\bfx(t),t)\rightarrow 0$ as $t\rightarrow\infty$ for the
closed loop system
\begin{eqnarray}\label{ONO}
\dot{\bfx}_1&=& \bff_1(\bfx)+\bfg_1(\bfx)u, \ \  \dot{\bfx}_2=
\bff_2(\bfx)+\nuvec(\bfx,\thetavec)+\bfg_2(\bfx)u \nonumber \\
%\dpsi&=&-\varphi(\psi) +
%z(\bfx,\thetavec,t) -
%z(\bfx,\hat{\thetavec},t)\nonumber \\
\dot{\hat\thetavec}&=&\Gamma
(\varphi(\psi)+\dpsi)\alphavec(\bfx,t), \ \Gamma>0
\end{eqnarray}
under Assumptions \ref{boundedness_psi},
\ref{certainty_equivalence_part1} -- \ref{alpha}. In addition, it
is possible to show that system (\ref{ONO}) has better performance
than that of the known schemes. This follows from the next theorem
{ (see also Proposition 1 below):}

\begin{thm}\label{L_2 dpsi_theorem} Let system (\ref{ONO}) be given and
Assumptions \ref{Regularity}--\ref{alpha} hold. Then for system
(\ref{ONO}) the following hold:

P1) $\varphi(\psi(t))\in L_2$, $\dpsi(t)\in L_2$;

P2)
$\|\hat{\thetavec}^{\ast}-\hat\thetavec(t)\|^{2}_{\Gamma^{-1}}$ is
non-increasing;

P3) $z((\bfx,\thetavec,t)-z(\bfx,\hat{\thetavec}(t),t))\in L_2$.

Furthermore,
\begin{eqnarray}\label{L_2_L_inf_performance}
 \|\varphi(\psi)\|_{2}^2 &\leq& 2
 Q(\psi)+\|\hat{\thetavec}(0)-\hat{\thetavec}^{\ast}\|^{2}_{(2D\Gamma)^{-1}},
 \ \
 \|\dpsi\|_{2}^2 \leq  2 Q(\psi)+\|\hat{\thetavec}(0)-\hat{\thetavec}^{\ast}\|^2_{(2D\Gamma)^{-1}} \nonumber \\
\|\psi\|_{\infty} &\leq& \Lambda
\left(Q(\psi)+\|\hat{\thetavec}(0)-\hat{\thetavec}^{\ast}\|^2_{(4D\Gamma)^{-1}}\right),
\end{eqnarray}
where
$Q(\psi)=\int_{0}^{\psi(\bfx(0),0)}\varphi(\varsigma)d\varsigma$
and $\Lambda(d)=\max_{|\psi|}\{|\psi| \ | \
\int_{0}^{|\psi|}\varphi(\varsigma)d\varsigma=d\}$.

If Assumption \ref{boundedness_psi} is satisfied and function
$z(\bfx,\hat{\thetavec},t)$ is locally bounded with respect to
$\bfx$, $\hat{\thetavec}$ and uniformly bounded with respect to
$t$, then

P4) trajectories of the system are bounded and
$\psi(\bfx(t))\rightarrow 0$ as $t\rightarrow\infty$;

If in addition functions $\varphi, z(\bfx,\thetavec,t)\in C^1$;
derivative $ \pd {z(\bfx,\thetavec,t)}/{\pd t}$ is uniformly
bounded in $t$; function $\alphavec(\bfx,t)$ is locally bounded
with respect to $\bfx$ and uniformly bounded with respect to $t$,
then

P5) $\dpsi\rightarrow 0$ as  $t\rightarrow \infty$;
$z((\bfx,\thetavec^{\ast},t)-z(\bfx,\hat{\thetavec}(t),t))\rightarrow
0$
 as  $t\rightarrow \infty$.
\end{thm}

The formal proof of the theorem is given in Appendix
2\footnote{%While the proof is similar to that of Theorem 1 in
%\cite{tpt2003_tac} for the error model
%$\dpsi=-\varphi(\psi)+z(\bfx,\thetavec)-z(\bfx,\hat{\thetavec})$,
%the properties P1)--P3), P5) were not explicitly formulated there.
Robustness and some other properties of algorithms (\ref{alg1_a})
are discussed in Appendix 1.}.

%\begin{rem}\label{rem_decr_uncert} \normalfont In addition to the fact that the $L_2$ and $L_{\infty}$
%norm bounds for $\dpsi$, $\dpsi$ reflect the transient behavior of
%system (\ref{ONO}) and that the algorithms are capable of dealing
%with nonconvexly parameterized plants, there are several
%additional advantages of the considered scheme. In particular,
%properties P2) and P3). Property P2) states that parametric
%uncertainty is not increasing with time. This provides us with an
%obvious improvement over the existing algorithms. Furthermore, it
%follows from P3) that
%$z((\bfx,\thetavec^{\ast},t)-z(\bfx,\hat{\thetavec}(t),t))\rightarrow
%0$ as $t\rightarrow\infty$ given the assumptions on function
%$z(\cdot)$. Taking into account equation (\ref{control}) the
%following identity results:
%$\lim_{t\rightarrow\infty}u(\bfx,\hat{\thetavec},t)=u(\bfx,\hat{\thetavec}^{\ast},t)$.
%This means not only that the control goal $\psi(\bfx,t)=0$ will
%asymptotically be reached but also that the control input as a
%function of time will asymptotically approache its ``ideal"
%values. Furthermore, if the control input
%$u(\bfx,\hat{\thetavec}^{\ast},t)\in L_2$ along the system
%solutions, then by property P5) the same applies to
%$u(\bfx,\hat{\thetavec},t)$ ($u(\bfx,\hat{\thetavec},t)\in L_2$ as
%the sum of two functions from $L_2$).
%\end{rem}
Notice that the function $\varphi(\psi)$ is nonlinear and that its
shape  influences the $L_2$ and $L_\infty$ norm bounds for $\psi$
and $\dpsi$. Because of this, according to
(\ref{L_2_L_inf_performance}) it is possible to improve the
performance of the system with respect to $L_2$ and $L_\infty$
bounds by varying the function $\varphi(\psi)$. The bounds
obtained for the $L_{\infty}$ norms may be improved further for
the case when the function $\varphi(\psi)$ is linear in $\psi$.
This is not too severe a restriction as the choice of function
$\varphi(\psi)$ is always up to the designer. Performance
characteristics of the system for this case are formulated in
Proposition \ref{Exponential_Convergence}.

\begin{prop}[Exponential Convergence]\label{Exponential_Convergence}
Let Assumptions \ref{Regularity}--\ref{alpha} hold and
$\varphi(\psi)=K \psi$, $K>0$. Then

P6) function $\psi(\bfx(t),t)$ converges exponentially fast into
the domain $|\psi(\bfx(t),t)|\leq
0.5\sqrt{\|\hat{\thetavec}(0)-\hat{\thetavec}^{\ast}\|^2_{(KD\Gamma)^{-1}}}$.
Specifically, the following holds:   $|\psi(\bfx(t),t)|\leq
 |\psi(\bfx(0),0)| e^{-K t} + 0.5\sqrt{\|\hat{\thetavec}(0)-\hat{\thetavec}^{\ast}\|^2_{(KD\Gamma)^{-1}}}$

Furthermore, let Assumption \ref{boundedness_psi} hold, function
$\alphavec(\bfx,t)$ be locally bounded with respect to $\bfx$ and
uniformly bounded in $t$; for any bounded $\bfx$ there exist
$D_1>0$ such
 that
%\label{least_growth_alpha}
$|z(\bfx,\hat{\thetavec},t)-z(\bfx,\hat{\thetavec}^{\ast},t)|\geq
D_1
|\alphavec(\bfx,t)^{T}(\hat{\thetavec}-\hat{\thetavec}^{\ast})|$,
%\end{eqnarray}
function $\alphavec(\bfx,t)$ is persistently exciting:
\begin{eqnarray}\label{PE}
\exists L>0, \ \delta>0: \
\int_{t}^{t+L}\alphavec(\bfx(\tau),\tau)\alphavec(\bfx(\tau),\tau)^{T}d\tau\geq
\delta I \ \forall t>0,
\end{eqnarray}
where $I\in R^{d\times d}$ -- identity matrix. Then

P7) both $\psi(\bfx(t),t)$ and
$\|\hat{\thetavec}-\hat{\thetavec}^{\ast}\|$ converge
exponentially fast to the origin.
%If in addition $\Gamma=\gamma I, \  \gamma \in R_+
%\ \gamma>0$ then
%\begin{eqnarray}\label{rates}
%& & |\psi(\bfx(t),t)|\leq \sqrt{2}|\psi(\bfx(0),0)| e^{-K t} +
%\|\hat{\thetavec}(0)-\hat{\thetavec}^{\ast}\|e^{- D_2 t}\nonumber \\
%& & \|\hat{\thetavec}(t)-\hat{\thetavec}^{\ast}\|\leq
%\|\hat{\thetavec}(0)-\hat{\thetavec}^{\ast}\|D_3
%e^{-\frac{D_1\delta}{L}t}.
%\end{eqnarray}
\end{prop}
It follows from Proposition \ref{Exponential_Convergence} that if
$\varphi(\psi)=K\psi$ then the estimate of the upper bound
$\sup_{t\geq t' }|\psi(\bfx(t),t|$ as a function of $t'$ for
$\psi(\bfx(t),t)$ in system (\ref{ONO}) exponentially converges
into the domain determined by the parametric uncertainty, the
values of  controller parameters $K$, and adaptation gain
$\Gamma$. Notice that this domain can be made arbitrary small,
subject to the choice of the values of $K$ and $\Gamma$. The rates
of convergence are given by P6). In the case of persistent
excitation an even stronger property is established. The system is
shown to be {\it exponentially stable} with respect to the target
manifold $\psi(\bfx,t)=0$ and point
$\hat{\thetavec}(t)=\hat{\thetavec}^{\ast}$.
% if additional
%assumption (\ref{least_growth_alpha}) holds for functions
%$z(\bfx,\thetavec,t)$ and $\alphavec(\bfx,t)$.

% Inequality
%(\ref{least_growth_alpha}) along with Assumption \ref{alpha} state
%should be bounded (both from above and below) by two functions
%that are linear in $\hat{\thetavec}-\thetavec$:
%$D\alphavec(\bfx,t)^{T}(\hat{\thetavec}-\thetavec)$ and
%$D_1\alphavec(\bfx,t)^{T}(\hat{\thetavec}-\thetavec)$ for any
%bounded $\bfx$.

Despite properties of algorithms (\ref{alg1_a}) such as improved
transient performance of the closed loop system and their ability
to deal with nonconvexly parameterized models, these algorithms
 are not realizable in the form of differential
equations, as they depend on unknown parameters explicitly. It was
proposed in \cite{tpt2003_tac} to use special filters to estimate
$\dpsi$. While the approach of \cite{tpt2003_tac} is acceptable
for systems with nonconvex parameterization, control system
performance may be suboptimal due to estimation errors. The
question is how to realize algorithms (\ref{alg1_a}) in a form
that depends neither on time-derivative $\dpsi$ nor on its
filtered estimate explicitly, nor on anything implying knowledge
of unknown parameters $\thetavec$.  Our solution, as mentioned in
Section 1, is to use the finite form (\ref{fin_form_ours}) of
adaptive algorithms instead of the differential form
(\ref{alg1_a}).
%:
%\begin{eqnarray}
%\hat{\thetavec}(\bfx(t),t)&=&\hat{\thetavec}_P(\bfx,t)+\hat{\thetavec}_I(t),  \ \dot{\hat{\thetavec}}_I=\mathcal{A}_2(\bfx,\hat{\thetavec},t),\nonumber\\
%\hat{\thetavec}_P(\bfx,t)&=&\mathcal{A}_1(\bfx,t)\nonumber
%\end{eqnarray}
In the next section we study under what conditions algorithms
(\ref{alg1_a}) can be represented in finite form
(\ref{fin_form_ours}).

\section{Adaptive Algorithms in Finite Form}

The outline of the section is as follows.  We start from a general
case and formulate the conditions ensuring the realization of
algorithm (\ref{alg1_a}) in finite form explicitly, i.e., without
any filters and further transformations of the closed-loop system.
The conditions we impose involve the existence of the solutions of
a system of partial differential equations. It is nontrivial to
check these assumptions for nonlinear model (\ref{plant_ins}).
That they hold, however, can be demonstrated for some special
combinations of plant models and goal functions $\psi(\bfx,t)$.
Further, we consider extension of the proposed method to a broader
class of nonlinear systems including systems with low-triangular
structure.

%We provide examples of those nonlinear systems that satisfy these
%criteria. Further, we consider cascades of systems such that each
%subsystem in the cascade has a stabilizing adaptive control
%algorithm in finite form. For these systems we provide an
%iterative procedure, resulting in a finite-form adaptive control
%algorithm for the whole cascade.

%Next, we consider the case when an extension of the plant dynamics
%makes it possible to achieve a finite-form realization of the
%adaptive control algorithms for the extended system. We suggest to
%embed the plant dynamics into the extended system for which the
%conditions sufficient for finite-form realization are always met.
%By doing so we eliminate the necessity to find a solution of
%partial differential equations to realize the algorithm.

\subsection{Explicit Realization}

Let us assume that in addition to the Assumptions
\ref{boundedness_psi}--\ref{alpha}, that are sufficient for system
(\ref{ONO}) to have the properties P1)--P7), the following hold
\begin{assume}[Explicit realization condition]\label{Realizability_1} For the given functions
$\alphavec(\bfx,t)$ and $\psi(\bfx,t)$ there exists function
$\Psi(\bfx)$ such that the following hold:
\begin{equation}\label{realizability_eq}
\Psi(\bfx): {\pd \Psi (\bfx,t)}/{\pd \bfx_2}=\psi(\bfx,t)({\pd
\alphavec(\bfx,t)}/{\pd \bfx_2})
\end{equation}
\end{assume}
Then realizations of the adaptive scheme described by equations
(\ref{alg1_a}) follow from the next theorem.

\begin{thm}\label{theorem_realization1} Let Assumption
\ref{Realizability_1} hold. Then there is a finite-form
realization of the algorithms (\ref{alg1_a}):
\begin{eqnarray}\label{fin_forms_ours_tr1}
\hat{\thetavec}(\bfx,t)&=&\Gamma(\hat{\thetavec}_P(\bfx,t)+\hat{\thetavec}_I(t));
\ \hat{\thetavec}_P(\bfx,t)=
\psi(\bfx,t)\alphavec(\bfx,t)-\Psi(\bfx,t)\nonumber \\
\dot{\hat{\thetavec}}_I&=&\varphi(\psi(\bfx,t))\alphavec(\bfx,t)+{\pd
\Psi(\bfx,t)}/{\pd t}-\psi(\bfx,t)({\pd
\alphavec(\bfx,t)}/{\pd t})-\nonumber\\
%-\frac{\pd \psi(\bfx(\tau),\tau)}{\pd
%\tau}\alphavec(\bfx(\tau),\tau)-\psi(\bfx(\tau),\tau)\frac{\pd
%\alphavec(\bfx(\tau),\tau)}{\pd \tau}\right)d\tau
& &  (\psi(\bfx,t)L_{\bff_1} \alphavec(\bfx,t)-L_{\bff_1}
\Psi(\bfx,t))-(\psi(\bfx,t)L_{\bfg_1}\alphavec(\bfx,t)-L_{\bfg_1}
\Psi(\bfx,t))u(\bfx,\hat{\thetavec},t)
\end{eqnarray}
\end{thm}
\begin{rem}\label{time_varying_K_realization} \normalfont It is
easy to see from (\ref{fin_forms_ours_tr1}) and the theorem proof
that realization of the algorithms
\begin{eqnarray}\label{alg_timevar}
\dot{\hat{\thetavec}}=\Gamma(\dpsi+\beta(\bfx,t))\alphavec(\bfx,t),
\end{eqnarray}
where $\beta(\bfx,t)$ is to guarantee at least the existence of
 solutions for the closed loop system,  is also possible.
Indeed, in order to realize these algorithms  it is sufficient to
replace equations for $\dot{\hat\thetavec}_I$ in
(\ref{fin_forms_ours_tr1}) by the following:
\begin{eqnarray}\label{fin_forms_ours_tr2}
\dot{\hat{\thetavec}}_I&=&\beta(\bfx,t)\alphavec(\bfx,t)+{\pd
\Psi(\bfx,t)}/{\pd t}-\psi(\bfx,t)({\pd
\alphavec(\bfx,t)}/{\pd t})-\nonumber\\
%-\frac{\pd \psi(\bfx(\tau),\tau)}{\pd
%\tau}\alphavec(\bfx(\tau),\tau)-\psi(\bfx(\tau),\tau)\frac{\pd
%\alphavec(\bfx(\tau),\tau)}{\pd \tau}\right)d\tau
& &  (\psi(\bfx,t)L_{\bff_1} \alphavec(\bfx,t)-L_{\bff_1}
\Psi(\bfx,t))-(\psi(\bfx,t)L_{\bfg_1}\alphavec(\bfx,t)-L_{\bfg_1}
\Psi(\bfx,t))u(\bfx,\hat{\thetavec},t),
\end{eqnarray}
One particular case of function
$\beta(\bfx,t)=(1+\delta(t))\psi(\bfx,t)$, $\delta: R_+\rightarrow
R_+$, $\delta\in C^0$ will be used later to show existence of the
adaptive control algorithms for nonlinearly parameterized plants
in the low-triangular form.

Theorem \ref{theorem_realization1} provides us with an answer to
the question of existence of realizable algorithms that satisfy
differential equations (\ref{alg1_a}), thus ensuring the
properties formulated in Theorem \ref{L_2 dpsi_theorem} and
Proposition \ref{Exponential_Convergence}.
% What is important is
%that the number of integrators for both algorithms (\ref{alg1_a})
%and (\ref{fin_forms_ours_tr1}) is the same.
 The disadvantage,
however, is that the functions $\Psi(\bfx,t)$ in Assumption
\ref{Realizability_1} are not easy to find. %, if they exist at all.
Existence of such functions itself is another nontrivial issue.
For instance, if $\dim{\bfx_2}=n$ and functions $\psi(\bfx,t)$,
$\alphavec(\bfx,t)$ do not depend explicitly on time $t$, then the
necessary conditions for the function $\Psi(\bfx)$ to exist is the
symmetry of all matrices $\frac{\pd}{\pd
\bfx}\left(\psi(\bfx)\frac{\pd \alpha_i(\bfx)}{\pd \bfx}\right), \
i\in\{1,\dots,n\}$. Nevertheless, despite difficulties in finding
those functions $\Psi(\bfx,t)$ that satisfy Assumption
\ref{Realizability_1}, there are several classes of dynamical
systems with certain structural properties that immediately reduce
Assumption \ref{Realizability_1} to more easily verifiable
requirements.
\end{rem}
\begin{cor}[Single-dimension uncertainty-dependent partition]\label{single_dim_par} Let $\dim(\bfx_2)=1$ and
function

\noindent $\psi(\bfx,t)\pd \alphavec(\bfx,t)/\pd x_n$  be
Riemann-integrable with respect to $x_n$, i.e., the following
integral exists
\begin{eqnarray}\label{Psi1}
\Psi(\bfx,t)=\int \psi(\bfx,t)\frac{\pd \alphavec(\bfx,t)}{\pd
x_n} d x_n
\end{eqnarray}
Then there is a finite-form realization of algorithms
(\ref{alg1_a}).
\end{cor}

\begin{rem}\label{rem_SDUP}\normalfont Corollary
\ref{single_dim_par} allows us to turn the problem of searching
for a function $\Psi(\bfx,t)$ satisfying equation
(\ref{realizability_eq})  into a problem of existence of the
indefinite integral of a function with respect to a single scalar
argument. It is clear from (\ref{Psi1}) that any one-dimensional
system with integrable $\psi(x)\frac{\pd \alpha(x)}{\pd x}$ has a
finite-form realization. An interesting  example is the class
of systems described by the following differential equations:
\begin{eqnarray}\label{example_t1}
\dot{x}_i&=&f_i(\bfx)+g_i(\bfx)u, \ \ \ \ i=1,\dots,n-1\nonumber \\
\dot{x}_n&=&f_n(\bfx)+\nu(\bfx,\thetavec)+g_n(\bfx)u,
\end{eqnarray}
where function $\nuvec(\bfx,\thetavec)$ satisfies Assumption
\ref{alpha}, which in turn is automatically satisfied if
$\nuvec(\bfx,\thetavec)$ linearly parameterized or
$\nu(\bfx,\thetavec)=\nuvec(\bfx^{T}\thetavec)$ and
$\nuvec(\cdot)$ is monotonic and belongs to a sector. In practice,
the indefinite integral in (\ref{Psi1}) can also be replaced by
$\Psi(\bfx,t)=\int_{x_n(0)}^{x_n(t)} \psi(\bfx,t)\frac{\pd
\alphavec(\bfx,t)}{\pd x_n} d x_n$.
%{\bf Don't we integrate over $x_n$, not $\xi$?}
\end{rem}

Equations of type (\ref{example_t1}) describe a class of dynamical
systems in which the uncertainties are concentrated in a single
equation.  There are many mechanical systems described by
equations (\ref{example_t1}) satisfying Assumption
\ref{Realizability_1} (e.g., the simple classical equations for
shaft dynamics with unknown load torque, rotating platforms and
pendulums).
%  {\bf Should also mention that some $g_i(\bfx)$ can be
%equal to zero.}
For the case where the uncertainty is a single scalar and function
$\alpha(\bfx)$ can be chosen as the goal function
$\psi(\bfx)=\alpha(\bfx)$, finite-form realization is also
possible ($\Psi(\bfx)=\frac{1}{2}\alpha^2(\bfx)$).
%Furthermore, it is
%desirable to notice that in order to satisfy the corollary
%assumptions it is enough to require that function
%\[
%\frac{\pd}{\pd \bfx_2}\psi(\bfx_1\oplus\bfx_2,t)\frac{\pd
%\alpha(\bfx_1\oplus\bfx_2,t)}{\pd \bfx_2}=\beta(\bfx_1,x_)
%\]
Another class of dynamical systems that automatically satisfy
Assumption \ref{Realizability_1} is given by the following
corollary.

\begin{cor}[Independence on partition $\mathcal{L}_2$]\label{independent_partition_cor} Let function $\alphavec(\bfx,t)$
be independent on $\mathcal{L}_2$, i. e., for any
$\bfx_2\in\mathcal{L}_2$ the following holds: ${\pd
\alphavec(\bfx,t)}/{\pd \bfx_2}={\pd
\alphavec(\bfx_1\oplus\bfx_2,t)}/{\pd \bfx_2}=0$ then there is a
finite-form realization of algorithms (\ref{alg1_a}).
\end{cor}
Corollary \ref{independent_partition_cor} conditions are
equivalent to the fact that the plant dynamics can be described by
system
\begin{eqnarray}\label{independent_partition_decomp}
\dot{\bfx}_1&=&\bff_1(\bfx_1\oplus\bfx_2)+\bfg_1(\bfx_1\oplus\bfx_2)u\nonumber \\
\dot{\bfx}_2&=&\bff_2(\bfx_1\oplus\bfx_2)+\nuvec(\bfx_1,\thetavec)+\bfg_2(\bfx_1\oplus\bfx_2)u
\end{eqnarray}
and that ${\pd \psi(\bfx_1\oplus\bfx_2,t)}/{\pd
\bfx_2}=\lambda(\bfx_1,t)$, where $\lambda: R^{m}\times
R\rightarrow R^{n-m}$ is a differentiable function with known
derivative $\pd \lambda(\bfx_1,t)/\pd t$. Therefore it is possible
to derive from Corollary \ref{independent_partition_cor} that
every error model:
%\begin{eqnarray}%\label{dpsi_pre_embedd}
$\dpsi=-\varphi(\psi)+z(\omegavec(t),\thetavec)-z(\omegavec(t),\hat{\thetavec})$,
%\end{eqnarray}
where $\omegavec(t): R\rightarrow R^n, \ \omegavec\in C^1$ is a
function with known time-derivatives $\dot{\omegavec}(t)$,
satisfies the sufficient conditions for realization of algorithm
(\ref{alg1_a}) in finite form. Indeed, this follows directly from
Assumption \ref{alpha}, as functions $\alphavec(\bfx,t)$ in this
case are independent of $\bfx$. Therefore if the derivatives
$\dot{\alphavec}(t)$ are known, the finite form realization
follows immediately from
\begin{eqnarray}
\hat{\thetavec}(\bfx,t)&=&\Gamma(\hat{\thetavec}_P(\bfx,t)+\hat{\thetavec}_I);
\ \
\hat{\thetavec}_P(\bfx,t)=\psi(\bfx,t)\alphavec(t)\nonumber \\
\dot{\hat{\thetavec}}_I&=&\varphi(\psi(\bfx,t))\alphavec(\bfx,t)-\psi(\bfx,t)\dot{\alphavec}(t)\nonumber
\end{eqnarray}
This fact,
%{\bf Which one?  Corollary 2?},
along with decomposition (\ref{independent_partition_decomp}),
will be used later in Section 3.2.
%when
%dealing with the approximate realizations of algorithms
%(\ref{alg1_a}).

%Assumption \ref{Realizability_1} can be relaxed for linearly
%parameterized systems (Fradkov, Stotsky).

So far, simplified conditions for the existence of the adaptive
algorithms in finite form  were derived from Theorem
\ref{theorem_realization1} for those classes of nonlinear systems
that have certain structural properties, such as single dimension
uncertainty-dependent partition (Corollary \ref{single_dim_par}
and equation (\ref{example_t1})) or independence of
$z(\bfx,\thetavec,t)$ on uncertainty-dependent partition $\bfx_2$
(Corollary \ref{independent_partition_cor}). These structural
properties allowed us to reduce Assumption \ref{Realizability_1}
to integrability of a function with respect to a single scalar
argument for a class of nonlinear systems.
%Based on these results it is possible to extend applicability of our
%method to a broader class of systems.
Taking these results into account, in the next section we present
a technique that allows us to extend our method to a broader class
of systems.

%This rather simple test, however, is only sufficient (but not
%necessary) to clarify existence of adaptive control algorithms
%with improved transient behavior and ability to deal with
%nonconvex parameterization. On the other hand, it is natural to
%expect that there are classes of systems that can be reduced to
%the considered cases, for which algorithms in finite form already
%exist. In the next section we present a technique that allows us
%to transform a nonlinear dynamical system into a form that obeys
%these sufficient conditions.

\subsection{Asymptotic Design via Embedding}

The main idea behind the extension of our results to a broader
class of nonlinear systems is as follows. Instead of trying to
find a general solution of  equation (\ref{realizability_eq}) in
Assumption \ref{Realizability_1} (which is a nontrivial task even
if such solution exists), we transform the original equations into
a form that satisfies much weaker requirements considered in
Corollaries \ref{single_dim_par} and
\ref{independent_partition_cor}. This transformation should not
necessarily be a one-to-one diffeomorphism, but the control goal
reaching in the new state space should guarantee reaching the
control goal of the original system. One way to assure this is to
embed the original system dynamics into one of a higher order, for
which a finite form realization of the adaptive control algorithms
is possible.

Let us represent the partitioned system (\ref{partitioned_plant})
in the following way:
\begin{eqnarray}\label{partitioned_plant_embed}
\dot{\bfx}_1&=&\bff_1(\bfx)+\bfg_1(\bfx)u\nonumber \\
\dot{\bfx}_2'&=&\bff_2'(\bfx)+\nuvec'(\bfx,\thetavec)+\bfg_2'(\bfx)u;
\ \
\dot{\bfx}_2''=\bff_2''(\bfx)+\nuvec''(\bfx,\thetavec)+\bfg_2''(\bfx)u,
\end{eqnarray}
where $\bfx_2'\oplus\bfx_2''=\bfx_2$, $\dim{\bfx_2'}=m_1, \
\dim{\bfx_2}''=n-m-m_1, \ 0 \leq m_1\leq n-m$.
% and functions
%$\bff_2',\bff_2'',\bfg_2',\bfg_2'',\nuvec',\nuvec''$ are defined,
%respectively as follows:
%\[
% \bff_2'=\left(\begin{array}{l}
%             f_m\\
%             f_{m+1}\\
%             \cdots\\
%             f_{m+m_1}
%             \end{array}
%             \bff_2=\left(\begin{array}{l}
%             f_{m+1}\\
%             f_{m+2}\\
%             \cdots\\
%             f_n
%             \end{array}
%             \right) \ \ \bfg_1=\left(\begin{array}{l}
%             g_1\\
%             g_2\\
%             \cdots\\
%             g_m
%             \end{array}
%             \right) \ \ \
%             \bfg_2=\left(\begin{array}{l}
%             g_{m+1}\\
%             g_{m+2}\\
%             \cdots\\
%             g_n
%             \end{array}
%\]
Using the notations above, we introduce the following assumption
\begin{assume}\label{assume_embedd} There exist

1) a partition of the state vector $\bfx$: $
\bfx=\bfx_1\oplus\bfx_2'\oplus\bfx_2''$, %\ \dim{\bfx_2'}=m_1, \
%\dim{\bfx_2}''=n-m-m_1, \ 0 \leq m_1\leq n-m

2) a system of differential equations
\begin{eqnarray}\label{embed_add}
\dot{\xivec}&=&\bff_{\xi}(\bfx,\xivec,t); \
\bfy_{\xi}=\bfh_{\xi}(\xivec),
\end{eqnarray}
\[
\xivec\in R^r, \ \bff_{\xi}: R^{n}\times R^{r}\times
R_+\rightarrow R^r, \ \bff_\xi\in C^1;  \bfh_{\xi}:
R^{r}\rightarrow R^{n-m-m_1}, \ \bfh_{\xi}\in C^1;
\]

3) a function $\Psi(\tilde{\bfx},t)\in C^1$,
$\tilde{\bfx}=\bfx_1\oplus\bfx_2'\oplus\bfh_{\xi}$ such that the
following conditions hold
\begin{eqnarray}
&{\pd \Psi(\tilde{\bfx},t)}/{\pd
\bfx_2'}=\psi(\tilde{\bfx},t)({\pd\alphavec(\tilde{\bfx},t)}/{\pd
\bfx_2' })& \label{realizability_eq2} \\
%\end{eqnarray}
%\begin{eqnarray}\label{boundedness_xi}
&\bfx\in L_{\infty}\Rightarrow \xivec\in L_{\infty}&
\label{boundedness_xi}
\end{eqnarray}
for any $\thetavec\in\Omega_{\theta}$ and $t\in R_+$ along the
solutions of the original system (\ref{plant_ins}).
\end{assume}

In addition to Assumption \ref{assume_embedd}, we would like to
formulate two alternative assumptions which, if satisfied, will
result in two different adaptation schemes with different
performance and robustness properties.

\begin{assume}\label{assume_L_infty_embed} Let system
(\ref{embed_add}) be given and
\begin{eqnarray}\label{L_2_embed}
z(\bfx,\thetavec,t)-z(\tilde{\bfx},\thetavec,t)\in L_\infty; \ \
(\psi(\bfx,t)-\psi(\tilde{\bfx},t))({\pd
\alphavec(\tilde\bfx,t)}/{\pd \bfx_2'})\nuvec'(\bfx,\thetavec)\in
L_\infty.
\end{eqnarray}
along  the solutions of (\ref{plant_ins}), (\ref{embed_add}).
%\begin{eqnarray}\label{L_2_embed_2}
%\end{eqnarray}
\end{assume}

\begin{assume}\label{assume_L_2_embed} Let system
(\ref{embed_add}) be given and $\pd \alphavec(\tilde{\bfx},t)/\pd
\bfx_2' \equiv 0, \ \ or \ \ \psi(\bfx,t)=\psi(\tilde\bfx,t).$
Furthermore, let
$z(\bfx,\thetavec,t)-z(\tilde{\bfx},\thetavec,t)\in L_2$ for any
$\thetavec\in\Omega_\theta$, $t>0$ along  the solutions of
 (\ref{plant_ins}), (\ref{embed_add}).
\end{assume}

Sufficient conditions for the desired embedding follow from the
next theorem.
\begin{thm}[Embedding Theorem]\label{theorem_embedd} Let function $\psi(\bfx,t)$ be given and Assumptions
\ref{boundedness_psi}--\ref{alpha_upper_bound},
\ref{assume_embedd} hold for system (\ref{plant_ins}). Then for
the extended system
\begin{eqnarray}\label{plant_ins_embed}
\dot{\bfx}&=&\bff(\bfx) + \varthetavec(\bfx,\thetavec) +
\bfg(\bfx)u\nonumber \\
\dot{\xivec}&=&\bff_{\xi}(\bfx,\xivec,t); \
\bfy_{\xi}=\bfh_{\xi}(\xivec),
\end{eqnarray}
there exists control function
$u(\bfx,\bfh_{\xi},\hat{\thetavec},t)$
\begin{eqnarray}\label{control_embedd}
u(\bfx,\bfh_{\xi},\thetavec,t)&=&(L_{\bfg(\bfx)}\psi(\bfx,t))^{-1}\left(-\varphi(\psi)-L_{\bff}\psi(\bfx,t)-L_{\varthetavec(\tilde{\bfx},\hat{\thetavec})}\psi(\tilde{\bfx},t)-{\pd\psi(\bfx,t)}/{\pd
t}\right)
\end{eqnarray}
and adaptation algorithms\footnote{The adaptation algorithms that
guarantee properties P8) and P9) are given by equations
(\ref{fin_forms_embedd}) and (\ref{fin_forms_embedd_1}),
(\ref{fin_forms_embedd_2}) respectively in Appendix 2.}:
$\hat{\thetavec}(\tilde{\bfx},t)=\Gamma(\hat{\thetavec}_{P}(\tilde{\bfx},t)+\hat{\thetavec}_I(t)),
\ \Gamma>0$ such that the following statements hold:

P8) if $|\varphi(\psi)|\geq K |\psi|, \  K>0$ and Assumption
\ref{assume_L_infty_embed} holds then
$\psi(\bfx,t),\bfx,\xivec,\hat{\thetavec}\in L_\infty$;

P9) if Assumption \ref{assume_L_2_embed} holds then
$\psi(\bfx,t)\in L_2\cap L_\infty,\dpsi\in L_2, \ \
z(\tilde{\bfx},\thetavec,t)-z(\tilde{\bfx},\hat{\thetavec},t)\in
L_2, \ \ \bfx,\xivec\in L_\infty$;

if in addition derivatives $\pd \psi(\bfx,t)/\pd \bfx$, $\pd
\psi(\bfx,t)/\pd t$ are uniformly bounded in $t$ and
$z(\bfx,\thetavec,t)-z(\tilde{\bfx},\thetavec,t)\in L_\infty$ then
$\dpsi\in
 L_\infty, \ \ z(\tilde{\bfx},\thetavec,t)-z(\tilde{\bfx},\hat{\thetavec},t)\in
 L_\infty, \ \ \lim_{t\rightarrow \infty}\psi(\bfx(t),t)= 0$.
\end{thm}

Theorem \ref{theorem_embedd} states not only existence of the
adaptive control algorithms but also provides us with exact
equations for the adaptive control function. These equations are
given by (\ref{control_embedd}), (\ref{fin_forms_embedd}), which
guarantee P8), and (\ref{control_embedd}),
(\ref{fin_forms_embedd_1}) or (\ref{fin_forms_embedd_2}) ensuring
P9) for $\pd\alphavec(\tilde\bfx)/\pd \bfx_2'\equiv 0$ or
$\psi(\tilde\bfx)=\psi(\bfx)$ respectively.

Although Theorem \ref{theorem_embedd} guarantees reaching of the
control goal and ascertains performance improvement (property
P9)), it does not ensure the same properties of adaptive control
algorithms as Theorem \ref{theorem_realization1} does. On the
other hand, the ability to deal with nonconvex parameterized
systems is preserved, except for cases that do not satisfy
Assumption \ref{alpha_upper_bound}. The drawbacks of this narrower
class of nonlinearly parameterized functions in the plant
right-hand side and a slight degradation in performance are
compensated by relaxing the requirement (\ref{realizability_eq})
of Assumption \ref{Realizability_1}. Notice also that the
difference in guaranteed  performance reflected in  P8) and P9)
has the consequence that the dimensions of vectors $\bfh_{\xi}$
are likely to be different in the both cases.
% followed by the
%difference in dimensions of vectors $\bfh_{\xi}$ for the both
%cases.
Indeed, to ensure equality $\pd \alphavec(\tilde\bfx,t)/\pd
\bfx_2'\equiv 0$ for arbitrary smooth function $\alphavec(\cdot)$,
we must replace the whole vector $\bfx_2$ by $\bfh_{\xi}(\xivec)$.
Therefore, in principle, embedding of the original system dynamics
into one of a higher order is desired if improved performance and
extended applicability are required.

%\begin{rem}\normalfont
Theorem \ref{theorem_embedd} offers a possible way to facilitate
the search for function $\Psi(\bfx,t)$ satisfying partial
differential equation (\ref{realizability_eq}) as defined in
Assumption \ref{Realizability_1}.  We replace the problem by one
of searching  for the embedding (\ref{plant_ins_embed}) which
satisfies Assumption \ref{assume_embedd} and
\ref{assume_L_infty_embed} or \ref{assume_L_2_embed}. The main
obstacle, finding a solution to equation (\ref{realizability_eq}),
is replaced with problem (\ref{realizability_eq2}), the
complexity\footnote{Here reduced complexity means that the number
of equations in the system is reduced.} of which should be
reduced, as $\dim{\bfx}_2'<\dim{\bfx}_2$ if embedding into
the higher-order dynamics is used.  %{\bf It's not obvious why the
%complexity will be lower.  Assuming a lucky higher-order
%embedding?}

Indeed, according to Assumption \ref{assume_embedd} and notations
introduced above, the dynamics of the extended system can be
described as
\begin{eqnarray}\label{extended_embed}
\dot{\bfx}_1&=&\bff_1(\bfx)+\bfg_1(\bfx)u,  \ \ \ \dot{\bfx}_{2'}=\bff_{2'}(\bfx)+\nuvec'(\bfx,\thetavec)+\bfg_{2'}(\bfx)u\nonumber \\
\dot{\bfh}_{\xi}&=&\frac{\pd \bfh_{\xi}}{\pd
\xivec}\bff_{\xi}(\bfx,\xivec,t), \ \ \ \
\dot{\bfx}_{2''}=\bff_{2''}(\bfx)+\nuvec''(\bfx,\thetavec)+\bfg_{2''}(\bfx)u,
\end{eqnarray}
where vector $\bfx_1\oplus\bfh_{\xi}$ stands for the
uncertainty-independent partition in the extended state space, and
vector $\bfx_2'$ is chosen to satisfy equation
(\ref{realizability_eq2}). Observe that function
$z(\tilde{\bfx},\thetavec,t)$ is independent of $\bfx_2''$ and
$\dim{\bfh_{\xi}=\dim{\bfx}_2''}$. Then for any $\bfh_{\xi}$:
$\dim{\bfh_{\xi}}>0$, we can conclude that
$\dim{\bfx}_2'<\dim{\bfx_2}=\dim{\bfx_2'\oplus\bfx_2''}$.

Notice also that, by the appropriate choice of the dimensions of
vectors $\xivec$ and $\bfh_{\xi}$
($\dim{\bfh_{\xi}=\dim{\bfx}_2''}$) in (\ref{embed_add}), the
dimension of vector $\bfx_2'$ can be reduced  to unity.
Alternatively, we may try to annihilate the partial derivative
$\frac{\pd \alphavec(\tilde{\bfx},t)}{\pd \bfx_2'}$  in
(\ref{realizability_eq2}). Hence, eventually  either Corollary
\ref{single_dim_par}  or  Corollary
\ref{independent_partition_cor} conditions will
 be satisfied for the extended system
(\ref{extended_embed}). This, in turn, implies that we can replace
assumption (\ref{realizability_eq2}) by a weaker requirement, such
as integrability of the function with respect to a single scalar
argument.

After obtaining computable function $\Psi(\tilde{\bfx},t)$, the
remaining problem is that we should be able to find an
extension (\ref{embed_add}) that guarantees properties
(\ref{L_2_embed}) and (\ref{boundedness_xi}) for the given
partition $\tilde{\bfx}=\bfx_1\oplus\bfx_2'\oplus\bfh_{\xi}$. If
such an extension exists, then Assumption \ref{assume_embedd} is
automatically satisfied, and adaptive control algorithms follow
immediately from Theorem \ref{theorem_embedd}.

Finding extension (\ref{embed_add}) that ensures boundedness  (and
square integrability) of the differences
$z(\bfx,\thetavec,t)-z(\tilde\bfx,\thetavec,t)$,
$(\psi(\bfx,t)-\psi(\tilde{\bfx},t))({\pd
\alphavec(\tilde\bfx,t)}/{\pd \bfx_2'})\nuvec'(\bfx,\thetavec)$ is
not an easy problem -- taking into account that partition
$\bfx_2''$ is also uncertainty-dependent.  It is possible to solve
it using specially designed {\it adaptive} or {\it high-gain}
auxiliary subsystems that track the reference signals $\bfx_2''$
with the desired performance: $z(\bfx,
\thetavec,t)-z(\tilde{\bfx},\thetavec,t),
(\psi(\bfx,t)-\psi(\tilde{\bfx},t))\frac{\pd
\alphavec(\tilde\bfx,t)}{\pd \bfx_2'}\nuvec'(\bfx,\thetavec)\in
L_2\cap L_\infty$. If, for example, partition $\bfx_2''$ is
linearly parameterized (i.e.,
$\nuvec''(\bfx,\thetavec)=\etavec''(\bfx)\thetavec)$, functions
$z(\bfx,\thetavec,t), \ \psi(\bfx,t)\frac{\pd
\alphavec(\tilde\bfx,t)}{\pd \bfx_2'}\nuvec'(\bfx,\thetavec)$ are
locally Lipshitz in $\bfx_2''$ and for any $\thetavec\in
\Omega_{\theta}$ the following inequalities hold:
$|z(\bfx,\thetavec,t)-z(\tilde\bfx,\thetavec,t)|\leq
\lambda_1(\bfx,\xivec,t)\|\bfx_2''-\bfh_\xi(\xivec)\|$,
$\|(\psi(\bfx,t)-\psi(\tilde{\bfx},t))\frac{\pd
\alphavec(\tilde\bfx,t)}{\pd \bfx_2'}\nuvec'(\bfx,\thetavec)\|\leq
\lambda_2(\bfx,\xivec,t)\|\bfx_2''-\bfh_\xi(\xivec)\|$ then the
suitable extension is defined by the following system:
\begin{eqnarray}\label{lin_extension}
\dot{\xivec}_1&=&\bff_2''(\bfx)+\etavec''(\bfx)\xivec_2+\bar{\lambda}(\bfx,\xivec,t)(\bfx_2''-\xivec_1)+\bfg_2''(\bfx)u\nonumber\\
\dot{\xivec}_2&=&\Gamma_1(\bfx_2''-\xivec_1)^{T}\etavec''(\bfx), \
\Gamma_1>0, \  \ \bfh_{\xi}(\xivec)=\xivec_1,
\end{eqnarray}
where $\xivec=\xivec_1\oplus\xivec_2$ and
$\bar{\lambda}(\bfx,\xivec,t)=\lambda_1^2(\bfx,\xivec,t)+\lambda_2^2(\bfx,\xivec,t)$.
To show this, it is sufficient to consider the following
Lyapunov's candidate:
$V(\bfx,\xivec)=0.5\|(\bfx_2''-\xivec_1)\|^2+0.5\|\thetavec-\xivec_2\|^{2}_{\Gamma^{-1}}$
and observe that $\dot{V}\leq - \bar{\lambda}(\bfx,\xivec,t)
\|\bfx_2''-\xivec_1\|^2\leq -
(z(\bfx,\thetavec,t)-z(\tilde{\bfx},\thetavec,t))^2\leq 0$.

For a class of nonlinear systems with low-triangular structure
\begin{eqnarray}\label{cascade__full_2}
\dot{x}_i&=&f_i(x_1,\dots,x_i,\thetavec_i)+x_{i+1}, \ \ i=1,\dots n-1,\nonumber\\
\dot{x}_n&=&f_n(x_1,\dots,x_n,\thetavec_n)+u + \varepsilon(t), \
\varepsilon(t)\in L_2,  \ \thetavec_i\in \Omega_{\theta}
\end{eqnarray}
the suitable extension is guaranteed by Lemma
\ref{lemma_embed_full_cascade} in Appendix 1. Then, combining the
results formulated in Theorem \ref{theorem_embedd} and Lemma
\ref{lemma_embed_full_cascade}, it is possible to show that  our
approach can be extended to a broad class of  systems like those
given by equations (\ref{cascade__full_2}). Let functions
$f_i(\cdot)$ in (\ref{cascade__full_2}) satisfy the following
assumption
\begin{assume}\label{assume_f_i} Let there exist
smooth functions $\bar{D}_i(\cdot): R^{i}\times R^{i}\times
R\rightarrow R$ such that for any $\thetavec_i\in \Omega_{\theta}$
the following holds:
%\begin{eqnarray}\label{growth_bound_f_i}
%& &
$(f_i(x_1,\dots,x_i,\thetavec_i)-f_i(x_1',\dots,x_i',\thetavec_i))^2\leq
\bar{D}_i^2(\bfx_i,\bfx_i')\|\bfx_i-\bfx_i'\|^2$,
$\bfx_i=(x_1,\dots,x_i)^{T}$, $\bfx_i'=(x_1',\dots,x_i')$.
%\end{eqnarray}
\end{assume}
It is clear that Assumption \ref{assume_f_i} holds for those
functions $f_i(x_1,\dots,x_i,\thetavec_i)$ that are, for example,
Lipshitz in $\bfx$. The results for low-triangular systems
(\ref{cascade__full_2}) are formulated in the next theorem
%On the other hand if $\Omega_{\theta}$ is compact then
%(\ref{growth_bound_f_i}) holds as well. To show this it is enough
%to apply Hadamar's lemma to the difference:
%\[
%f_i(x_1,\dots,x_i,\thetavec_i)-f_i(x_1',\dots,x_i',\thetavec_i)=F(\bfx,\bfx',\thetavec_i)(\bfx_i-\bfx_i')
%\]
%and then Lemma 2.1 from \cite{Lin_smooth} to the
\begin{thm}[Finite Forms for Low-Triangular Systems]\label{theorem_full_cascade} Let system
(\ref{cascade__full_2}) and goal function $\psi(x_1)=0$ be given,
and there exist functions $\alphavec_i(x_1,\dots,x_i)$ such that
Assumptions \ref{alpha}, \ref{alpha_upper_bound} hold for the
functions $f_i(x_1,\dots,x_i,\thetavec_i)$ in
(\ref{cascade__full_2}) respectively. Furthermore, let
$f_i(x_1,\dots,x_i,\thetavec_i)$ satisfy Assumption
\ref{assume_f_i},  $\alphavec_i(x_1,\dots,x_i)$,
$f_i(x_1,\dots,x_i,\thetavec_i)$, $i=1,\dots,n$, $\psi_1(x_1)$ be
smooth and the following condition holds: $\psi(x_1)\in L_\infty
\Rightarrow x_1\in L_\infty$.

%the
%following integrals exist
%\begin{eqnarray}\label{Psi_i_cascade_full}
%\Psi_i(x_i,t)=\int x_i\frac{\pd
%\alphavec_i(x_1,\dots,x_i,\omegavec_i(t))}{\pd x_i} dx_i
%\end{eqnarray}
Then there exist an auxiliary system
\begin{eqnarray}\label{embed_track_sys1}
\dot{\xivec}&=&\bff_{\xi}(\bfx,\xivec,\nuvec), \ \xivec_0\in
R^{n}, \ \ \dot{\nuvec}=\bff_{\nu}(\bfx,\xivec,\nuvec), \ \nuvec_0
\in R^m,
\end{eqnarray}
as well as smooth  functions  $\psi_i(x_i,t)$, $i=1,\dots,n$,
$\hat{\thetavec}_{P}(\bfx,\xivec)$, control
$u(\bfx,\hat{\thetavec},\xivec,\nuvec)$, and adaptation algorithm
\begin{eqnarray}
\hat{\thetavec}(\bfx,\xivec,\hat{\thetavec}_{I})&=&\gamma
(\hat{\thetavec}_{P}(\bfx,\xivec)+\hat{\thetavec}_{I}), \
\gamma>0,  \  \
\dot{\hat{\thetavec}}_{I}=\bff_{\hat\theta}(\bfx,\hat{\thetavec},\xivec,\nuvec),\nonumber
\end{eqnarray}
such that

1) $\psi_i(x_i,t),\psi \in L_2\cap L_{\infty}$, $\dpsi,\dpsi_i\in
L_2$, $i=1,\dots,n$

2) $\hat{\thetavec}\in L_\infty$ and
$u(\bfx,\hat{\thetavec},\xivec,\nuvec)-u(\bfx,\thetavec_n,\xivec,\nuvec)\in
L_2$

3) $\bfx,\xivec,\nuvec \in L_\infty$

4) if $\varepsilon(t)\in L_\infty$ then $\dpsi,\dpsi_i\in
L_\infty$,  and  $\lim_{t\rightarrow \infty}\psi(x_1(t))=0, \ \
\lim_{t\rightarrow\infty}\psi_i(x_i(t),t)=0$, $i=1,\dots,n$.

%Moreover, if
%
%\[
%f_i(0,\dots,0,\thetavec_i)=0,  \
%\lim_{t\rightarrow\infty}\psi(x_1(t))=0\Rightarrow
%\lim_{t\rightarrow\infty}x_1(t)=0, \  \dot{\varepsilon}(t)\in
%L_\infty
%\]
%5) then $x_i\rightarrow 0$ as $t\rightarrow\infty$, $i=1,\dots,n$,
%$u(\bfx(t),\hat{\thetavec}(t),\xivec(t),\nuvec(t))\rightarrow 0$
%as $t\rightarrow \infty$.
\end{thm}

%The theorem shares some of the statements of Theorem
%{\ref{theorem_cascade_finite_forms}}, Part I, as formulated for
%the class of systems given by  (\ref{cascade_2}), Part I.
Theorem \ref{theorem_full_cascade} extends the applicability of
algorithms in finite form to systems described by equation
(\ref{cascade__full_2}). Relying entirely on Lemma
\ref{lemma_embed_full_cascade} (Appendix 1) and  Theorem
\ref{theorem_embedd}, Theorem \ref{theorem_full_cascade} allows us
to design adaptive control algorithms for cascades with nonlinear
parameterization without the need for damping nonlinearities.
However, performance is weaker. For instance, decrease
(non-increase) of the term
$\|\thetavec-\hat{\thetavec}(t)\|^{2}_{\Gamma^{-1}}$ is not
guaranteed in this case.
%{\bf Is it our case?}
Nevertheless, adaptive control algorithms in finite form, in
addition to their ability to deal with nonlinear parameterization,
still guarantee certain improvements in performance. For instance,
square integrability of the control effort due to adaptation
(statement 2) of the theorem)
%{\bf Why do we require $L_2\cap L_\infty$?  I thought $\in L_2$ also implies $\in L_\infty$.}
and $\psi_i(x_i,t),\dpsi_i\in L_2\cap L_\infty$ are ensured.
%Moreover, as can be derived from Lemma \ref{L_2_theorem_eps}, the
%upper bounds of $L_2$ norms for functions $\varphi(\psi)$ and
%$\dpsi$ have the same value for both $\varphi(\psi)$ and $\dpsi$.
In the next section we illustrate our method with the examples.

\section{Examples}

Let us consider the following system:
\begin{eqnarray}\label{example_plant2}
\dot{x}_1&=&x_1^2\theta_0+x_2;  \ \ \dot{x}_2=x_1\theta_1 + x_2
\theta_2 + u,
\end{eqnarray}
where parameters $\theta_0$,$\theta_1$ and $\theta_2$ are assumed
to be unknown. The control goal is to steer the system towards the
following manifold: $x_1-1=0$. To design adaptive algorithms in
finite form for system (\ref{example_plant2}), we follow the steps
of Theorem \ref{theorem_full_cascade} proof:

1) {\it Intermediate control design}. Derive control function
$u_1(x_1,\hat{\theta}_0)$ such that for the reduced system
\begin{eqnarray}
\dot{x}_1=x_1^2\theta_0 + u_1(x_1,{\hat\theta}_0) +
\varepsilon_1(t), \ \varepsilon_1(t)\in L_2; \ \
{\hat\theta}_0={\hat\theta}_{0,P}(x_1)+{\hat\theta}_{0,I}(t)\nonumber
\end{eqnarray}
reaching of the control goal is guaranteed:
$\psi(x_1(t))=x_1(t)-1\rightarrow 0$ as $t\rightarrow\infty$.
Moreover, function
$u_1(x_1,{\hat\theta}_0(x_1,{\hat\theta}_{0,I}))$ should ensure
that $\psi,\dpsi\in L_2$.

2) {\it Embedding}. Extend the system dynamics (or embed it into)
with auxiliary system
\begin{eqnarray}\label{ex2_emb1}
\dot{\xi}=f_{\xi}(\bfx,\xi,\nu); \ \dot{\nu}=f_{\nu}(\bfx,\xi,\nu)
\end{eqnarray}
in order to guarantee that
\begin{eqnarray}\label{ex2_L_21}
u(x_1,{\hat\theta}_0(x_1,{\hat\theta}_{0,I}))-u(\xi,{\hat\theta}_0(\xi,{\hat\theta}_{0,I}))\in
L_2, \ \ x_1-\xi\in L_2
\end{eqnarray}

3) {\it Control function design}. Introduce new goal function
$\psi_2(x_2,t)= x_2 -
u_1(\xi,{\hat\theta}_0(\xi,{\hat\theta}_{0,I}))$ and derive
control function $u(x_1,x_2,\xi,t)$ such that $\dpsi_2\in L_2$,
$\psi_2\in L_2\cap L_\infty$. The last automatically implies that
$\dot{x}_1=x_1^2\theta_0 +
x_2=x_1^2\theta_0+u_1(x_1,{\hat\theta}_0(x_1,{\hat\theta}_{0,I}))+\mu(t)$,
where
$\mu(t)=x_2-u_1(x_1,{\hat\theta}_0(x_1,{\hat\theta}_{0,I}))=(x_2-u_1(\xi,{\hat\theta}_0(\xi,{\hat\theta}_{0,I})))+(u_1(\xi,{\hat\theta}_0(\xi,{\hat\theta}_{0,I}))-u_1(x_1,{\hat\theta}_0(x_1,{\hat\theta}_{0,I})))\in
L_2$. Therefore, according to the choice of function
$u_1(x_1,{\hat\theta}_0(x_1,{\hat\theta}_{0,I}))$, control
$u(x_1,x_2,\xi,t)$ guarantees that $\psi(x_1(t))\rightarrow 0$ as
$t\rightarrow\infty$, $\psi,\dpsi\in L_2$.

We begin by determining the function
$u(x_1,{\hat\theta}_0(x_1,{\hat\theta}_{0,I}))$. Let
$u_1(x_1,{\hat\theta}_0)=-(x_1-1)-\hat{\theta}_0 x_1^2$, where
$\hat{\theta}_0$ satisfies the following differential equation:
\begin{eqnarray}\label{ex2_alg1}
\dot{\hat{\theta}}_0=(x_1-1+\dot{x}_1)x_1^2
\end{eqnarray}
It follows from Lemma \ref{L_2_theorem_eps} that control function
$u_1(x_1,{\hat\theta}_0)$ with algorithm (\ref{ex2_alg1})
guarantee that $\psi,\dpsi\in L_2$, $\psi(x_1(t))\rightarrow 0$ as
$t\rightarrow\infty$.
% On the other hand,
According to Theorem \ref{theorem_realization1}, finite form
realization of (\ref{ex2_alg1}) can be given as follows:
$\hat{\theta}_{0}(x_1,\hat{\theta}_{0,I}(t))=1/3 x_1^3
+\hat{\theta}_{0,I}(t); \ \
\dot{\hat{\theta}}_{0,I}=(x_1-1)x_1^2$. Substituting this into
$u_1(x_1,{\hat\theta}_0)$ we get the following expression for
$u_1(\cdot)$:
\begin{eqnarray}\label{ex2_contr2}
u_1(x_1,{\hat\theta}_0(x_1,{\hat\theta}_{0,I}))=-(x_1-1)-{1}/{3}x_1^{5}-x_1^{2}{\hat\theta}_{0,I}(t);
\ \ \dot{\hat{\theta}}_{0,I}=\psi(x_1)\alpha_1(x_1)=(x_1-1)x_1^2.
\end{eqnarray}
Thus step 1 is completed.

Let us design system (\ref{ex2_emb1}) which guarantees that
(\ref{ex2_L_21}) holds for function (\ref{ex2_contr2}). First
consider the difference:
\begin{eqnarray}\label{ex2_emb2}
u(x_1,{\hat\theta}_0(x_1,{\hat\theta}_{0,I}))-u(\xi,{\hat\theta}_0(\xi,{\hat\theta}_{0,I}))=-(x_1-\xi)(1+(x_1+\xi)\hat{\theta}_{0,I}+
{1}/{3}(x_1^4+x_1^3\xi+x_1^2\xi^2+x_1\xi^{3}+\xi^{4}))
\end{eqnarray}
and denote
$F(x_1,\xi,\hat{\theta}_{I,0})=(1+(x_1+\xi)\hat{\theta}_{0,I}+
\frac{1}{3}(x_1^4+x_1^3\xi+x_1^2\xi^2+x_1\xi^{3}+\xi^{4}))$. It
follows from Lemma \ref{lemma_embed_full_cascade} that there
exists system (\ref{ex2_emb1}) such that condition
(\ref{ex2_L_21}) holds. In fact, this system can be given by the
following equation
\begin{eqnarray}\label{ex2_emb3}
\dot{\xi}&=&(x_1-\xi)(F^2(x_1,\xi,\hat{\theta}_{0,I})+1)+x_1^2{\hat\theta}_{\xi}+x_2
\end{eqnarray}
where ${\hat\theta}_{\xi}$ satisfies the following differential
equation
%\begin{eqnarray}\label{ex2_alg3}
$\dot{\hat{\theta}}_{\xi}=(x_1-\xi+\dot{x}_1-\dot{\xi})x_1^2$.
%\end{eqnarray}
Finite form  realization  of this algorithm\footnote{Introduction of algorithms
(\ref{ex2_alg4}) is not necessary here because the original system is
linearly parameterized, and condition (\ref{ex2_L_21}) can be
satisfied even with conventional (gradient) adaptation schemes.
Nevertheless, we would like to keep the consistency of our current
calculations with those steps made in the proof of Theorem
\ref{theorem_full_cascade} in order to illustrate what would
happen if the right hand sides are nonlinearly parameterized.}
follows from Theorem \ref{theorem_realization1}, and it can be
written as:
\begin{eqnarray}\label{ex2_alg4}
{\hat\theta}_{\xi}&=&{1}/{3}x_1^3+{\hat\theta}_{\xi,I};  \ \
 \dot{\hat{\theta}}_{\xi,I}= (x_1-\xi)x_1^2-x_1^{2}((x_1-\xi)(F^2(x_1,\xi,\hat{\theta}_{0,I})+1)+x_1^2{\hat\theta}_{\xi}+x_2)
\end{eqnarray}
Taking into account (\ref{ex2_alg4}) and (\ref{ex2_emb3}) system
(\ref{ex2_emb1}) which ensures (\ref{ex2_L_21}) can be represented
as follows
\begin{eqnarray}\label{ex2_emb4}
\dot{\xi}&=&(x_1-\xi)(F^2(x_1,\xi,\hat{\theta}_{0,I})+1)+\frac{1}{3}x_1^5+\hat{\theta}_{\xi,I}(t)x_1^2+x_2\nonumber\\
\dot{\hat{\theta}}_{\xi,I}&=&(x_1-\xi)x_1^2-x_1^{2}((x_1-\xi)(F^2(x_1,\xi,\hat{\theta}_{0,I})+1)+\frac{1}{3}x_1^5+\hat{\theta}_{\xi,I}(t)x_1^2+x_2).
\end{eqnarray}
Therefore, step 2 is completed as well. To conclude the controller
design let us consider new target manifold
$x_2-u_1(\xi,{\hat\theta}_0(\xi,{\hat\theta}_{0,I}))=0$ and goal
function
$\psi_2(x_2,t)=x_2-u_1(\xi,{\hat\theta}_0(\xi,{\hat\theta}_{0,I}))=x_2+\xi-1+\frac{1}{3}\xi^5+\hat{\theta}_{0,I}\xi^2$.
Let us write function $\psi_2(\cdot)$ derivative with respect to
time $t$:
\begin{eqnarray}
\dpsi_2&=&\dot{x}_2-\frac{\pd
u_1(\xi,{\hat\theta}_0(\xi,{\hat\theta}_{0,I}))}{\pd \xi}
\dot{\xi} - \frac{\pd
u_1(\xi,{\hat\theta}_0(\xi,{\hat\theta}_{0,I}))}{\pd
\hat{\theta}_{0,I}}\dot{\hat\theta}_{0,I}
= x_1\theta_1+x_2\theta_2 + u + \xi^2(x_1-1)x_1^2+\nonumber \\
& & (1+\frac{5}{3}\xi^4+2\xi
\hat{\theta}_{0,I})((x_1-\xi)(F^2(x_1,\xi,\hat{\theta}_{0,I})+1)+\frac{1}{3}x_1^5+\hat{\theta}_{\xi,I}(t)x_1^2+x_2)\nonumber
\end{eqnarray}
Therefore, control function
\begin{eqnarray}\label{ex2_contr3}
u&=&-\xi \hat{\theta}_1 - x_2 \hat{\theta}_2 - \xi^2(x_1-1)x_1^2 -
(x_2+\xi-1+\frac{1}{3}\xi^5+\hat{\theta}_{0,I}\xi^2) -\nonumber\\
& &  (1+\frac{5}{3}\xi^4+2\xi
\hat{\theta}_{0,I})((x_1-\xi)(F^2(x_1,\xi,\hat{\theta}_{0,I})+1)+\frac{1}{3}x_1^5+\hat{\theta}_{\xi,I}(t)x_1^2+x_2)
\end{eqnarray}
results in the following error model:
%\begin{eqnarray}%\label{ex2_err1}
$\dpsi_2=-\psi_2(x_2,t)+x_1\theta_1 + x_2 \theta_2 -
x_1\hat{\theta}_1 - x_2 \hat{\theta}_2$.
%\end{eqnarray}
Taking into account condition (\ref{ex2_L_21}), we can rewrite
derivative $\dpsi_2$ as
%\begin{eqnarray}%\label{ex2_err2}
$\dpsi_2=-\psi_2(x_2,t)+\xi \theta_1 + x_2 \theta_2 - \xi
\hat{\theta}_1 - x_2 \hat{\theta}_2 + \varepsilon(t)$,
%\end{eqnarray}
where $\varepsilon(t)=(x_1-\xi)\theta_1\in L_2$. It follows from
Lemma \ref{L_2_theorem_eps} that adaptation algorithm
\begin{eqnarray}\label{ex2_alg5}
& &{\dot{\hat\theta}}_1=(\psi_2(x_2,t)+\dpsi_2)\alpha_1(\xi); \ \
{\dot{\hat\theta}}_2=(\psi_2(x_2,t)+\dpsi_2)\alpha_2(x_2), \
\alpha_1(\xi)=\xi, \ \alpha_2(x_2)=x_2 \nonumber
\end{eqnarray}
guarantees $\psi_2\in L_2\cap L_\infty$ and $\dpsi_2\in L_2$.
Realization of algorithms (\ref{ex2_alg5}) can be obtained from
Theorem \ref{theorem_realization1}:
\begin{eqnarray}\label{ex2_alg6}
\hat{\theta}_1(x_2,\xi,\hat{\theta}_{0,I},t)&=&(x_2+\xi-1+\frac{1}{3}\xi^5+\hat{\theta}_{0,I}\xi^2)\xi+\hat{\theta}_{1,I}(t);
\ \
\dot{\hat{\theta}}_{1,I}=(x_2+\xi-1+\frac{1}{3}\xi^5+\hat{\theta}_{0,I}\xi^2)(\xi-\dot{\xi})\nonumber
\\
\hat{\theta}_2(x_2,\xi,\hat{\theta}_{0,I},t)&=&\frac{x_2^2}{2}+\hat{\theta}_{2,I}(t);
\ \
\dot{\hat{\theta}}_{2,I}=(x_2+\xi-1+\frac{1}{3}\xi^5+\hat{\theta}_{0,I}\xi^2)x_2+\frac{\pd
\Psi_2}{\pd \xi}\dot{\xi}+\frac{\pd \Psi_2}{\pd
\hat{\theta}_{0,I}}\dot{\hat{\theta}}_{0,I},
\end{eqnarray}
where $\Psi_2(x_2,\xi,\hat{\theta}_{0,I})=\int
\psi_2(x_2,t)\frac{\pd \alpha_2(x_2)}{\pd x_2}d x_2 =
\frac{x_2^2}{2}+(\xi-1+\frac{1}{3}\xi^5+\hat{\theta}_{0,I}\xi^2)x_2$.

We would also like to compare performance of the proposed
adaptation scheme with adaptive backstepping control algorithms.
Adaptive backstepping design for system (\ref{example_plant2})
according to \cite{Kanellakopoulos91} results in
%following
control algorithm:
\begin{eqnarray}\label{example2_backstepping_2}
u_1&=&-2x_2-(x_1-1)-\hat{\theta}_3
x_1^2-x_1^4(x_1-1)-2{\hat\theta}_3 x_1 x_2-(x_1^2+2{\hat\theta}_3
x_1^3)\hat\theta-x_1{\hat\theta}_1-x_2{\hat\theta}_2\nonumber \\
\dot{\hat\theta}&=&(x_2+x_1-1+{\hat\theta}_3 x_1^2)
x_1^2(1+2{\hat\theta}_3 x_1); \ \
\dot{\hat\theta}_1=(x_2+x_1-1+{\hat\theta}_3 x_1^2)x_1\nonumber
\\
\dot{\hat\theta}_2&=&(x_2+x_1-1+{\hat\theta}_3x_1^2)x_2; \ \
\dot{\hat\theta}_3=(x_1-1)x_1^2
\end{eqnarray}
Adaptive backstepping with tuning functions \cite{Krstic92}
results in
\begin{eqnarray}\label{example2_backstepping_22}
u_1&=&-(x_2+x_1-1+x_1^2\hat\theta)-(x_1-1)-(1+2x_1\hat\theta)(x_2+\hat\theta
x_1^2)-x_1^2 \tau-x_1{\hat\theta}_1-x_2{\hat\theta}_2\nonumber \\
\dot{\hat\theta}&=&\tau; \ \ \tau = (x_1-1) x_1^2 + (x_2+x_1-1+x_1^2{\hat\theta})x_1^2(1+2 x_1\hat\theta) \nonumber \\
\dot{\hat\theta}_1&=& (x_2+x_1-1+x_1^2\hat\theta)x_1; \ \
\dot{\hat\theta}_2=(x_2+x_1-1+x_1^2\hat\theta)x_2
\end{eqnarray}
We simulated the adaptive system dynamics for the following set of
parameters and initial conditions $x_1(0)=2, \ x_2(0)=0.2,  \
\theta=1, \ \theta_3(0)=\theta(0)=3, \
\theta_{1}(0)=\theta_2(0)=-2, \ \xi_2(0)=0, \ \xi_1(0)=0, \ k=10$.
 Initial conditions for
${\hat\theta}_{1,I}(0),{\hat\theta}_{2,I}(0)$ and
${\hat\theta}_{3,I}(0)$ where chosen to satisfy
$\hat\theta_1(0)=\hat\theta_2(0)=-2, \ \hat\theta_3(0)=3$.
Parameters $\theta_1,\ \theta_2$: $\theta_1=1$, $\theta_2=0.5$. %{\bf Let's try other initial and target conditions
%just to make sure that ours is always better.}
As an additional measure of performance, we introduced the
variable $\Delta\hat\theta$ which indicates the distance in the
controller parameter space between the estimates and real values
of the parameters. Simulation results are presented in Figure 1.
In Figure 1 thick solid lines show the system dynamics with
algorithm (\ref{ex2_contr3}),(\ref{ex2_alg6}), dotted lines show
the system behavior with algorithm (\ref{example2_backstepping_2}),
and dash-dotted lines correspond to algorithm
(\ref{example2_backstepping_22}). It turns out that
system (\ref{example_plant2}) with algorithm
(\ref{example2_backstepping_22}) also reaches the goal manifold,
but after around $400$ seconds of modeling time. We can see again
that transient performance of the adaptive algorithms proposed in
the paper is better than that of conventional algorithms. In
addition, we calculated the integral $I=\int_0^{T}
u_1^2(\tau)d\tau, \ T=500$, for every controller along the system
solutions. The values of the functional $I$ indicate how much energy
is spent to achieve the control goal. For control function
(\ref{ex2_contr3}),(\ref{ex2_alg6}) $I=627.10$, for adaptive
backstepping controller (\ref{example2_backstepping_2}) $I=13
329.28$, for controller (\ref{example2_backstepping_22}) $I=263
872.58$.
\begin{figure}
\begin{center}
\includegraphics[width=300pt]{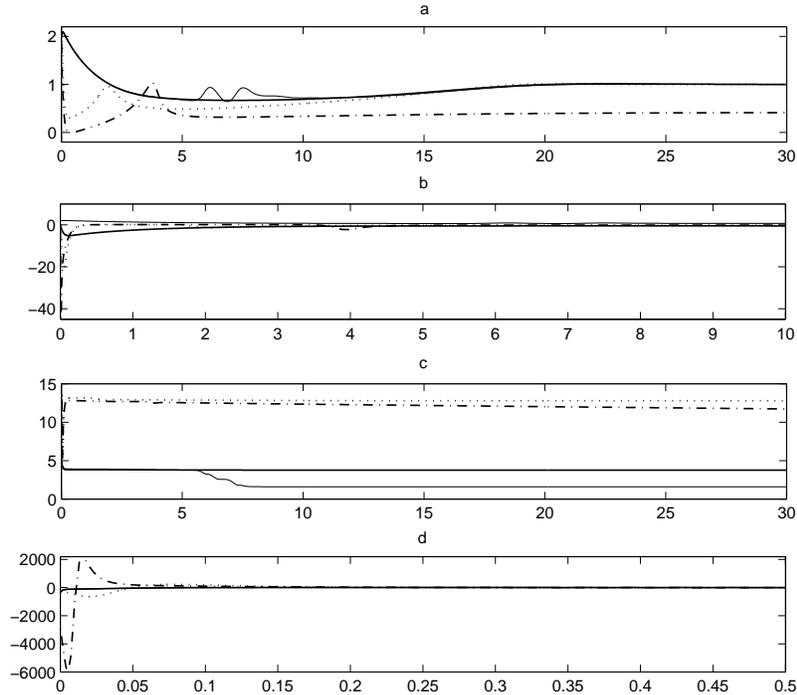}
\end{center}
\begin{center}
\caption{Plots of system (\ref{example_plant2}),
(\ref{example_plant3}) trajectories with control functions
(\ref{ex2_contr3}),(\ref{ex2_alg6}) (thick solid lines),
(\ref{example2_backstepping_2}) (dotted line), (\ref{example2_backstepping_22}) (dash-dotted line),
and (\ref{example_control3}) (thin solid line). Plot $a$ -- $x_1$ as a function of
time, $b$ -- $x_2$  as a function of time, $c$ --
$\Delta\hat{\theta}$ as a function of time, $d$ -- $u$ as a
function of time.}
\end{center}
\end{figure}
To illustrate the ability of our algorithms to deal with nonlinear
parameterization, we change (\ref{example_plant2}) to
\begin{eqnarray}\label{example_plant3}
\dot{x}_1&=&x_1^2\theta+x_2; \  \dot{x}_2= 5\tanh(x_1\theta_1 +
x_2 \theta_2) + u_2
\end{eqnarray}
Nonlinearity $\tanh(x_1\theta_1 + x_2 \theta_2)$ satisfies
Assumption \ref{alpha} with respect to function
$\alphavec(\bfx)=(x_1,x_2)^T$ and, in addition, Assumption
\ref{alpha_upper_bound} is also satisfied for any bounded $x_1$
and $x_2$. Then according to Theorem \ref{theorem_full_cascade},
control function
\begin{eqnarray}\label{example_control3}
u&=&-5\tanh(\xi \hat{\theta}_1 + x_2 \hat{\theta}_2) -
\xi^2(x_1-1)x_1^2 -
(x_2+\xi-1+\frac{1}{3}\xi^5+\hat{\theta}_{0,I}\xi^2) -\nonumber\\
& &  (1+\frac{5}{3}\xi^4+2\xi
\hat{\theta}_{0,I})((x_1-\xi)(F^2(x_1,\xi,\hat{\theta}_{0,I})+1)+\frac{1}{3}x_1^5+\hat{\theta}_{\xi,I}(t)x_1^2+x_2)
\end{eqnarray}
along with  (\ref{ex2_alg6}) guarantees that $\psi_1,
\psi_2,\dpsi_1,\dpsi_2\in L_2\cap L_\infty$. The simulation
results of system (\ref{example_plant3}) with control algorithm
(\ref{example_control3}), (\ref{ex2_alg6}) are given in Figure 1
(thin solid lines). The value of functional $I$ for this case is
$3186.83$.

\section{Conclusions}

The method proposed in this paper suggests a new methodology to
design adaptive control algorithms. Our method to design of the
adaptation schemes is consistent with recent trends in adaptive
control, for instance, \cite{Krstic94}, where nonlinear
controllers are proposed to adaptively stabilize
 linear plants. Indeed, when
derived for linear systems algorithms in finite form
will also result in nonlinearities in
the controller. These nonlinearities are to be introduced, in
particular, to improve the performance of the adaptive system.  In
contrast to \cite{Krstic94}, we show  not only that the $L_2$ and
$L_\infty$ norm bounds are computable for the state vector, but
also that properties P1)--P7) are ensured. The method, however, is
different from conventional approaches, as it is not restricted
by realizability issues. While in conventional
parametric adaptive control the realizability
%(or
%uncertainty-independence)
of adaptation schemes in differential form determines the
properties of the resulting systems (including poor performance
and restricted applicability), in our method we first determine
the desired properties of the controller (Theorem \ref{L_2
dpsi_theorem}, Proposition \ref{Exponential_Convergence} and
Lemmas \ref{lemma_dist_feedback_alg} -- \ref{L_2_theorem_eps}) and
only then deal with the realizability problem. In order to realize
the adaptive algorithms in finite form explicitly, i.e. without
extension of the system state space, special restrictions
formulated in Assumption \ref{Realizability_1} are to be satisfied
(Theorem \ref{theorem_realization1} and Corollaries
\ref{single_dim_par}, \ref{independent_partition_cor}).

To realize the adaptation algorithms that do not satisfy the
explicit realizability conditions formulated in Assumption
\ref{Realizability_1}, we embed the original system into a system
of higher order. This system  should  satisfy a-priori certain
structural conditions that are formulated in Corollaries
\ref{single_dim_par} and \ref{independent_partition_cor}. These
two ideas (design of an algorithm aiming for its best properties,
not its realizability, and design of an embedding for realization)
result in a new method, which is shown to be applicable to a
sufficiently large class of systems with nonlinear
parameterization, e.g., like those given by systems
(\ref{cascade__full_2}).  It is very important that no damping or
discontinuities are injected directly into the control
function in contrast to \cite{Kojic_2002,Lin,Lin_2002_smooth}.

In the present article we hope to have extended the scope of
applicability and performance of adaptive control algorithms. Our results
to date are applicable to the full-state feedback case only.
Extension of the results to the output-feedback case remains a future study topic.

%\section{Acknowledgments}

%The authors are thankful to Professor Terekhov Valery
%Alexandrovich, Saint-Petersburg State Electrical Engineering
%University for his numerous remarks and helpful suggestions during
%the preparation of the manuscript.

%\begin{figure}
%\begin{center}
%\includegraphics[width=\textwidth]{cascade1par_4contr.eps}
%\end{center}
%\end{figure}

\bibliographystyle{plain}
\bibliography{fffinal_reduced}

%\bibliographystyle{IEEEtran}
%\bibliography{fffinal_reduced}
%\small

\section{Appendix 1}

%\small

%\baselineskip 6mm

In this section we consider auxiliary and technical results that
are used in the paper. Let the system dynamics with respect to the
function $\psi(\bfx,t)$ be described as follows:
\begin{eqnarray}\label{plant2_psi_ins_eps}
\dpsi=-\varphi(\psi) + z(\bfx,\thetavec,t) -
z(\bfx,\hat{\thetavec},t)+\varepsilon(t),
\end{eqnarray}
where function $\varepsilon: R_+\rightarrow R$, $\varepsilon\in
C^0$ models unknown disturbances due to unmodeled dynamics or
measurement errors. In addition, we  will assume that the
adaptation algorithms are affected by a disturbance:
\begin{eqnarray}\label{alg_dist0}
\dot{\hat{\thetavec}}&=&\Gamma
((\dpsi+\varphi(\psi(\bfx,t)))\alphavec(\bfx,t)+\delta(t)), \ \
\delta:R_+\rightarrow R^d, \  \delta\in C^0 \nonumber
\end{eqnarray}
\begin{lem}\label{lemma_dist_feedback_alg} Let error model
(\ref{plant2_psi_ins_eps}) be given, $\delta,\varepsilon\in
L_\infty$, $|\varphi(\psi)|>K|\psi|$, $K>0$ and Assumptions
\ref{boundedness_psi},
\ref{certainty_equivalence_part1}--\ref{alpha_upper_bound} hold
for $\varepsilon\equiv 0$. Then $\psi(\bfx,t)$ and
$\hat{\thetavec}$, $\bfx(t)$ are bounded for the error model
(\ref{plant2_psi_ins_eps}) with algorithm
\begin{eqnarray}\label{alg_dist}
\dot{\hat{\thetavec}}&=&(\Gamma
(\dpsi+\varphi(\psi(\bfx,t)))\alphavec(\bfx,t)+\delta(t)-\lambda\hat{\thetavec}),
\  \lambda >0.
\end{eqnarray}
%If in addition $\varepsilon(t),\delta(t)\rightarrow 0$ as
%$t\rightarrow \infty$ then $\psi(\bfx,t)$ converges into
%\begin{eqnarray}
% |\varphi(\psi)|\leq
% \Delta_0+\sqrt\frac{\lambda}{{(D-D_1)}}\|{\hat\thetavec}^{\ast}\|,
%\end{eqnarray}
%where $\Delta_0>0$ is arbitrary small.
%
\end{lem}
{\it Lemma \ref{lemma_dist_feedback_alg} proof.} Denote
$Q(\psi)=\int_0^{\psi}\varphi(\varsigma)d\varsigma$ and consider
the following function
\begin{eqnarray}\label{stand_lf}
V(\psi,\hat{\thetavec},{\hat\thetavec}^{\ast})=2(D-D_1)Q(\psi) +
0.5\|\hat\thetavec-\hat{\thetavec}^{\ast}\|^2_{\Gamma^{-1}}.
\end{eqnarray}
Its derivative satisfies the following ($z(\bfx,\thetavec,t) -
z(\bfx,\hat{\thetavec},t)+ \varepsilon(t)=\varphi(\psi)+\dpsi$ due
to equation (\ref{plant2_psi_ins_eps})):
\begin{eqnarray}\label{ldist_eq0}
\dot{V}&=&-2(D-D_1) \varphi(\psi)\dpsi +
(\hat{\thetavec}-\hat{\thetavec}^{\ast})^{T}((\varphi(\psi)+\dpsi)\alphavec(\bfx,t)+\delta(t)+\lambda\hat{\thetavec})=
-2(D-D_1)\varphi\dpsi - \nonumber \\
& & (z(\bfx,\thetavec,t) +
z(\bfx,\hat{\thetavec},t))(\hat{\thetavec}-{\hat\thetavec}^{\ast})^{T}\alphavec(\bfx,t)+\varepsilon(t)(\hat{\thetavec}-{\hat\thetavec}^{\ast})^{T}\alphavec(\bfx,t)+(\hat{\thetavec}-{\hat\thetavec}^{\ast})^{T}(\delta(t)-\lambda
\hat{\thetavec})\nonumber
\end{eqnarray}
From Assumptions \ref{alpha}, \ref{alpha_upper_bound} it follows
that
\begin{eqnarray}\label{ldist_eq01}
& & -(z(\bfx,\hat{\thetavec},t) -
z(\bfx,{\thetavec},t))(\hat{\thetavec}-{\hat\thetavec}^{\ast})^{T}\alphavec(\bfx,t)+\varepsilon(t)(\hat{\thetavec}-{\hat\thetavec}^{\ast})^{T}\alphavec(\bfx,t)-\frac{D_1\varepsilon^2(t)}{4}+\frac{D_1\varepsilon^2(t)}{4}\leq\nonumber\\
&\leq & - D (z(\bfx,\thetavec,t) - z(\bfx,\hat{\thetavec},t))^2 +
 D_1 |\varepsilon(t)| |z(\bfx,\thetavec,t) - z(\bfx,\hat{\thetavec},t)|
 -\frac{D_1\varepsilon^2(t)}{4}+\frac{D_1\varepsilon^2(t)}{4}\nonumber\\
&=& -(D-D_1)(z(\bfx,\thetavec,t) - z(\bfx,\hat{\thetavec},t))^2 -
D_1\left(|z(\bfx,\thetavec,t) -
z(\bfx,\hat{\thetavec},t)|-\frac{|\varepsilon(t)|}{2}\right)^2+\frac{D_1\varepsilon^2(t)}{4}
\end{eqnarray}
Then
\begin{eqnarray}\label{ldist_eq1}
\dot{V}& \leq & 2(D-D_1)\varphi(\psi)\dpsi -
(D-D_1)(z(\bfx,\thetavec,t) - z(\bfx,\hat{\thetavec},t))^2 - D_1
\left(|z(\bfx,\thetavec,t) -
z(\bfx,\hat{\thetavec},t)|-\frac{|\varepsilon(t)|}{2}\right)^{2}\nonumber
\\& & + \frac{D_1\varepsilon^2(t)}{4}-(\hat{\thetavec}-{\hat\thetavec}^{\ast})^{T}(\delta(t)+\lambda
\hat{\thetavec})\leq
-(D-D_1)(\varphi^{2}(\psi)+\dpsi^{2}-2(\varphi(\psi)+\dpsi)\varepsilon(t)+\varepsilon^{2}(t))+\nonumber
\\ & & + \frac{D_1\varepsilon^2(t)}{4} - \lambda
\|\hat{\thetavec}-{\hat{\thetavec}}^{\ast}\|^{2}+\|\hat\thetavec-{\hat\thetavec}^{\ast}\|(\|\delta(t)\|+\lambda\|{\hat\thetavec}^{\ast}\|)\nonumber
\\ &\leq&
-(D-D_1)\left(1-\frac{1}{\Delta_1^2}\right)(\varphi^{2}(\psi)+\dpsi^{2})
+ 2(D-D_1)\Delta_1^{2}\varepsilon^{2}(t)+
\frac{D_1\varepsilon^2(t)}{4}-\lambda\left(1-\frac{1}{\Delta_2^2}\right)\|\hat{\thetavec}-\hat{\thetavec}^{\ast}\|^2+\nonumber
\\ & &
\Delta_2^{2}\frac{(\|\delta(t)\|+\lambda\|{\hat\thetavec}^{\ast}\|)^2}{4\lambda},
\end{eqnarray}
where $\Delta_1$, $\Delta_2>1$. In the lemma conditions
$\delta,\varepsilon\in L_\infty$. Therefore, taking into account
estimate (\ref{ldist_eq1}) and inequality $|\varphi(\psi)|> K
|\psi|$, we conclude that derivative $\dot{V}$ is
negative-definite for any $\psi,\hat{\thetavec}$ that belong to
the following set:
\begin{eqnarray}\label{l1_1}
& & \Omega_{t>0}=\left\{ \psi,\hat{\thetavec}\left|
(D-D_1)\left(1-\frac{1}{\Delta_1^2}\right)\varphi^{2}(\psi)+\lambda\left(1-\frac{1}{\Delta_2^2}\right)\|\hat{\thetavec}-\hat{\thetavec}^{\ast}\|^2\geq\right.\right.
\nonumber \\
& &
\left.\|\varepsilon(t)\|_{\infty}^{2}\left(2(D-D_1)\Delta_1^{2}+\frac{D_1}{4}\right)+
\Delta_2^{2}\frac{(\|\delta(t)\|_{\infty}+\lambda\|{\hat\thetavec}^{\ast}\|)^2}{4\lambda}\right\}\nonumber
\end{eqnarray}
Hence $\hat\thetavec$, $\psi(\bfx,t)$ are bounded.{\it The lemma
is proven.}

%To prove the rest of the lemma notice that if
%$\varepsilon,\delta\rightarrow 0$ as $t\rightarrow\infty$, then
%there exist $t_1>0$ and arbitrary small $\Delta>0$ such that for
%any $t>t_1$ set $\Omega_{t>t_1}$ can be written as follows
%\begin{eqnarray}\label{l1_2}
%& & \Omega_{t>t_1}=\left\{ \psi,\hat{\thetavec}\left|
%(D-D_1)\left(1-\frac{1}{\Delta_1^2}\right)\varphi^{2}(\psi)+\lambda\left(1-\frac{1}{\Delta_2^2}\right)\|\hat{\thetavec}-\hat{\thetavec}^{\ast}\|^2\geq\right.\right.
%\left.\Delta+
%\Delta_2^{2}\frac{(\lambda\|{\hat\thetavec}^{\ast}\|)^2}{4\lambda}\right\}\nonumber
%\end{eqnarray}
%Select $\Delta_1=\Delta_2$, then
%\begin{eqnarray}\label{l1_3} & &
%\Omega_{t>t_1}\subset\left\{ \psi,\hat{\thetavec}\left|
%(D-D_1)\varphi^{2}(\psi)+\lambda\|\hat{\thetavec}-\hat{\thetavec}^{\ast}\|^2\geq\right.\right.
%\left.\frac{\Delta}{\Delta_1^2-1}+
%\frac{1}{\Delta_1^{2}-1}\frac{(\lambda\|{\hat\thetavec}^{\ast}\|)^2}{4\lambda}\right\}\nonumber\subset\nonumber\\
%& & \left\{ \psi,\hat{\thetavec}\left|
%(D-D_1)\varphi^{2}(\psi)+\lambda\|\hat{\thetavec}-\hat{\thetavec}^{\ast}\|^2\geq\right.\right.
%\left.\frac{\Delta}{\Delta_1^2-1}+
%\lambda\|{\hat\thetavec}^{\ast}\|^2\right\}\subset\nonumber \\
%& & \left\{ \psi,\hat{\thetavec}\left| \ \
%|\varphi(\psi)|\geq\right.\right. \left.
%\sqrt{\frac{\Delta}{(D-D1)(\Delta_1^2-1)}}+
%\sqrt\frac{\lambda}{{D-D_1}}\|{\hat\thetavec}^{\ast}\|\right\}\nonumber
%\end{eqnarray}
%%{\bf Why $\ge$, not $\le$ (see sign in (69))?}
%{\it The lemma is proven.}

Let $\varepsilon\in L_2$ and $\delta(t)\equiv 0$.  In this case it
is possible to show that the control goal is reached in the closed
loop system with slightly modified version of algorithm
(\ref{alg1_a}).
\begin{lem}\label{L_2_theorem_eps} Let the following error model
be given
\begin{eqnarray}\label{plant2_psi_ins_var}
\dpsi=-\varphi(\psi)(1+F(t)) + z(\bfx,\thetavec,t) -
z(\bfx,\hat{\thetavec},t)+\varepsilon(t),
\end{eqnarray}
where  $F: R_+\rightarrow R_+$, $F(t)\in C^0$, $\varepsilon(t)\in
C^0$, $\varepsilon(t)\in L_2$, Assumptions
\ref{certainty_equivalence_part1}-\ref{alpha_upper_bound} hold for
$F(t)\equiv 0$, $\varepsilon(t)\equiv 0$ and adaptation algorithm
satisfy equation
\begin{eqnarray}\label{alg1_b}
\dot{\hat\thetavec}=\Gamma(\dpsi+\varphi(\psi)(1+F(t))\alphavec(\bfx,t).
\end{eqnarray}
Then

1) $\psi(\bfx,t)\in L_\infty$, $\varphi(\psi(\bfx,t))\in L_2\cap
L_\infty$, $\sqrt{F(t)}\varphi(\psi(\bfx,t))\in L_2$;
$\hat{\thetavec}\in L_\infty$

2) $z(\bfx,\thetavec,t)-z(\bfx,\hat{\thetavec},t)\in L_2$.\\
If $F(t)\in L_\infty$ then

3) $\dpsi\in L_2$.\\
If in addition functions $\varepsilon(t)\in L_\infty$ and function
$z(\bfx,\hat{\thetavec},t)$ is locally bounded with respect to
$\bfx$, $\hat{\thetavec}$, uniformly bounded with respect to $t$
then

4) $\psi(\bfx,t)\rightarrow 0$ as $t\rightarrow\infty$.

\end{lem}
{\it Lemma \ref{L_2_theorem_eps} proof.} Function
$\varepsilon(t)\in L_2$, therefore integral
$\int_t^{\infty}\varepsilon^2(\tau)d\tau<\infty$. Consider the
following function
\[
V_{\hat{\thetavec}}(\hat{\thetavec},\hat{\thetavec}^{\ast},t)=\frac{D_1}{4}\int_t^{\infty}\varepsilon^2(\tau)d\tau
+
\frac{1}{2}\|\hat{\thetavec}-\hat{\thetavec}^{\ast}\|^2_{\Gamma^{-1}}.
\]
Its time-derivative can be written as follows:
\begin{eqnarray}
\dot{V}_{\hat{\thetavec}}&=&-1/4{D_1}\varepsilon^{2}(t)+(\varphi(\psi)(1+F(t))+\dpsi)(\hat{\thetavec}-\hat{\thetavec}^{\ast})^{T}\alphavec(\bfx,t)\nonumber
\\ &=&-1/4{D_1}\varepsilon^{2}(t)+(z(\bfx,\thetavec,t)-z(\bfx,\hat{\thetavec},t)+\varepsilon(t))(\hat{\thetavec}-\hat{\thetavec}^{\ast})^{T}\alphavec(\bfx,t)\nonumber
\end{eqnarray}
Taking into account inequality (\ref{ldist_eq01}) we can write the
following estimate for $\dot{V}_{\hat\theta}$:
\begin{eqnarray}\label{d_theta}
\dot{V}_{\hat\theta}\leq
-(D-D_1)(z(\bfx,\thetavec,t)-z(\bfx,\hat{\thetavec},t))^{2}=-(D-D_1)(\varphi(\psi)(1+F(t))+\dpsi-\varepsilon(t))^{2}
\end{eqnarray}

It follows from (\ref{d_theta}) that
$z(\bfx,\hat{\thetavec},t)-z(\bfx,\hat{\thetavec}^{\ast},t)\in
L_2$. Let us denote
$\mu(t)=\varepsilon(t)+z(\bfx(t),\hat{\thetavec}(t),t)-z(\bfx(t),\hat{\thetavec}^{\ast},t)$.
Taking this equality into account, error model
(\ref{plant2_psi_ins_var}) can be written as follows
$\dpsi=-(1+F(t))\varphi(\psi)+\mu(t)$, where function $\mu(t)\in
L_2$ as a sum of the functions from $L_2$. Consider the following
nonnegative function $V_1(\psi,t)$:
\[
V_1(\psi,t)=\int_0^\psi\varphi(\xi)d\xi+\frac{1}{4}\int_t^\infty\mu^2(\tau)d\tau
\]
Its time-derivative is:
\begin{eqnarray}\label{eq_ins_lemma2}
\dot{V}_1&=&-\varphi(\psi)(1+F(t))\varphi(\psi)+\varphi(\psi)\mu(t)-\frac{1}{4}\mu^2(t)
-F(t)\varphi^2(\psi)-\left(\varphi(\psi)-\frac{1}{2}\mu(t)\right)^{2}
\end{eqnarray}
It follows from inequality (\ref{eq_ins_lemma2}) that
$\psi(\bfx,t), \ \varphi(\psi(\bfx(t),t))\in L_\infty$.
Furthermore, $\sqrt{F(t)}\varphi(\psi(\bfx(t),t))\in L_2$ and
$(\varphi(\psi(\bfx(t),t))-\mu(t)/2)\in L_2$. Given that
$\mu(t)\in L_2$ it is clear that $\varphi(\psi(\bfx(t),t))\in
L_2$. Hence statements 1) and 2) of the lemma are proven. Let
$F(t)\in L_\infty$ then $(1+F(t))\varphi(\psi(\bfx(t),t))\in L_2$
and therefore according to (\ref{plant2_psi_ins_var}) $\dpsi\in
L_2$ as well. Thus statement 3) is proven. To show that 4) holds
it is sufficient to notice that $\hat{\thetavec}$ is bounded due
to (\ref{d_theta}). According to Assumption \ref{boundedness_psi}
state $\bfx$ is bounded as $\psi(\bfx,t)$ is bounded. Then $\dpsi$
is bounded if $\varepsilon(t)$ is bounded and function
$\bfz(\bfx,\thetavec,t)$ is locally bounded. Hence applying
Barbalat's lemma we conclude that $\psi(\bfx,t)\rightarrow 0$ as
$t\rightarrow\infty$. {\it The lemma is proven.}

\begin{lem}\label{lemma_embed_full_cascade} Let system:
\begin{eqnarray}\label{embed_ref_sys}
\dot{x}_i&=&f_i(x_1,\dots,x_i,\thetavec_i)+\beta_i(\bfx,t), \ \ \
i=1,\dots,n
\end{eqnarray}
and smooth function $u(\bfx,\bfz,\thetavec_0): R^{n}\times
R^m\times R^d \rightarrow R$, be given. Let us assume that
$\thetavec_0\in \Omega_0$, $\Omega_0$ be bounded and there exist
smooth functions $\bar{F}(\bfx,\bfx',\bfz)$,
$\bar{D}_i(\bfx,\bfx')$, $i=1,\dots,n$ such that the following
properties hold:
\[
(u(\bfx,\bfz,\thetavec_0)-u(\bfx',\bfz,\thetavec_0))^2\leq
\|\bfx-\bfx'\|^2 \bar{F}^2(\bfx,\bfx',\bfz) \ \ \ \forall
\thetavec_0\in \Omega_0, \bfx,\bfx'\in R^n
\]
\begin{eqnarray}
(f_i(x_1,\dots,x_i,\thetavec_i)-f_i(x_1',\dots,x_i',\thetavec_i))^2&\leq&
\|\tilde{\bfx}_i-\tilde{\bfx}_i'\|^2 \bar{D}_i^2(\bfx_i,\bfx_i') \
\ \ \forall \thetavec_i\in \Omega_\theta, \bfx_i,\bfx_i'\in
R^n\nonumber
\\
\tilde{\bfx}_i=(x_1,\dots,x_i,0,\dots,0)^{T}, \ & &
\tilde{\bfx}_i'=(x_1',\dots,x_i',0,\dots,0)^{T} \nonumber
\end{eqnarray}

Let us also assume that there exist $\alphavec_i(\bfx)$ such that
Assumptions \ref{alpha}, \ref{alpha_upper_bound} hold for the
functions $f_i(x_1,\dots,x_i,\thetavec_i)$ respectively.

Then there exist $\xivec(t): R\rightarrow R^n$, $\nuvec(t):
R\rightarrow R^m$, smooth functions $\bff_{\xi}(\cdot)$,
$\bff_{\nuvec}(\cdot)$ and corresponding system:
\begin{eqnarray}\label{embed_track_sys}
\dot{\xivec}&=&\bff_{\xi}(\bfx,\xivec,\bfz,\nuvec),
\ \xivec_0\in R^{n} \nonumber \\
\dot{\nuvec}&=&\bff_{\nu}(\bfx,\xivec,\bfz,\nuvec), \ \nuvec_0 \in
R^m
\end{eqnarray}
such that

1) $u(\bfx,\bfx,\thetavec_0)-u(\bfq_i,\bfz,\thetavec_0)\in L_2, \
i=1,\dots,n$
\[
\bfq_i=(\xi_1,\dots,\xi_{i},x_{i+1},\dots,x_n)^{T};
\]

2)
$f_i(x_1,\dots,x_i,\thetavec_i)-f_i(\xi_1,\dots,\xi_{i-1},x_i,\thetavec_i)\in
L_2 , \ i=2,\dots,n$;

3) $\bfx\in L_\infty \Rightarrow \xivec,\nuvec\in L_\infty$.

%Furthermore, if $\bar{F}(\bfx,\bfx',\bfz)\in C^{r}$,
%$\bar{D}_i(\bfx,\bfx')\in C^{r}$,
%$f_i(x_1,\dots,x_i,\thetavec_i)\in C^{r}$, $\beta_i(\bfx,t)\in
%C^r$, $i=1,\dots,n$ then
\end{lem}
{\it Lemma \ref{lemma_embed_full_cascade} proof.} For the sake of
notational convenience we would like to use the following
notations:
\[
f_i(x_1,\dots,x_i,\thetavec_i)=f_i(\bfx,\thetavec_i), \
\psi_{\xi_i}=x_i-\xi_i
\]
\begin{eqnarray}\label{eq_8.5_lemma_embed}
\varepsilon_i(t)=f_i(\bfq_{i-2},\thetavec_i)-f_i(\bfq_{i-1},\thetavec_i),
\ i=2,\dots,n.
\end{eqnarray}

Consider the following system of differential equations:
\begin{eqnarray}\label{eq_9_lemma_embed}
\dot{\xi}_i&=&((\bar{F}_i^2(\bfq_{i-1},\bfq_{i},\bfz)+\sum_{j=i}^{k}\bar{D}_{j+1}^2(\bfq_{i-1},\bfq_{i}))+1)(x_i-\xi_i)+f_i(\bfq_{i-1},{\hat\thetavec}_{\xi_i})+\beta_i(\bfx,t)\\
\dot{\hat{\thetavec}}_{\xi_i}&=&\gamma_{\xi_i}(\psi_{\xi_i}((\bar{F}_i^2(\bfq_{i-1},\bfq_i,\bfz)+\sum_{j=i}^{k}\bar{D}_{j+1}^2(\bfq_{i-1},\bfq_i))+1)+\dpsi_{\xi_i})\alphavec_i(\bfq_{i-1}),
\ \gamma_{\xi_i}>0 \nonumber, \ i=1,\dots,k,
\end{eqnarray}
where
$\bar{F}_{i}(\bfq_{i-1},\bfq_{i},\bfz)=\bar{F}(\bfq_{i-1},\bfq_{i},\bfz)$.
Taking into account (\ref{embed_ref_sys}) and
(\ref{eq_9_lemma_embed}) let us write the following error model:
\begin{eqnarray}\label{eq_10_lemma_embed}
\dot{\psi}_{\xi_i}&=&-((\bar{F}_i^2(\bfq_{i-1},\bfq_i,\bfz)+\sum_{j=i}^{k}\bar{D}_{j+1}^2(\bfq_{i-1},\bfq_{i}))+1)\psi_{\xi_i}-f_i(\bfq_{i-1},{\hat\thetavec}_{\xi_i})+f_i(\bfx,\thetavec_i)\\
\dot{\hat{\thetavec}}_{\xi_i}&=&\gamma_{\xi_i}(\psi_{\xi_i}((\bar{F}_i^2(\bfq_{i-1},\bfq_i,\bfz)+\sum_{j=i}^{k}\bar{D}_{j+1}^2(\bfq_{i-1},\bfq_{i}))+1)+\dpsi_{\xi_i})\alphavec_i(\bfq_{i-1}),
\ \gamma_{\xi_i}>0 \nonumber
\end{eqnarray}
It is clear that trajectories of system (\ref{eq_9_lemma_embed})
for $k=1$ satisfy the following condition
$u(\bfx,\bfz,\thetavec_0)-u(\bfq_1,\bfz,\thetavec_0)\in L_2$ (this
follows directly from Lemma \ref{L_2_theorem_eps}). Consider the
case when $k=2$. Taking into account the equations for
$\dpsi_{\xi_1}$ and $\dot{\hat{\thetavec}}_{\xi_1}$ we can derive
from Lemma \ref{L_2_theorem_eps} that
\[
\sqrt{\bar{F}_1^2(\bfq_0,\bfq_1,\bfz)+\bar{D}_{2}^2(\bfq_0,\bfq_1)}(x_1-\xi_1)\in
L_2, \ x_1-\xi_1\in L_{\infty},  \  \hat{\thetavec}_{\xi_1}\in
L_\infty
\]
Hence $\bar{F}_1(\bfq_0,\bfq_1,\bfz)(x_1-\xi_1)\in L_2, \
\bar{D}_{2}(\bfq_0,\bfq_1)(x_1-\xi_1)\in L_2$ as
$\sqrt{\bar{F}_1^2(\bfq_0,\bfq_1,\bfz)+\bar{D}_{2}^2(\bfq_0,\bfq_1)}\geq
|\bar{F}_1(\bfq_0,\bfq_1,\bfz)|$ and
$\sqrt{\bar{F}_1^2(\bfq_0,\bfq_1,\bfz)+\bar{D}_{2}^2(\bfq_0,\bfq_1)}\geq
|\bar{D}_{2}(\bfq_0,\bfq_1)|$. Therefore, we can conclude that
\[
u(\bfq_0,\bfz,\thetavec_0)-u(\bfq_1,\bfz,\thetavec_0)\in L_2, \
f_2(\bfq_0,\thetavec_2)-f_2(\bfq_1,\thetavec_2)\in L_2
\]
as $|u(\bfq_0,\bfz,\thetavec_0)-u(\bfq_1,\bfz,\thetavec_0)|\leq
|\bar{F}_1(\bfq_0,\bfq_1,\bfz)(x_1-\xi_1)|\in L_2$, $
|f_2(\bfq_0,\thetavec_2)-f_2(\bfq_1,\thetavec_2)|\leq
|\bar{D}_{2}(\bfq_0,\bfq_1)(x_1-\xi_1)|\in L_2$. Notice that
$\varepsilon_2(t)=f_2(\bfq_0,\thetavec_2)-f_2(\bfq_1,\thetavec_2)=f_2(\bfx,\thetavec_2)-f_2(\bfq_1,\thetavec_2)$.
Therefore, the equations for $\dpsi_{\xi_2}$ become as follows:
\begin{eqnarray}\label{eq_11_lemma_embed}
\dot{\psi}_{\xi_2}&=&-(\bar{F}_2^2(\bfq_1,\bfq_2,\bfz)+1)\psi_{\xi_2}-f_2(\bfq_{1},{\hat\thetavec}_{\xi_2})+f_2(\bfq_{1},\thetavec_2)+\varepsilon_2(t)\\
\dot{\hat{\thetavec}}_{\xi_2}&=&\gamma_{\xi_2}(\psi_{\xi_2}(\bar{F}_2^2(\bfq_1,\bfq_2,\bfz)+1)+\dpsi_{\xi_2})\alphavec_2(\bfq_{1}),
\ \gamma_{\xi_2}>0 \nonumber, \ \varepsilon(t)\in L_2
\end{eqnarray}
Hence, applying Lemma \ref{L_2_theorem_eps} to system
(\ref{eq_10_lemma_embed}), (\ref{eq_11_lemma_embed}) and taking
into account that
$u(\bfq_1,\bfz,\thetavec_0)-u(\bfq_2,\bfz,\thetavec_0)|\leq|\bar{F}_2(\bfq_1,\bfq_2,\bfz)(x_2-\xi_2)|$
we can conclude that
\begin{eqnarray}
& &u(\bfq_0,\bfz,\thetavec_0)-u(\bfq_1,\bfz,\thetavec_0)\in L_2, \
x_1-\xi_1\in L_\infty, \ \hat{\thetavec}_{\xi_1}\in L_\infty
\nonumber \\
& & u(\bfq_1,\bfz,\thetavec_0)-u(\bfq_2,\bfz,\thetavec_0) \in
L_2,\ x_2-\xi_2\in L_\infty, \ \hat{\thetavec}_{\xi_2}\in L_\infty
\nonumber \\
& & f_2(\bfq_0,\thetavec_2)-f_2(\bfq_1,\thetavec_2)\in L_2,
\nonumber
\end{eqnarray}
and, subsequently,
$u(\bfq_0,\bfz,\thetavec_0)-u(\bfq_2,\bfz,\thetavec_0)\in L_2$ as
a sum of two signals from $L_2$.

Let us now consider arbitrary $2<k\leq n$. It follows from Lemma
\ref{L_2_theorem_eps} that for the error model with respect to
function $\psi_{\xi_1}$:
\[
\dot{\psi}_{\xi_1}=-((\bar{F}_1^2(\bfq_0,\bfq_1,\bfz)+\sum_{j=1}^{k}\bar{D}_{j+1}^2(\bfq_0,\bfq_1))+1)\psi_{\xi_1}-f_1(\bfq_{0},{\hat\thetavec}_{\xi_1})+f_1(\bfx,\thetavec_1)
\]
and corresponding subsystem
\[
\dot{\hat{\thetavec}}_{\xi_1}=\gamma_{\xi_1}(\psi_{\xi_1}((\bar{F}_1^2(\bfq_0,\bfq_1,\bfz)+\sum_{j=1}^{k}\bar{D}_{j+1}^2(\bfq_0,\bfq_1))+1)+\dpsi_{\xi_1})\alphavec_1(\bfq_{0}),
\]
one can derive that
$\sqrt{\bar{F}_1^2(\bfq_0,\bfq_1,\bfz)+\sum_{j=1}^{k}\bar{D}_{j+1}^2(\bfq_0,\bfq_1)}(x_1-\xi_1)\in
L_2$. This in consequence implies that
\begin{eqnarray}
&u(\bfx,\bfz,\thetavec_0)-u(\bfq_1,\bfz,\thetavec_0)\in
L_2&\nonumber\\
&f_i(\bfx,\thetavec_i)-f_i(\bfq_1,\thetavec_i)\in L_2&,
i=2,\dots,k \nonumber
\end{eqnarray}
Hence we can write the error model for $\psi_{\xi_2}$ in
(\ref{eq_10_lemma_embed}) in the following form
\[
\dot{\psi}_{\xi_2}=-((\bar{F}_2^2(\bfq_1,\bfq_2,\bfz)+\sum_{j=3}^{k}\bar{D}_{j+1}^2(\bfq_1,\bfq_2))+1)\psi_{\xi_2}-f_2(\bfq_{1},{\hat\thetavec}_{\xi_2})+f_2(\bfq_1,\thetavec_2)+\varepsilon_2(t)
\]
where $\varepsilon_2(t)\in L_2$. It follows from Lemma
\ref{L_2_theorem_eps} that
\begin{eqnarray}
&u(\bfq_1,\bfz,\thetavec_0)-u(\bfq_2,\bfz,\thetavec_0)\in
L_2&\nonumber\\
&f_i(\bfq_1,\thetavec_i)-f_i(\bfq_2,\thetavec_i)\in L_2&,
i=3,\dots,k \nonumber
\end{eqnarray}
in system (\ref{eq_9_lemma_embed}). Notice also that
$f_3(\bfx,\thetavec_3)-f_3(\bfq_2,\thetavec_3)\in L_2$ as a sum of
two functions from $L_2$. By the similar reasoning it is can be
shown that for any $2\leq i \leq n$ we can represent the error
model system (\ref{eq_10_lemma_embed}) as follows
\[
\dot{\psi}_{\xi_i}=-((\bar{F}_i^2(\bfq_{i-1},\bfq_i,\bfz)+\sum_{j=i+1}^{k}\bar{D}_{j+1}^2(\bfq_{i-1},\bfq_i)))+1)\psi_{\xi_i}-f_i(\bfq_{i-1},{\hat\thetavec}_{\xi_i})+f_i(\bfq_{i-1},\thetavec_i)+\varepsilon_i(t),
\]
where $\varepsilon_i(t)\in L_2$. Therefore, using Lemma
\ref{L_2_theorem_eps} again we can conclude that
\begin{eqnarray}
&u(\bfq_{j-1},\bfz,\thetavec_0)-u(\bfq_{j},\bfz,\thetavec_0)\in
L_2&\nonumber\\
&f_i(\bfq_{j-1},\thetavec_i)-f_i(\bfq_{j},\thetavec_i)\in L_2,& \
x_i-\xi_i\in L_\infty, \ \ \hat{\thetavec}_{\xi_i}\in L_\infty, \
\ i=j,\dots, n, \nonumber
\end{eqnarray}
The last, however, implies that
$u(\bfx,\bfz,\thetavec_0)-u(\bfq_i,\bfz,\thetavec_0)\in L_2$. In
order to complete the proof we have to make sure that system
(\ref{eq_9_lemma_embed}) is physically realizable. In particular,
realization of subsystems
\begin{eqnarray}\label{eq_12_lemma_embed}
\dot{\hat{\thetavec}}_{\xi_i}=\gamma_{\xi_i}(\psi_{\xi_i}((\bar{F}_i^2(\bfq_{i-1},\bfq_i,\bfz)+\sum_{j=i}^{k}\bar{D}_{j+1}^2(\bfq_{i-1},\bfq_{i}))+1)+\dpsi_{\xi_i})\alphavec_i(\bfq_{i-1}),
\ \gamma_{\xi_i}>0, \ i=1,\dots,k,
\end{eqnarray}
shell not be dependent on any uncertainties $\thetavec_i$. It
follows, however, from Theorem \ref{theorem_realization1} that
there are realizations of algorithms (\ref{eq_12_lemma_embed}) in
finite form:
\begin{eqnarray}\label{auxilary_alg21}
\hat{\thetavec}_{\xi_i}(\bfq_{i-1},\xi_i,t) &=&\gamma_{\xi_i}(\hat{\thetavec}_{\xi_i,P}(\bfq_{i-1},\xi_i)+\hat{\thetavec}_{\xi_i,I}(t)), \ \gamma_{\xi_i}>0\nonumber\\
\hat{\thetavec}_{\xi_i,P}(\bfq_{i-1},\xi_i)&=&\psi_{\xi_i}(x_i,\xi_i)\alphavec_i(\bfq_{i-1})-\Psi_{\xi_i}(\bfq_{i-1},\xi_i)\nonumber\\
\dot{\hat{\thetavec}}_{\xi_i,I}&=&((\bar{F}_i^2(\bfq_{i-1},\bfq_i,\bfz)+\sum_{j=i}^{k}\bar{D}_{j+1}^2(\bfq_{i-1},\bfq_i))+1)\psi_{\xi_i}(x_i,\xi_i)\alphavec_i(\bfq_{i-1})+\nonumber \\
& & \sum_{j=1}^{i}\frac{\pd \Psi_{\xi_i}(\bfq_{i-1},\xi_i)}{\pd
\xi_j}\dot{\xi}_j-\sum_{j=1}^{i-1}\psi_{\xi_i}(x_i,\xi_i)\frac{\pd
\alphavec_i(\bfq_{i-1})}{\pd
\xi_j}\dot{\xi}_j \nonumber \\
& &
\Psi_{\xi_i}(\bfq_{i-1},\xi_i)=\int_{x_i(0)}^{x_i(t)}\psi_{\xi_i}(x_i,\xi_i)\frac{\pd
\alphavec_i(\bfq_{i-1}) }{\pd x_i}d x_i
\end{eqnarray}

Notice also that if $\bfx\in L_{\infty}$ then $\xivec\in
L_{\infty}$ and hence
$\hat{\thetavec}_{\xi_i,P}(\bfq_{i-1},\xi_i)\in L_\infty$ as
$\hat{\thetavec}_{\xi_i,P}(\bfq_{i-1},\xi_i)$ is smooth. Given
that
$\hat{\thetavec}_{\xi_i}=\gamma_{\xi_i}(\hat{\thetavec}_{\xi_i,P}(\bfq_{i-1},\xi_i)+\hat{\thetavec}_{\xi_i,I})$
and both $\hat{\thetavec}_{\xi_i},
\hat{\thetavec}_{\xi_i,P}(\bfq_{i-1},\xi_i)\in L_\infty$ then we
can conclude that $\hat{\thetavec}_{\xi_i,I}\in L_\infty$ for
$\bfx\in L_\infty$.

It is easy to see that denoting $\nuvec=\hat{\thetavec}_{\xi_i,I}$
we can transform system (\ref{eq_9_lemma_embed}),
(\ref{auxilary_alg21}) into (\ref{embed_track_sys}) which
satisfies statements 1)--3) of the lemma. {\it The lemma is
proven.}

\section{Appendix 2}

{\it Theorem \ref{L_2 dpsi_theorem} proof.} Let us consider the
following positive-definite function:
$V_{\hat{\thetavec}}(\hat{\thetavec},\hat{\thetavec}^{\ast})=
\frac{1}{2}\|\hat{\thetavec}-\hat{\thetavec}^{\ast}\|^2_{\Gamma^{-1}}$.
Its time-derivative according to equations (\ref{ONO}) can be
derived as follows:
$\dot{V}_{\hat{\thetavec}}(\hat{\thetavec},\hat{\thetavec}^{\ast})=(\varphi(\psi)+\dpsi)(\hat{\thetavec}-\hat{\thetavec}^{\ast})^{T}\alphavec(\bfx,t)$.
According to Assumption \ref{alpha} and equality
(\ref{plant2_psi_ins}) it is easy to see that
\begin{eqnarray}\label{parameric_deviation_derivative}
\dot{V}_{\hat{\thetavec}}(\hat{\thetavec},\hat{\thetavec}^{\ast})=-(z(\bfx,\hat{\thetavec},t)-z(\bfx,\thetavec,t))(\hat{\thetavec}-\hat{\thetavec}^{\ast})^{T}\alphavec(\bfx,t)
\leq-D(z(\bfx,\hat{\thetavec},t)-z(\bfx,\thetavec,t))^2=-D(\varphi(\psi)+\dpsi)^{2}\leq0
\end{eqnarray}
Therefore $V_{\hat{\thetavec}}$ is non-increasing (property P2) is
proven). Furthermore, integration of $\dot{V}_{\hat{\thetavec}}$
with respect to time results in
\[
V_{\hat{\thetavec}}(\hat{\thetavec}(0),\hat{\thetavec}^{\ast})-V_{\hat{\thetavec}}(\hat{\thetavec}(t),\hat{\thetavec}^{\ast})\geq
D\int_{0}^{t}(\dpsi(\tau)+\varphi(\psi(\tau)))^2 d\tau \geq 0.
\]
Function $V_{\hat{\thetavec}}$ is non-increasing and bounded from
below as  $V_{\hat{\thetavec}}\geq 0$, therefore
\[
D\int_{0}^{t}(\dpsi(\tau)+\varphi(\psi(\tau)))^2  d\tau \leq
V_{\hat{\thetavec}}(\hat{\thetavec}(0),\hat{\thetavec}^{\ast})<\infty.
\]
Hence
$(\varphi(\psi)+\dpsi)=(z(\bfx,\thetavec,t)-z(\bfx,\hat{\thetavec},t))=(z(\bfx,\hat{\thetavec}^{\ast},t)-z(\bfx,\hat{\thetavec},t))\in
L_2$ (property P3)).

To prove property P1) let us consider the following function:
$V(\psi,\hat{\thetavec},\hat{\thetavec}^{\ast})=2D Q(\psi)
+V_{\hat{\thetavec}}(\hat{\thetavec},\hat{\thetavec}^{\ast})$,
where $Q(\psi)=\int_0^{\psi}\varphi(\varsigma)d\varsigma$.
Function $V(\psi,\hat{\thetavec})$ is positive-definite with
respect to $\psi(\bfx,t)$ and
$\hat{\thetavec}-\hat{\thetavec}^{\ast}$ because of Assumption
\ref{certainty_equivalence_part2}. Its time-derivative obeys
inequality:
$\dot{V}(\psi,\hat{\thetavec},\hat{\thetavec}^{\ast})\leq2D\varphi(\psi)\dpsi
- D(\dpsi+\varphi(\psi))^{2}=-D\varphi^{2}(\psi)-D\dpsi^{2} \leq
0$.

Therefore, function
$V(\psi,\hat{\thetavec},\hat{\thetavec}^{\ast})$ is bounded and
non-increasing. Furthermore
\begin{eqnarray}\label{t1_L2}
\infty>V(\psi(\bfx(0),0),\hat{\thetavec}(0),\hat{\thetavec}^{\ast})\geq
V(\psi(\bfx(0),0),\hat{\thetavec}(0),\hat{\thetavec}^{\ast})-V(\psi(\bfx(t),t),\hat{\thetavec}(t),\hat{\thetavec}^{\ast})&\geq&
D\int_{0}^t\varphi^2(\psi(\bfx(\tau),\tau))d\tau\geq 0\nonumber \\
\infty>V(\psi(\bfx(0),0),\hat{\thetavec}(0),\hat{\thetavec}^{\ast})\geq
V(\psi(\bfx(0),0),\hat{\thetavec}(0),\hat{\thetavec}^{\ast})-V(\psi(\bfx(t),t),\hat{\thetavec}(t),\hat{\thetavec}^{\ast})&\geq&
D\int_{0}^t\dpsi^2(\tau)d\tau\geq 0.
\end{eqnarray}
or, equivalently, $\dpsi(t)\in L_2$, $\varphi(\psi(t))\in L_2$.
Hence, property P1) is proven as well. The $L_2$ norm bounds
(\ref{L_2_L_inf_performance}) for $\varphi(\psi)$ and $\dpsi$
follow immediately from inequality (\ref{t1_L2}):
\begin{eqnarray}
\|\varphi(\psi)\|_2^2 & \leq&
{D}^{-1}V(\psi(\bfx(0),0),\hat{\thetavec}(0),\hat{\thetavec}^{\ast}),
\ \|\dpsi\|_2^2  \leq
{D}^{-1}V(\psi(\bfx(0),0),\hat{\thetavec}(0),\hat{\thetavec}^{\ast})\nonumber
\end{eqnarray}
The $L_\infty$ norm bound for $\psi(\bfx(t),t)$ results from the
inequality:
$V(\psi(\bfx(0),0),\hat{\thetavec}(0),\hat{\thetavec}^{\ast})-V(\psi(\bfx(t),t),\hat{\thetavec}(t),\hat{\thetavec}^{\ast})\geq0$.
Consider function $\Lambda$ defined as
$\Lambda(d)=\max_{|\psi|}\{|\psi| \ | \
\int_{0}^{|\psi|}\varphi(\varsigma)d\varsigma=d\}$ and notice that
it is monotonic and nondecreasing. Therefore, given that
$\int_{0}^{\psi(\bfx(t),t)}\varphi(\varsigma)d\varsigma\leq
\frac{1}{2D}V(\psi(\bfx(0),0),\hat{\thetavec}(0),\hat{\thetavec}^{\ast})$
we can conclude that $|\psi|\leq
\Lambda\left(\frac{1}{2D}V(\psi(\bfx(0),0),\hat{\thetavec}(0),\hat{\thetavec}^{\ast})\right)$.
To prove property P4) notice that function
$V(\psi(\bfx(t),t),\hat{\thetavec}(t),\hat{\thetavec}^{\ast})$ is
bounded. Hence by Assumption \ref{certainty_equivalence_part2}
function $\psi(\bfx(t),t)$ is bounded as well. According to
Assumption \ref{boundedness_psi} boundedness of $\psi(\bfx(t),t)$
implies boundedness of the sate $\bfx$. In addition it is assumed
that $z(\bfx,\hat{\thetavec},t)$ is locally bounded with respect
to $\bfx,\hat\thetavec$ and uniformly bounded in $t$. Therefore
the difference $z(\bfx,\thetavec,t)-z(\bfx,\hat\thetavec,t)$ is
bounded. Furthermore, by Assumption
\ref{certainty_equivalence_part2} function $\varphi(\psi)\in C^0$
and therefore it is bounded as well given that $\psi$ is bounded.
Hence $\dpsi$ is bounded and by applying Barbalat's lemma one can
show that $\psi(\bfx(t),t)\rightarrow 0$ at $t\rightarrow\infty$.
%Property P4) is proven.

To compete the proof of the theorem consider the difference
$z(\bfx,\thetavec,t)-z(\bfx,\hat\thetavec,t)$. Let function
$\varphi\in C^1$, function $z(\bfx,\thetavec,t)$ is differentiable
in $\bfx$, $\thetavec$; derivative $ \pd
{z(\bfx,\thetavec,t)}/{\pd t}$ is bounded uniformly in $t$;
function $\alphavec(\bfx,t)$ is locally bounded with respect to
$\bfx$ and uniformly bounded with respect to $t$, then $d/dt
(z(\bfx,\thetavec,t)-z(\bfx,\hat\thetavec,t))$ is bounded. On the
over hand there exists the following limit
\[
\lim_{t\rightarrow\infty}\int_0^{t}
(z(\bfx,\thetavec,\tau)-z(\bfx,\hat\thetavec,\tau))^{2}=\int_0^{\infty}
(z(\bfx,\thetavec,\tau)-z(\bfx,\hat\thetavec,\tau))^{2} \leq
\frac{1}{D}V_{\hat{\thetavec}}(\hat{\thetavec}(0),\hat{\thetavec}^{\ast})
\]
as $\int_0^{t}
(z(\bfx,\thetavec,\tau)-z(\bfx,\hat\thetavec,\tau))^{2}$ is
non-decreasing and bounded from above. Hence by Barbalat's lemma
it follows that
$z(\bfx,\thetavec,\tau)-z(\bfx,\hat\thetavec,\tau)\rightarrow 0$
as $t\rightarrow\infty$. Notice also that
$\psi(\bfx(t),t)\rightarrow 0$ as $t\rightarrow \infty$. Then
$\dpsi\rightarrow 0$ as $t\rightarrow\infty$. {\it The theorem is
proven.}

{\it Proof of Proposition \ref{Exponential_Convergence}.} Consider
the following integral\footnote{That we substitute the arguments
of the functions $\dpsi(\cdot)$ and $\psi(\cdot)$ with $t$ means
that we consider them as functions of time.}
$\int_0^{t}(\dpsi(\tau)+\varphi(\psi(\tau))^2 d\tau$. It was shown
in Theorem \ref{L_2 dpsi_theorem} proof that
$\int_0^{t}(\dpsi(\tau)+\varphi(\psi(\tau))^2 d\tau\leq
\frac{1}{2D}\|\hat{\thetavec}(0)-\hat{\thetavec}^{\ast}\|^{2}_{\Gamma^{-1}}$
along system (\ref{ONO}) solutions. Let us define
$\mu(t)=\dpsi(t)+\varphi(\psi(t))$. In the other words
\begin{eqnarray}\label{prop_exp_conv_dpsi}
\dpsi=-\varphi(\psi)+\mu(t),
\end{eqnarray}
where $\int_{0}^{\infty}\mu^2(\tau)d\tau\leq
\frac{1}{2D}\|\hat{\thetavec}(0)-\hat{\thetavec}^{\ast}\|^{2}_{\Gamma^{-1}}$.
According to the proposition conditions, $\varphi(\psi)=K\psi$, it
is possible to derive the solution of equation
(\ref{prop_exp_conv_dpsi}) as follows
$\psi(t)=\psi(0)e^{-Kt}+\int_0^{t}e^{-K(t-\tau)}\mu(\tau)d\tau$.
Hence
\begin{eqnarray}\label{prop_exp_conv_psi}
|\psi(t)|&\leq& |\psi(0)|e^{- Kt} + \sqrt{\left(\int_{0}^t
e^{-K(t-\tau)}\mu(\tau)d\tau\right)^2}\leq |\psi(0)|e^{- Kt} +
\sqrt{\int_{0}^t
e^{-2K(t-\tau)}d\tau\int_0^{t}\mu^{2}(\tau)d\tau}\nonumber \\
&\leq & |\psi(0)|e^{-
Kt}+\frac{1}{2}\sqrt{\frac{1}{KD}\|\hat{\thetavec}(0)-\hat{\thetavec}^{\ast}\|^{2}_{\Gamma^{-1}}}.
\end{eqnarray}
Property P6) is thus proven. In order to prove property P7)
consider
\[
\dot{\hat{\thetavec}}=\Gamma(\dpsi+\varphi(\psi))\alphavec(\bfx,t)=\Gamma
(z(\bfx,\thetavec,t)-z(\bfx,\hat{\thetavec},t))\alphavec(\bfx,t).
\]
Function
\begin{eqnarray}
& &
D_1|\alphavec(\bfx,t)^{T}(\hat{\thetavec}-\hat{\thetavec}^{\ast})|\leq|z(\bfx,\thetavec,t)-z(\bfx,\hat{\thetavec},t))|\leq
D|\alphavec(\bfx,t)^{T}(\hat{\thetavec}-\hat{\thetavec}^{\ast})|\nonumber\\
& &
\alphavec(\bfx,t)^{T}(\hat{\thetavec}-\hat{\thetavec}^{\ast})(z(\bfx,\hat\thetavec,t)-z(\bfx,\hat{\thetavec}^{\ast},t))>0
\  \forall \ z(\bfx,\hat\thetavec^{\ast},t)\neq
z(\bfx,\hat{\thetavec},t).\nonumber
\end{eqnarray}
Therefore, there exists $D_1\leq\kappa(t)\leq D$ such that
\[
\dot{\hat{\thetavec}}=-\kappa(t)\Gamma
\alphavec(\bfx,t)^{T}(\hat{\thetavec}-\hat{\thetavec}^{\ast})\alphavec(\bfx,t)=-\kappa(t)\Gamma
\alphavec(\bfx,t)\alphavec(\bfx,t)^{T}
(\hat{\thetavec}-\hat{\thetavec}^{\ast}).
\]
Hence
\begin{eqnarray}\label{prop_exp_conv_dtheta}
\hat{\thetavec}(t)-\hat{\thetavec}^{\ast}=e^{-\Gamma
\int_{0}^{t}\kappa(\tau)\alphavec(\bfx(\tau),\tau)\alphavec(\bfx(\tau),\tau)^{T}d\tau
}(\hat{\thetavec}(0)-\hat{\thetavec}^{\ast})
\end{eqnarray}
Consider the integral
$\Gamma\int_{0}^{t}\kappa(\tau)\alphavec(\bfx(\tau),\tau)\alphavec(\bfx(\tau),\tau)^{T}d\tau$
for $t>L$
\[
\Gamma\int_{0}^{t}\kappa(\tau)\alphavec(\bfx(\tau),\tau)\alphavec(\bfx(\tau),\tau)^{T}d\tau\geq
\Gamma D_1
\int_{0}^{t}\alphavec(\bfx(\tau),\tau)\alphavec(\bfx(\tau),\tau)^{T}d\tau,
\]
where $\alphavec(\bfx(t),t)$ is persistently exciting. For any
$t>L$ there exists integer $n\geq 0$ such that $t=nL+r$, $r\in R,
0 \leq r < L$. Therefore
\[
\Gamma D_1
\int_{0}^{t}\alphavec(\bfx(\tau),\tau)\alphavec(\bfx(\tau),\tau)^{T}d\tau\geq
\Gamma D_1 n \delta I \geq \left(\frac{\Gamma D_1
\delta}{L}t-I\right).
\]
Then taking into account (\ref{prop_exp_conv_dtheta}) one can
write
\begin{eqnarray}\label{prop_exp_conv_norm_theta}
\|\hat{\thetavec}(t)-\hat{\thetavec}^{\ast}\|\leq
\|e^{\left(-\frac{\Gamma D_1
\delta}{L}t+I\right)}\|\|\hat{\thetavec}(0)-\hat{\thetavec}^{\ast}\|,
\end{eqnarray}
i. e. $\hat\thetavec(t)$ converges to $\hat{\thetavec}^{\ast}$
exponentially fast. It means that there exist positive constants
$\lambda>0$, $\lambda\neq K$ and $D_{\hat{\thetavec}}>0$ such that
$\|\hat{\thetavec}(t)-\hat{\thetavec}^{\ast}\|\leq e^{-\lambda t}
\|\hat{\thetavec}(0)-\hat{\thetavec}^{\ast}\|
D_{\hat{\thetavec}}$. It follows from Theorem \ref{L_2
dpsi_theorem} that $\psi(\bfx(t),t)$ is bounded. In addition due
to Assumption \ref{boundedness_psi} we can conclude that $\bfx$ is
bounded as well. By the proposition assumptions function
$\alphavec(\bfx,t)$ is locally bounded with respect to $\bfx$ and
uniformly bounded in $t$. Therefore, there exists
$D_{\alphavec}>0$ such that
$|\alphavec(\bfx,t)^{T}(\hat\thetavec(t)-\hat{\thetavec}^{\ast})|\leq
D_{\alphavec}\|\hat{\thetavec}(t)-\hat{\thetavec}^{\ast}\|$.
Taking into account that
$z(\bfx,\thetavec,t)-z(\bfx,\hat\thetavec,t)=\mu(t)$ and
$|z(\bfx,\thetavec,t)-z(\bfx,\hat\thetavec,t)|\leq D
|\alphavec(\bfx,t)^{T}({\hat\thetavec}(t)-{\hat\thetavec}^{\ast})|$
 one can derive from (\ref{prop_exp_conv_dpsi}) the following
 estimate
\begin{eqnarray}
|\psi(t)|\leq |\psi(0)|e^{-K t} +
\|\hat{\thetavec}(0)-\hat{\thetavec}^{\ast}\|
D_{\hat{\thetavec}}D_{\alphavec} D
\int_0^{t}e^{-K(t-\tau)}e^{-\lambda\tau}d\tau  \leq |\psi(0)|e^{-K
t} + \frac{D_{\hat{\thetavec}}D_{\alphavec}D}{K-\lambda}
\|\hat{\thetavec}(0)-\hat{\thetavec}^{\ast}\| e^{-\lambda t}
\end{eqnarray}
{\it The proposition is proven.}

{\it Proof of Theorem \ref{theorem_realization1}.} The theorem
proof is quite straightforward and follows from explicit
differentiation of function $\hat{\thetavec}(\bfx,t)$ with respect
to time:
$\dot{\hat{\thetavec}}(\bfx,t)=\Gamma({\dot{\hat{\thetavec}}_{P}}+\dot{\hat\thetavec}_I)=\Gamma(\dpsi\alphavec(\bfx,t)+\psi\dot{\alphavec}(\bfx,t)-\dot{\Psi}(\bfx,t)+\dot{\hat\thetavec}_I)$.
Notice that
\begin{eqnarray}\label{t2_1}
&&\psi\dot{\alphavec}(\bfx,t)-\dot{\Psi}(\bfx,t)+\dot{\hat{\thetavec}}_I=\psi(\bfx,t)\frac{\pd
\alphavec(\bfx,t)}{\pd \bfx_1}\dot{\bfx}_1+\psi(\bfx,t)\frac{\pd
\alphavec(\bfx)}{\pd \bfx_2}\dot{\bfx}_2 + \psi(\bfx,t)\frac{\pd
\alphavec(\bfx,t)}{\pd t}-\nonumber\\
& & \frac{\pd \Psi(\bfx,t)}{\pd \bfx_1}\dot{\bfx}_1-\frac{\pd
\Psi(\bfx,t)}{\pd \bfx_2}\dot{\bfx}_2-\frac{\pd \Psi(\bfx,t)}{\pd
t}+\dot{\hat\thetavec}_I
\end{eqnarray}
According to Assumption \ref{Realizability_1},  $\frac{\pd
\Psi(\bfx,t)}{\pd \bfx_2}=\psi(\bfx,t)\frac{\pd
\alphavec(\bfx,t)}{\pd \bfx_2}$. Then taking into account
(\ref{t2_1}), we can obtain
\begin{eqnarray}\label{t2_2}
& &
\psi\dot{\alphavec}(\bfx,t)-\dot{\Psi}(\bfx,t)+\dot{\hat{\thetavec}}_I=\left(\psi(\bfx,t)\frac{\pd
\alphavec(\bfx,t)}{\pd \bfx_1}-\frac{\pd \Psi}{\pd \bfx_1
}\right)\dot{\bfx}_1+\psi(\bfx,t)\frac{\pd \alphavec(\bfx,t)}{\pd
t}-\frac{\Psi(\bfx,t)}{\pd t}+\dot{\hat{\thetavec}}_I
\end{eqnarray}
Notice that according to the proposed notations we can rewrite the
term $\left(\psi(\bfx,t)\frac{\pd \alphavec(\bfx,t)}{\pd
\bfx_1}-\frac{\pd \Psi}{\pd \bfx_1 }\right)\dot{\bfx}_1$ in the
following form: $\left(\psi(\bfx,t)L_{\bff_1}
\alphavec(\bfx,t)-L_{\bff_1} \Psi(\bfx,t)\right)+
\left(\psi(\bfx,t)L_{\bfg_1} \alphavec(\bfx,t)-L_{\bfg_1}
\Psi(\bfx,t)\right)u(\bfx,\hat{\thetavec},t)$. Hance it follows
from (\ref{fin_forms_ours_tr1}) and (\ref{t2_2}) that
$\psi\dot{\alphavec}(\bfx,t)-\dot{\Psi}(\bfx,t)+\dot{\hat{\thetavec}}_I=\varphi(\psi)\alphavec(\bfx,t)$.
Therefore
$\dot{\hat\thetavec}(\bfx,t)=\Gamma(\dpsi+\varphi(\psi))\alphavec(\bfx,t)$.
{\it The theorem is proven.}

{\it Proof of Theorem \ref{theorem_embedd}.} To prove the theorem,
first notice that control function (\ref{control_embedd}) provides
the following error model dynamics
\begin{eqnarray}\label{t5_1}
\dpsi&=&-\varphi(\psi)+z(\bfx,\thetavec,t)-z(\tilde\bfx,\hat\thetavec,t),
\
z(\tilde{\bfx},\hat{\thetavec},t)=L_{\nuvec(\tilde{\bfx},\hat{\thetavec})}\psi(\tilde{\bfx},t).
\end{eqnarray}
By adding and subtracting the function
$z(\tilde{\bfx},\thetavec,t)$ from the right-hand side of
(\ref{t5_1}) we get the following:
\[
\dpsi=-\varphi(\psi)+z(\bfx,\thetavec,t)-z(\tilde{\bfx},\thetavec,t)+z(\tilde{\bfx},\thetavec,t)-z(\tilde\bfx,\hat\thetavec,t),
\]
where the difference
$z(\bfx,\thetavec,t)-z(\tilde{\bfx},\thetavec,t)$ is bounded due
to Assumption \ref{assume_L_infty_embed}. Denote
$\varepsilon(t)=z(\bfx,\thetavec,t)-z(\tilde{\bfx},\thetavec,t)$,
then
\begin{eqnarray}\label{t_embed_error_model}
\dpsi=-\varphi(\psi)+\varepsilon(t)+z(\tilde{\bfx},\thetavec,t)-z(\tilde\bfx,\hat\thetavec,t),
\end{eqnarray}
where $\varepsilon \in L_\infty$. Denote
\begin{eqnarray}\label{t_embed_notations}
\tilde{\bff}&=&\bff_1(\bfx)\oplus\bff_2'(\bfx)\oplus\frac{\pd
\bfh_{\xi}}{\pd \xivec}\bff_{\xi}(\bfx,\xivec,t)\nonumber \\
\tilde{\bfg}&=&\bfg_1(\bfx)\oplus\bfg_2'(\bfx)\oplus\mathbf{0}_{\bfx_2''},
\ \
\mathbf{0}_{\bfx_2''}=(\underbrace{0,\dots,0}_{\dim{\bfx_2''}})^{T}\nonumber.
\end{eqnarray}
 Let us consider the
following adaptation algorithm:
\begin{eqnarray}\label{fin_forms_embedd}
\hat{\thetavec}(\bfx,\tilde{\bfx},t)&=&\Gamma(\hat{\thetavec}_{P}(\bfx,\tilde{\bfx},t)+\hat{\thetavec}_I(t)),
\ \Gamma>0; \ \
\hat{\thetavec}_P(\bfx,\tilde{\bfx},t)=\psi({\bfx},t)\alphavec(\tilde{\bfx},t)-\Psi(\tilde{\bfx},t);
\nonumber
\\
\dot{\hat{\thetavec}}_{I}&=&\varphi(\psi({\bfx},t))\alphavec(\tilde{\bfx},t)+\frac{\pd
\Psi(\tilde{\bfx},t)}{\pd t}-\psi(\bfx,t)\frac{\pd
\alphavec(\tilde{\bfx},t)}{\pd
t}-\left(\psi({\bfx},t)L_{\tilde\bff}\alphavec(\tilde{\bfx},t)-L_{\tilde\bff}\Psi(\tilde{\bfx},t)\right)-\nonumber \\
& &
\left(\psi({\bfx},t)L_{{\tilde\bfg}}\alphavec(\tilde{\bfx},t)-L_{{\tilde\bfg}}\Psi(\tilde{\bfx},t)\right)u(\bfx,\bfh_{\xi},\hat{\thetavec},t)-\lambda\hat{\thetavec}(\tilde\bfx,t),
\end{eqnarray}
where $\lambda>0$. Differentiation of function $\hat{\thetavec}_P$
with respect to time leads to:
\begin{eqnarray}\label{t_embed_dtheta_P}
\dot{\hat{\thetavec}}_P(\bfx,\tilde{\bfx},t)&=&\dpsi({\bfx},t)\alphavec(\tilde{\bfx},t)+\psi(\bfx,t)\dot{\alphavec}(\tilde\bfx,t)-\dot{\Psi}(\tilde{\bfx},t)=\dpsi({\bfx},t)\alphavec(\tilde{\bfx},t)+\psi(\bfx,t)\frac{\pd
\alphavec(\tilde\bfx,t)}{\pd t}-\frac{\pd \Psi(\tilde\bfx,t)}{\pd
t}\nonumber \\
& & + \left(\psi(\bfx,t)\frac{\pd \alphavec(\tilde\bfx,t)}{\pd
\bfx_1}-\frac{\pd \Psi(\tilde\bfx,t)}{\pd
\bfx_1}\right)\dot{\bfx}_1 + \left(\psi(\bfx,t)\frac{\pd
\alphavec(\tilde\bfx,t)}{\pd \bfx_2'}-\frac{\pd
\Psi(\tilde\bfx,t)}{\pd \bfx_2'}\right)\dot{\bfx}_2'+\nonumber \\
& & \left(\psi(\bfx,t)\frac{\pd \alphavec(\tilde\bfx,t)}{\pd
\bfh_\xi}-\frac{\pd \Psi(\tilde\bfx,t)}{\pd
\bfh_\xi}\right)\dot{\bfh}_\xi
\end{eqnarray}
Taking into account (\ref{t_embed_notations}) and
(\ref{partitioned_plant_embed}) we can rewrite
(\ref{t_embed_dtheta_P}) as follows:
\begin{eqnarray}\label{t_embed_dtheta_P1}
\dot{\hat{\thetavec}}_P(\bfx,\tilde{\bfx},t)&=&\dpsi({\bfx},t)\alphavec(\tilde{\bfx},t)-\frac{\pd
\Psi(\tilde{\bfx},t)}{\pd t}+\psi(\bfx,t)\frac{\pd
\alphavec(\tilde{\bfx},t)}{\pd
t}+\left(\psi({\bfx},t)L_{\tilde\bff}\alphavec(\tilde{\bfx},t)-L_{\tilde\bff}\Psi(\tilde{\bfx},t)\right)+ \\
& &
\left(\psi({\bfx},t)L_{{\tilde\bfg}}\alphavec(\tilde{\bfx},t)-L_{{\tilde\bfg}}\Psi(\tilde{\bfx},t)\right)u(\bfx,\bfh_{\xi},\hat{\thetavec},t)+
\left(\psi(\bfx,t)\frac{\pd \alphavec(\tilde\bfx,t)}{\pd
\bfx_2'}-\frac{\pd \Psi(\tilde\bfx,t)}{\pd
\bfx_2'}\right)\nuvec'(\bfx,\thetavec)\nonumber
\end{eqnarray}
Notice also that according to Assumption \ref{assume_embedd}:
\[
\frac{\pd \Psi(\tilde\bfx,t)}{\pd
\bfx_2'}=\psi(\tilde{\bfx},t)\frac{\pd
\alphavec(\tilde{\bfx},t)}{\pd \bfx_2}
\]
and
$(\psi(\bfx,t)-\psi(\tilde\bfx,t))L_{\nuvec'(\bfx,\thetavec)}\alphavec(\tilde{\bfx},t)\in
L_\infty$ due to Assumption \ref{assume_L_infty_embed}. Denoting
$(\psi(\bfx,t)-\psi(\tilde\bfx,t))L_{\nuvec'(\bfx,\thetavec)}\alphavec(\tilde{\bfx},t)=\delta(t)$
and using equalities (\ref{t_embed_dtheta_P1}) and
(\ref{fin_forms_embedd}) we can derive that
\begin{eqnarray}\label{t5_2}
\dot{\hat\thetavec}=\Gamma(\dot{\hat\thetavec}_{P}+\dot{\hat\thetavec}_I)=\Gamma((\dpsi+\varphi(\psi(\bfx,t)))\alphavec(\tilde{\bfx},t)+\delta(t)-\lambda\hat{\thetavec}),
\end{eqnarray}
where function $\delta(t)$ is bounded. Let us define the extended
state space vector $\bfq=\bfx\oplus\xivec$. Furthermore, we define
$z_{\bfq}(\bfq,\thetavec,t)=z(\tilde{\bfx},\thetavec,t)$,
$\alphavec_{\bfq}(\bfq,t)=\alphavec(\tilde{\bfx},t)$,
$\psi_\bfq(\bfq,t)=\psi(\bfx,t)$. Given the chosen notations,
Algorithm (\ref{t5_2}) can be written as follows:
%\begin{eqnarray}%\label{t5_3}
$\dot{\hat{\thetavec}}=\Gamma((\dpsi_{\bfq}+\varphi(\psi_{\bfq}))\alphavec_{\bfq}(\bfq,t)+\delta(t)-\lambda\hat\thetavec)$.
%\end{eqnarray}
Moreover instead of equation (\ref{t_embed_error_model}) we can
write
%\begin{eqnarray}%\label{t5_4}
$\dpsi=-\varphi(\psi)+z_{\bfq}(\bfq,\thetavec,t)-z_{\bfq}(\bfq,\hat{\thetavec},t)+\varepsilon(t)$.
%\end{eqnarray}

It is easy to see that Assumptions \ref{alpha} and
\ref{alpha_upper_bound} hold for the extended system. Assumption
\ref{boundedness_psi} is also satisfied with respect to the goal
function $\psi_\bfq(\bfq,t)$  due to hypothesis
(\ref{boundedness_xi}) in Assumption \ref{assume_embedd}. Indeed,
$\psi_{\bfq}(\bfq,t)=\psi(\bfx,t)\in L_\infty \Rightarrow \bfx\in
L_{\infty} \Rightarrow \xivec\in L_\infty \Rightarrow \bfq\in
L_{\infty}$. Therefore, according to Assumption
\ref{assume_L_infty_embed} and Lemma
\ref{lemma_dist_feedback_alg}, we can conclude that $\psi(\bfx,t)$
is bounded and furthermore trajectories $\bfx,\xivec$ are bounded
as well.
% If, in addition
%$\varepsilon(t),\delta(t)\rightarrow 0$ as $t\rightarrow\infty$
%then $\psi(\bfx,t)$ converges into
%\[
% |\varphi(\psi)|\leq
% \Delta_0+\sqrt\frac{\lambda}{{(D-D_1)}}\|{\hat\thetavec}^{\ast}\|,
%\]
%where $\Delta_0>0$ is arbitrary small constant.
Thus property P8) is proven. To prove property P9) it is
sufficient to notice that $\delta(t)=0$ either due to the equality
${\pd \alphavec(\tilde\bfx,t)}/{\pd \bfx_2'}\equiv 0$ or
$\psi(\bfx,t)=\psi(\tilde\bfx,t)$. Let ${\pd
\alphavec(\tilde\bfx,t)}/{\pd \bfx_2'}\equiv 0$, then P9) follows
explicitly from Assumption \ref{assume_L_2_embed} and Lemma
\ref{L_2_theorem_eps} applied to (\ref{t_embed_error_model}) with
algorithm
\begin{eqnarray}\label{fin_forms_embedd_1}
& &
\hat{\thetavec}(\bfx,\tilde{\bfx},t)=\Gamma(\hat{\thetavec}_{P}(\bfx,\tilde{\bfx},t)+\hat{\thetavec}_I(t)),
\ \Gamma>0; \ \
\hat{\thetavec}_P(\bfx,\tilde{\bfx},t)=\psi({\bfx},t)\alphavec(\tilde{\bfx},t);\nonumber
\\
& &
\dot{\hat{\thetavec}}_{I}=\varphi(\psi({\bfx},t))\alphavec(\tilde{\bfx},t)-\psi(\bfx,t){\pd
\alphavec(\tilde{\bfx},t)}/{\pd
t}-\psi({\bfx},t)L_{\tilde\bff}\alphavec(\tilde{\bfx},t)-
\left(\psi({\bfx},t)L_{{\tilde\bfg}}\alphavec(\tilde{\bfx},t)\right)u(\bfx,\bfh_{\xi},\hat{\thetavec},t).
\end{eqnarray}
which is in fact algorithm (\ref{fin_forms_embedd_2}) for
$\lambda=0$ and $\Psi(\tilde\bfx,t)\equiv0$. If
$\psi(\bfx,t)=\psi(\tilde\bfx,t)$, then according to Lemma
\ref{L_2_theorem_eps}, algorithm (\ref{fin_forms_embedd}) with
$\lambda=0$:
\begin{eqnarray}\label{fin_forms_embedd_2}
\hat{\thetavec}(\bfx,\tilde{\bfx},t)&=&\Gamma(\hat{\thetavec}_{P}(\bfx,\tilde{\bfx},t)+\hat{\thetavec}_I(t)),
\ \Gamma>0, \ \
\hat{\thetavec}_P(\bfx,\tilde{\bfx},t)=\psi({\bfx},t)\alphavec(\tilde{\bfx},t)-\Psi(\tilde{\bfx},t)\nonumber
\\
\dot{\hat{\thetavec}}_{I}&=&\varphi(\psi({\bfx},t))\alphavec(\tilde{\bfx},t)+{\pd
\Psi(\tilde{\bfx},t)}/{\pd t}-\psi(\bfx,t){\pd
\alphavec(\tilde{\bfx},t)}/{\pd
t}-\left(\psi({\bfx},t)L_{\tilde\bff}\alphavec(\tilde{\bfx},t)-L_{\tilde\bff}\Psi(\tilde{\bfx},t)\right)-\nonumber \\
& &
\left(\psi({\bfx},t)L_{{\tilde\bfg}}\alphavec(\tilde{\bfx},t)-L_{{\tilde\bfg}}\Psi(\tilde{\bfx},t)\right)u(\bfx,\bfh_{\xi},\hat{\thetavec},t).
\end{eqnarray}
ensures P9) as well. {\it The theorem is proven.}

{\it Proof of Theorem \ref{theorem_full_cascade}.} We will prove
the theorem by induction from order $1$ to $n$ for system
(\ref{cascade__full_2}). According to the theorem conditions
functions $f_1(x_1,\thetavec_1)$, $\alphavec_1(x_1)$ are smooth
and therefore according to Lemma \ref{L_2_theorem_eps} and Theorem
\ref{theorem_realization1}, there exists smooth function
$u(x_1,\hat{\thetavec}_1)$:
\begin{eqnarray}
u(x_1,\hat{\theta}_1)&=&-f_1(x_1,\hat{\thetavec}_1)-\varphi_1(\psi(x_1)),
\ \
 \hat{\thetavec}_1=\gamma_1(\hat{\thetavec}_{1,P}(x_1)+\hat{\thetavec}_{1,I}(t)),
\ \gamma_1>0 \nonumber\\
\hat{\thetavec}_{1,P}(x_1)&=&\psi(x_1)\alphavec_1(x_1)-\Psi(x_1),
\ \ \Psi(x_1)=\int_{x_1(0)}^{x_1(t)}\psi(x_1)\frac{\pd
\alphavec_{1}(x_1)}{\pd x_1} d x_1 \nonumber \\
\dot{\hat{\thetavec}}_{1,I}&=&\varphi_1(\psi(x_1))\alphavec_1(x_1)\nonumber
\end{eqnarray}
such that $\psi_1(x_1,t)\in L_{2}\cap L_{\infty}$, $\dpsi_1\in
L_2$ for the system of the following type: $
\dot{x}_1=f_1(x_1,\thetavec_1)+u(x_1,\hat{\thetavec}_1) +
\varepsilon_1(t)$, $\varepsilon_1(t)\in L_2$.  Hence, the basis of
induction is proven.

Let us assume that the theorem statements hold true for the
systems of order $i$, i.e. there exists such smooth function
$u_i(\bfx_i,\hat{\thetavec}_i,\xivec_i,\nuvec_i)$,
$\bfx_i,\xivec_i\in R^i$, $\bfx_i=(x_1,\dots,x_i)^{T}$,
$\xivec_i=(\xi_1,\dots,\xi_i)^T$ and the corresponding goal
functions $\psi_j(x_j,t)$, $j=1,\dots,i$ such that
$\psi_j(x_j,t)\in L_2\cap L_\infty$, $\dpsi_j\in L_2$ for system
(\ref{cascade__full_2}) of order $i$:
\begin{eqnarray}\label{cascade_full_4}
\dot{x}_j&=&f_j(x_1,\dots,x_j,\thetavec_j)+x_{j+1}, \ j\in\{1,\dots,i-1\}, \nonumber \\
\dot{x}_i&=&f_i(x_1,\dots,x_i,\thetavec_i)+u_i+\varepsilon_i(t), \
\varepsilon_i(t)\in L_2.
\end{eqnarray}
Therefore, in order to prove the theorem it is enough to show that
its statements hold for system (\ref{cascade__full_2}) of order
$i+1$ given that it holds for the systems like
(\ref{cascade_full_4}).

According to the inductive assumption function
$u_i(\bfx_i,\hat{\thetavec}_i,\xivec_i,\nuvec_i)$ is smooth. Then
by Hadamar's lemma there exists such
$F(\bfx_i,\bfx_i',\hat{\thetavec}_i,\xivec_i,\nuvec_i)$ that
$u_i(\bfx_i,\hat{\thetavec}_i,\xivec_i,\nuvec_i)-u_i(\bfx_i',\hat{\thetavec}_i,\xivec_i,\nuvec_i)=F(\bfx_i,\bfx_i',\hat{\thetavec}_i,\xivec_i,\nuvec_i)(\bfx_i-\bfx_i')$.
Let us denote
$\bar{F}^2_{i+1}(\bfx_i,\bfx_i',\hat{\thetavec}_i,\xivec_i,\nuvec_i)=\|F(\bfx_i,\bfx_i',\hat{\thetavec}_i,\xivec_i,\nuvec_i)\|^2$.
Furthermore, due to the theorem conditions functions
$f_j(\bfx_j,\thetavec_j)$, $j=1,\dots,i+1$ satisfy the following
additional assumptions:
\[
(f_j(\bfx_j,\thetavec_j)-f_j(\bfx_j',\thetavec_j))^2\leq
\|{\bfx}_j-{\bfx}_j'\|^2 \bar{D}_j^2(\bfx_j,\bfx_j')  \ \ \forall
\ \thetavec_j\in \Omega_{\theta}
\]
Therefore, it follows from Lemma \ref{lemma_embed_full_cascade}
that there exists a system of differential equations
\begin{eqnarray}
\dot{\xivec}_{i+1}&=&\bff_{\xi_{i+1}}(\xivec_{i+1},\bfx,\bfz,\nuvec_{i+1}),
\xivec\in R^{i}\nonumber \\
\dot{\nuvec}_{i+1}&=&\bff_{\nu_{i+1}}(\xivec_{i+1},\bfx,\bfz), \ \
\bfz= \hat{\thetavec}_i\oplus\xivec_i\oplus\nuvec_i\nonumber
\end{eqnarray}
such that
\begin{eqnarray}
& &
u_i(\bfx_i,\hat{\thetavec}_i(\bfx_i,\xivec_i,\hat{\thetavec}_{I,i}),\xivec_i,\nuvec_i)-u_i(\xivec_{i+1},\hat{\thetavec}_i(\xivec_{i+1},\xivec_i,\hat{\thetavec}_{I,i}),\xivec_i,\nuvec_i)\in
L_2\nonumber \\
& & f_{i+1}(\bfx_{i+1},\thetavec_{i+1})-
f_{i+1}(\xivec_{i+1}\oplus x_{i+1},\thetavec_{i+1})\in L_2,
\nonumber
% & & \tilde{\xivec}_{i+1}=(\xivec_{i+1}\oplus x_{i+1})^{T}
\end{eqnarray}
Let us introduce new goal function
$\psi_{i+1}(x_{i+1},\xivec_{i+1},\xivec_{i},\nuvec_{i},\hat{\thetavec}_{I,i})=x_{i+1}-u_i(\xivec_{i+1},\hat{\thetavec}_i(\xivec_{i+1},\xivec_i,\hat{\thetavec}_{I,i}),\xivec_i,\nuvec_i)$
and consider its time-derivative $\dpsi_{i+1}$:
\begin{eqnarray}\label{cascade_full_5}
\dpsi_{i+1}&=&f_{i+1}(\bfx_{i+1},\thetavec_{i+1}) + u_{i+1} -
L_{\bff_{\xi_i}}u_i(\xivec_{i+1},\hat{\thetavec}_i(\xivec_{i+1},\xivec_i,\hat{\thetavec}_{I,i}),\xivec_i,\nuvec_i)-\nonumber
\\
& &
L_{\bff_{\nu_i}}u_i(\xivec_{i+1},\hat{\thetavec}_i(\xivec_{i+1},\xivec_i,\hat{\thetavec}_{I,i}),\xivec_i,\nuvec_i)-L_{\bff_{\xi_{i+1}}}u_i(\xivec_{i+1},\hat{\thetavec}_i(\xivec_{i+1},\xivec_i,\hat{\thetavec}_{I,i}),\xivec_i,\nuvec_i)
-\nonumber \\
& &
L_{\bff_{\hat{\theta}_i}}u_i(\xivec_{i+1},\hat{\thetavec}_i(\xivec_{i+1},\xivec_i,\hat{\thetavec}_{I,i}),\xivec_i,\nuvec_i).
\end{eqnarray}
Denote $\varepsilon_{i+1}(t)=f_{i+1}(\bfx_{i+1},\thetavec_{i+1})-
f_{i+1}(\xivec_{i+1}\oplus x_{i+1},\thetavec_{i+1})$ and rewrite
(\ref{cascade_full_5}) in the following way:
\begin{eqnarray}\label{cascade_full_6}
\dpsi_{i+1}&=&f_{i+1}(\xivec_{i+1}\oplus x_{i+1},\thetavec_{i+1})
+ u_{i+1}  -
L_{\bff_{\xi_i}}u_i(\xivec_{i+1},\hat{\thetavec}_i(\xivec_{i+1},\xivec_i,\hat{\thetavec}_{I,i}),\xivec_i,\nuvec_i)-\nonumber
\\
& &
L_{\bff_{\nu_i}}u_i(\xivec_{i+1},\hat{\thetavec}_i(\xivec_{i+1},\xivec_i,\hat{\thetavec}_{I,i}),\xivec_i,\nuvec_i)-L_{\bff_{\xi_{i+1}}}u_i(\xivec_{i+1},\hat{\thetavec}_i(\xivec_{i+1},\xivec_i,\hat{\thetavec}_{I,i}),\xivec_i,\nuvec_i)
-\nonumber \\
& &
L_{\bff_{\hat{\theta}_i}}u_i(\xivec_{i+1},\hat{\thetavec}_i(\xivec_{i+1},\xivec_i,\hat{\thetavec}_{I,i}),\xivec_i,\nuvec_i+
\varepsilon_{i+1}(t).
\end{eqnarray}
Let us select input $u_{i+1}$ as follows
\begin{eqnarray}\label{cascade_full_7}
u_{i+1}&=&-\varphi_{i+1}(\psi_{i+1}(x_{i+1},\xivec_{i+1},\xivec_{i},\nuvec_{i},\hat{\thetavec}_{I,i}))+
L_{\bff_{\xi_i}}u_i(\xivec_{i+1},\hat{\thetavec}_i(\xivec_{i+1},\xivec_i,\hat{\thetavec}_{I,i}),\xivec_i,\nuvec_i)+\nonumber
\\
& &
L_{\bff_{\nu_i}}u_i(\xivec_{i+1},\hat{\thetavec}_i(\xivec_{i+1},\xivec_i,\hat{\thetavec}_{I,i}),\xivec_i,\nuvec_i)+L_{\bff_{\xi_{i+1}}}u_i(\xivec_{i+1},\hat{\thetavec}_i(\xivec_{i+1},\xivec_i,\hat{\thetavec}_{I,i}),\xivec_i,\nuvec_i)
+\nonumber \\
& &
L_{\bff_{\hat{\theta}_i}}u_i(\xivec_{i+1},\hat{\thetavec}_i(\xivec_{i+1},\xivec_i,\hat{\thetavec}_{I,i}),\xivec_i,\nuvec_i)-f_{i+1}(\xivec_{i+1}\oplus
x_{i+1},\hat{\thetavec}_{i+1})
\end{eqnarray}
Denoting
$\psi_{i+1}(x_{i+1},\xivec_{i+1},\xivec_{i},\nuvec_{i},\hat{\thetavec}_{I,i})=\psi_{i+1}(x_{i+1},t)$
(as $\xivec_{i+1},\xivec_{i},\nuvec_{i}, \hat{\thetavec}_{I,i}$
are functions of time $t$) and substituting (\ref{cascade_full_7})
and (\ref{cascade_full_6}) into (\ref{cascade_full_5}) we can
write the following expression for $\dpsi_{i+1}$
\begin{eqnarray}\label{cascade_full_8}
\dpsi_{i+1}&=&-\varphi_{i+1}(\psi_{i+1}(x_{i+1},t))+f_{i+1}(\xivec_{i+1}\oplus
x_{i+1},\thetavec_{i+1})- f_{i+1}(\xivec_{i+1}\oplus
x_{i+1},\hat{\thetavec}_{i+1})+\varepsilon_{i+1}(t)
\end{eqnarray}
It follows from the theorem conditions that there exits such
function $\alphavec_{i+1}(\xivec_{i+1}\oplus x_{i+1})$ that
Assumptions \ref{alpha} and \ref{alpha_upper_bound} are satisfied
for the function $f_{i+1}(\xivec_{i+1}\oplus
x_{i+1},{\thetavec}_{i+1})$. Consider the following adaptation
algorithm:
\begin{eqnarray}\label{cascade_full_9}
\dot{\hat{\thetavec}}_{i+1}=\gamma_{i+1}(\dpsi_{i+1}+\varphi_{i+1}(\psi_{i+1}(x_{i+1},t)))\alphavec_{i+1}(\xivec_{i+1}\oplus
x_{i+1}), \ \gamma_{i+1}>0
\end{eqnarray}
Realization of algorithms (\ref{cascade_full_9})  is guaranteed by
Theorem \ref{theorem_realization1} and can be given as follows:
\begin{eqnarray}\label{cascade_full_10}
\hat{\thetavec}_{i+1}(\xivec_{i+1}\oplus
x_{i+1},t)&=&\gamma_{i+1}(\hat{\thetavec}_{i+1,P}(\xivec_{i+1}\oplus
x_{i+1},t)+\hat{\thetavec}_{i+1,I}(t)), \ \gamma_{i+1}>0 \nonumber
\\
\hat{\thetavec}_{i+1,P}(\xivec_{i+1}\oplus x_{i+1},t)&=&
\psi_{i+1}(x_{i+1},t)\alphavec_{i+1}(\xivec_{i+1}\oplus x_{i+1})-\Psi_{i+1}(\xivec_{i+1}\oplus x_{i+1},t)\nonumber\\
\dot{\hat{\thetavec}}_{i+1,I}&=&\varphi_{i+1}(\psi(x_{i+1},t))\alphavec_{i+1}(\xivec_{i+1}\oplus
x_{i+1}) - L_{\bff_{\xi}}\alphavec_{i+1}(\xivec_{i+1}\oplus
x_{i+1})+\nonumber
\\
& & L_{\bff_{\xi}}\Psi_{i+1}(\xivec_{i+1}\oplus
x_{i+1},t)+\frac{\pd \Psi_{i+1}(\xivec_{i+1}\oplus x_{i+1},t)}{\pd
t}
\nonumber \\
\Psi_{i+1}(\xivec_{i+1}\oplus
x_{i+1},t)&=&\int_{x_{i+1}(0)}^{x_{i+1}(t)}\psi_{i+1}(x_{i+1},t)\frac{\pd
\alphavec_{i+1}(\xivec_{i+1}\oplus x_{i+1})}{\pd x_{i+1}}d
x_{i+1}.
\end{eqnarray}
It follows from Lemma \ref{L_2_theorem_eps} that for the error
model (\ref{cascade_full_8}) with adaptation algorithm
(\ref{cascade_full_9}) and its realization (\ref{cascade_full_10})
the following statements hold true:
$\hat{\thetavec}_{i+1}(\xivec_{i+1}\oplus x_{i+1},t)\in L_\infty$,
$\dpsi_{i+1}\in L_2$, $\varphi_{i+1}(\psi_{i+1}(x_{i+1},t))\in
L_2\cap L_{\infty}$. Given that $\varepsilon_{i+1}(t)\in L_2$ we
can conclude that
\[
u_{i+1}(\bfx_i,\hat{\thetavec}_{i+1},\xivec_{i+1},\xivec_{i},\nuvec_{i+1},\nuvec_{i},\hat{\thetavec}_{I,i})-u_{i+1}(\bfx_i,\thetavec_{i+1},\xivec_{i+1},\xivec_{i},\nuvec_{i+1},\nuvec_{i},\hat{\thetavec}_{I,i})\in
L_2
\]
Let us denote
$u_{i+1}(\bfx_i,\hat{\thetavec}_{i+1},\xivec_{i+1},\tilde{\nuvec}_{i+1})=u_{i+1}(\bfx_i,\hat{\thetavec}_{i+1},\xivec_{i+1},\xivec_{i},\nuvec_{i+1},\nuvec_{i},\hat{\thetavec}_{I,i})$,
$\tilde{\nuvec}_{i+1}=\xivec_{i}\oplus\nuvec_{i+1}\oplus\nuvec_{i}\oplus\hat{\thetavec}_{I,i}$.
According to the introduced notations it is easy to see that
statement 2) of the theorem holds. In addition to this notice that
the choice of appropriate function $\varphi_{i+1}(\cdot)$ in
(\ref{cascade_full_7}) is up to the designer. Therefore choosing
$\varphi_{i+1}(\cdot): \ |\varphi_{i+1}(\cdot)|\geq
k_{k+1}|\psi_{i+1}|, \ k_{i+1}>0$ we can guarantee that
$\psi_{i+1}(x_{i+1},t)\in L_2\cap L_\infty$.

The last, however, according to the inductive hypothesis implies
that $\psi_k(x_k,t)\in L_2\cap L_\infty$ and $\dpsi_k\in L_2$ for
any $k=1,\dots,i$, $\psi\in L_2\cap L_\infty$ and $\dpsi\in L_2$.
Hence statement 1) of the theorem is proven as well.

Let us prove statement 3). According to the inductive hypothesis
$\bfx_i,\xivec_i,\nuvec_i,\hat{\thetavec}_{i}$ are bounded.
Furthermore, $\hat{\thetavec}_{I,i}(t)$  is bounded as
$\hat{\thetavec}_{P,i}(\bfx_i,\xivec_i)$ is smooth function and
$\hat{\thetavec}_i=\gamma_i(\hat{\thetavec}_{P,i}(\bfx_i,\xivec_i)+\hat{\thetavec}_{I,t})$.
Then taking into account Lemma \ref{lemma_embed_full_cascade} we
can conclude that $\xivec_{i+1},\nuvec_{i+1}$ are bounded. Hence
$\tilde{\nuvec}_{i+1}$ is bounded. Let us show that $x_{i+1}$ is
bounded as well. First notice that the difference
$\varepsilon_{i}(t)=u_i(\bfx_i,\hat{\thetavec}_i(\bfx_i,\xivec_i,\hat{\thetavec}_{I,i}),\xivec_i,\nuvec_i)-u_i(\xivec_{i+1},\hat{\thetavec}_i(\xivec_{i+1},\xivec_i,\hat{\thetavec}_{I,i}),\xivec_i,\nuvec_i)$
is bounded as $u_i$ is smooth and its arguments are bounded. On
the other hand we have just shown that
$\psi_{i+1}=x_{i+1}-u_i(\xivec_{i+1},\hat{\thetavec}_i(\xivec_{i+1},\xivec_i,\hat{\thetavec}_{I,i}),\xivec_i,\nuvec_i)$
is bounded. Therefore, $x_{i+1}$ is bounded. Hence statement 3) is
proven.

Derivatives $\dot{\psi}_j$, $j=1,\dots,i$ are bounded as
$\varepsilon_i(t)$ is bounded (according to the inductive
hypothesis the theorem holds for any $j=1\dots,i$). If, however,
$\varepsilon(t)$ is bounded then $\dot{\psi}_{i+1}$ is bounded as
well as $u_{i+1}(\cdot)$, $f_{i+1}(\cdot)$ are smooth and
$\bfx_{i+1},\xivec_{i+1},\nuvec_{i+1},\hat{\thetavec}_{i+1}$ are
bounded. Therefore, it follows from Lemma \ref{L_2_theorem_eps}
that $\psi_{i+1}\rightarrow 0$ as $t\rightarrow\infty$. Thus
statement 4) is proven.
%Functions $f_{i+1}(\cdot)$, $u_{i+1}(\cdot)$ and $u_i(\cdot)$ are
%smooth. Let $\dot{\varepsilon}(t)\in L_\infty$. Therefore
%derivative $\ddot{\psi}_{i+1}$ is bounded as $\bfx_{i+1}$,
%$\xivec_{i+1}$, $\nuvec_{i+1}$, $\xivec_i$, $\nuvec_i$ are bounded
%and $\bff_{\xi_{i+1}}$, $\bff_{xi_i}$, $\bff_{\nu_{i+1}}$,
%$\bff_{nu_i}$ are locally bounded. On the other hand we have that
%$\dpsi_{i+1}\in L_2$. Then it follows from Barbalatt' lemma that
%$\dpsi_{i+1}\rightarrow 0$ as $t\rightarrow \infty$. Notice also
%that $\dot{\varepsilon}_i(t)\in L_\infty$ and therefore, according
%to the inductive hypothesis $u_i(\cdot)\rightarrow 0$ as
%$t\rightarrow \infty$. Hence $x_{i+1}\rightarrow 0$ and moreover
{\it The theorem is proven.}

\end{document}